# Blind Hyperspectral and Multispectral Images Fusion: A Unified Tensor Fusion Framework from Coupled Inverse Problem Perspective *

Ying Gao†, Michael K. Ng‡, and Chunfeng Cui§

**Abstract.** Hyperspectral and multispectral images fusion aims at integrating a low-resolution hyperspectral image (LR-HSI) and a high-resolution multispectral image (HR-MSI) to construct a high-resolution hyperspectral image (HR-HSI). It is generally assumed that spatial blurring operator and spectral response operator are prior-known. However, such an assumption is extremely restrictive in practice. To overcome this limitation, this paper formulates blind fusion as a coupled inverse problem, integrating blind deconvolution in the spatial domain with blind unmixing in the spectral domain. From this novel perspective, we propose a unified tensor fusion framework capable of flexible self-adjustment and real-time fusion without pre-training. We further introduce an optimization model for the joint estimation of the target HR-HSI, the spatial point spread function, and the spectral response function. To solve this model, we devise a partially linearized alternating direction method of multipliers (ADMM) algorithm with Moreau envelope smoothing, accompanied by the rigorous convergence analysis. An initialization estimator tailored to the specific characteristics of the fusion problem is proposed. Numerical comparisons with state-of-the-art methods on both synthetic and real-world datasets demonstrate the compelling performance of the proposed method.

## 1. Introduction.

The trade-off between spatial resolution and spectral resolution represents a fundamental constraint [46, 51] in spectral sensing, particularly in the field of remote sensing. Due to the limited total energy of sunlight reflected from the scene and the fact that photons emitted by the sun are distributed across multiple spectral bands, the spatial resolution must be reduced to ensure that the number of photons in each band remains above a minimum threshold. This inherent trade-off leads to hyperspectral images typically exhibiting low spatial resolution, whilst multispectral sensors capture multispectral images with higher spatial resolution. Therefore, the fusion of HSI and MSI has gained increasing research interest to reconstruct the remote sensing images with high spectral and spatial resolutions (see, e.g., [14, 15, 19, 20]). In the context of hyperspectral image fusion, the observed HSI and MSI are respectively modeled as spatially and spectrally degraded versions of the target HR-HSI. Accordingly, the fusion is considered non-blind when the degradation process is known in advance; otherwise, it is termed a blind fusion.

---

*

**Funding:** This work was supported by the GDSTC: Guangdong and Hong Kong Universities "1+1+1" Joint Research Collaboration Scheme UICR0800008-24, and the National Natural Science Foundation of China (Nos. 12471282, 12131004)

†School of Mathematical Sciences, Beihang University, Beijing, 100191, China (yinggao@buaa.edu.cn).

‡Department of Mathematics, Hong Kong Baptist University, Hong Kong (michael-ng@hkbu.edu.hk).

§Corresponding author. School of Mathematical Sciences, Beihang University, Beijing, 100191, China (chunfengcui@buaa.edu.cn).

Over the past decade, we have witnessed a flurry of research activities on non-blind fusion. Existing model-based methods can be classified into two categories: matrix-based methods and tensor-based methods. The first line of research aims to enhance the spatial resolution of HSI by statistical priors [21, 33, 45] or subspace techniques [18, 25]. A notable matrix-based method incorporating Bayesian priors is the work [21], which introduced a maximum a posteriori (MAP) estimator for the HSI-MSI fusion problem. Later, within the Bayesian framework, Wei et al. [45] proposed a fast multispectral image fusion method by solving a Sylvester equation, which was further generalized to incorporate various priors such as Gaussian, sparse representation, and total variation. Furthermore, Nezhad et al. [33] introduced a fusion framework combining spectral unmixing and sparse coding to enhance the spatial resolution of HSI while preserving its spectral fidelity. The method first extracts spectral endmembers from the HSI through unmixing techniques and then refines their spatial distribution using a co-registered MSI. Wei et al. [44] developed an unmixing-based fusion scheme that incorporated the non-negativity and sum-to-one constraints as prior information to address the fusion and unmixing problems. Matrix-based methods with subspace techniques typically begin by unfolding the LR-HSI tensor into matrix form, which is then factorized into basis and coefficient matrices. A notable advancement was made by Kawakami et al. [25], who introduced a two-stage approach that first estimates the spectral basis using a sparse prior and subsequently computes the coefficient matrix via sparse coding. Comparative analysis demonstrates that the joint optimization framework, as proposed in [27, 52], which alternately updates the spectral basis and coefficient matrix, outperforms the sequential two-stage framework in terms of reconstruction accuracy. Afterwards, in [18], a non-negative dictionary-learning algorithm was devised to estimate the spectral basis and the coefficient matrix by a structured sparse coding method.

The tensor-based methods for non-blind fusion primarily decompose the target HR-HSI using various low-rank decomposition approaches, then construct the optimization model based on the known spatial or spectral degradation. Notably, Kanatsoulis et al. [24] adapted the CP decomposition to characterize the low-rankness of HR-HSI and proposed the coupled CP decomposition model for HSI-MSI fusion. Li et al. [30] advanced the field by employing Tucker decomposition and established a coupled sparse tensor factorization model. An interpretable block-term tensor model was proposed by Ding et al. [17] for hyperspectral super-resolution, where the latent factors of decomposition correspond to physically meaningful properties of spectral image constituents. In [36], Prévost et al. adopted LL1 decomposition and blind unmixing to address the hyperspectral super-resolution problem with spectral variability. In addition, the tensor train decomposition, tensor ring decomposition, and tensor triple decomposition have also been employed for non-blind tasks (see, e.g., [9, 13, 47] and references therein). Dian et al. [16] recently proposed the generalized tensor nuclear norm (GTNN) to overcome the limitations of the traditional tensor nuclear norm (TNN) in non-blind hyperspectral image fusion. Unlike TNN, which applies the Fourier transform solely along the third mode of a tensor, GTNN extends this by performing Fourier transforms along all three modes, thereby capturing more comprehensive correlations and reducing sensitivity to tensor mode permutations.

In contrast, significantly less attention has been devoted to blind fusion. The limited number of existing blind methods can be classified into three categories: hypersharpening-based methods, model-based methods, and learning-based methods. Hypersharpening-based methods (see, e.g., [4, 35, 38]) first employs regression-based methods to determine a weight vector linking the degraded MSI bands to each HSI band. It then applies the weighted combination of all high spatial-resolution MSI bands to enhance the spatial resolution of each HSI band. Despite their advantages in simplicity and speed, hypersharpening-based methods suffer from an inherent limitation that restricts solution quality and lacks a solid theoretical foundation. Model-based methods (see, e.g., [39, 50]) focus on characterizing the relationship between the observed HSI-MSI pair and the target HR-HSI to establish the degradation model, and introduce additional regularization to constrain the solution space and enable accurate reconstruction. A prominent model-based method is Hysure [39], which establishes a two-stage framework: it initially estimates spatial and spectral degradation operators from input images, followed by target HR-HSI reconstruction through optimization. Error accumulation in the two-stage framework potentially restricts the high-precision reconstruction. Recently, Yang et al. [50] proposed an innovative subspace-based triple decomposition for the blind fusion model, in which the target HR-HSI is decomposed into a coefficient tensor and a dictionary, and a novel triple decomposition is employed on the coefficient tensor. Unsupervised and self-supervised learning methods (see, e.g., [12, 29, 34]) have attracted much attention for the blind fusion task, as they operate without external training data. By integrating the estimation of degradation operators into the network itself, these methods learn directly from the input within an end-to-end framework. While achieving super-resolution quality, deep neural networks confront inherent challenges arising from their high-dimensional parameter spaces, including prohibitive computational demands and unstable performance when processing cross-sensor data in real-world applications.

**1.1. Contributions.** This research focuses on developing a novel model-based method to address the blind fusion problem, which is capable of offering high-quality solutions with theoretical guarantees and computational efficiency. We formulate blind fusion as a coupled inverse problem that integrates blind deconvolution in the spatial domain with blind unmixing in the spectral domain. Building on this perspective, we propose a unified tensor fusion framework guided by the sensor information, and further introduce an optimization model. Fundamentally, our framework operates through the integration of sensor information, physical degradation model, and customized algorithm, all of which are essential for achieving high-fidelity reconstruction of HR-HSI. The main contributions of this paper are summarized as three aspects.

(i) A unified tensor fusion framework integrating blind deconvolution and blind unmixing. Our framework jointly estimates the target HR-HSI, the point spread function (PSF) for deconvolution, and the spectral response function (SRF) for unmixing. This effectively circumvents the error propagation inherent in the classical two-stage blind model-based method [39].

(ii) A structural degradation model with physical interpretability. The intrinsic sensor information embedded in the observed HSI-MSI data provides the physically-grounded downsampling spatial operator and the spectral response band ranges for our framework. On the other hand, we further introduce several principled regularizations or constraints to restrict the solution space: (a) the transformed tubal nuclear norm

Table 1.1

*Several related works of convergence results about ADMM-type algorithms for solving problem (1.1). The convergence is obtained under different conditions with respect to the subblocks.*

| Algorithm | $\{\psi_i\}_{i=1}^{n-1}$ | $\psi_n$ | $\varphi_1$ | $\{\varphi_j\}_{j=2}^{m}$ | Convergence (rate) |
|---|---|---|---|---|---|
| ADMM [22] | $L$-smooth | $L$-smooth | Nonsmooth | - | subseq./$\mathcal{O}(\epsilon^{-1})$ |
| ADMM [43] | Nonsmooth | Nonsmooth | $L$-smooth | - | subseq./seq. |
| proximal-ADMM [7] | Nonsmooth | - | Nonsmooth | - | subseq./seq. |
| gradient-based ADMM [23] | Nonsmooth | - | - | - | $\mathcal{O}(\epsilon^{-2})$ |
| partially linearized ADMM (Ours) | Nonsmooth | Nonsmooth | Nonsmooth | Nonsmooth | subseq.*/seq.*/$\mathcal{O}(\epsilon^{-4})$ |

**Note:** $n$ and $m$ are the number of blocks of $\boldsymbol{x}$ and $\boldsymbol{y}$, respectively. For simplicity, "subseq." and "seq." denote subsequential and sequential convergence to a stationary point of the original problem, respectively; $\mathcal{O}(\cdot)$ denotes the iteration complexity with respect to the original problem; "subseq.*" and "seq.*" denote subsequential and sequential convergence to an $\epsilon$ stationary point of the original problem, respectively, under certain conditions; A dash (-) indicates the information is not required.

(TTNN) with non-negativity constraints for the target HR-HSI reconstruction; and (b) unit simplex constraints applied to both the spatial blur kernel for PSF and the spectral response vector for SRF.

(iii) A customized algorithm with provable convergence guarantees. The resulting minimization is a nonconvex nonsmooth problem with linear equation constraint. However, its inherent multiblock and nonseparable structure poses a challenge for existing ADMM-type algorithms, rendering them inefficient. Therefore, we devise a partially linearized ADMM algorithm with Moreau envelope smoothing to address the above difficulties. In addition, we provide a rigorous theoretical convergence of the proposed algorithm. Specifically, we demonstrate that the sequence generated by our proposed partially linearized ADMM algorithm converges to a stationary point of the Moreau smoothed problem. Furthermore, we establish iteration complexities of $\mathcal{O}(\epsilon^{-2})$ and $\mathcal{O}(\epsilon^{-4})$ to reach an $\epsilon$-stationary point of the Moreau smoothed problem and the original problem, respectively.

**1.2. Related work.** To clarify our algorithmic contribution, we consider the general minimization as follows,

$$\min_{\boldsymbol{x},\boldsymbol{y}} \Psi(\boldsymbol{x},\boldsymbol{y}) := f(\boldsymbol{x}_1,\ldots,\boldsymbol{x}_n) + g(\boldsymbol{y}_1,\ldots,\boldsymbol{y}_m) + \sum_{i=1}^{n}\psi_i(\boldsymbol{x}_i) + \sum_{j=1}^{m}\varphi_j(\boldsymbol{y}_j) \quad \text{s.t. } \boldsymbol{A}\boldsymbol{x} + \boldsymbol{B}\boldsymbol{y} = 0, \tag{1.1}$$

where $f$ and $g$ are Lipschitz continuously differentiable with respect to each block and typically represent loss functions; the $\psi_i$ and $\varphi_j$ denote the regularization terms imposed on certain blocks.

Table 1.1 summarizes several ADMM-type algorithms designed for structural nonconvex and nonsmooth minimization, including those that are multiblock and nonseparable. The main advantage of this work is that the convergence of our algorithm is guaranteed under weaker assumptions than those required by existing methods. More specifically, a critical condition for establishing convergence is the boundedness of the Lagrange multipliers. The boundedness, in turn, requires the objective function $\Psi(\mathbf{x}, \mathbf{y})$ to be Lipschitz differentiable with respect to at least one subblock (see, e.g., [7, 23, 43]). However, the resulting optimization model is equivalent to imposing nonsmooth regularizations on each block, which consequently fails to satisfy the required condition for convergence. Motivated by the Moreau envelope smoothing technique in [28, 53], our proposed algorithm, supported by its convergence analysis, is able to circumvent the aforementioned limitation, and is directly applicable to the resulting model. In terms of convergence results, we establish the subsequential and sequential convergence to a stationary point of the Moreau-smoothed problem and derive an iteration complexity bound for the original problem.

**1.3. Organization.** The rest of this paper is organized as follows. Section 2 outlines the preliminaries necessary for the following discussions. In Section 3, we introduce the background and analyze the structure of the blind fusion problem. Furthermore, we present the unified tensor framework followed by our proposed model. The algorithm is presented in Section 4, followed by its convergence in Section 5. The experimental results are given and discussed in Section 6. Finally, we conclude some remarks in Section 7.

## 2. Preliminaries.

In this paper, scalars are denoted by small letters $(a, b, \ldots)$, vectors are written in bold small letters $(\boldsymbol{a}, \boldsymbol{b}, \ldots)$, matrices correspond to bold capital letters $(\boldsymbol{A}, \boldsymbol{B}, \ldots)$, and tensors are written as bold calligraphic letters $(\boldsymbol{\mathcal{X}}, \boldsymbol{\mathcal{Y}}, \ldots)$.

Given an $n$-dimensional vector $\boldsymbol{x} \in \mathbb{R}^n$, let $\|\boldsymbol{x}\|_p := (\sum_{i=1}^n |x_i|^p)^{1/p}$ denote the $\ell^p$-norm of $\boldsymbol{x}$. Particularly, $\|\boldsymbol{x}\| := \|\boldsymbol{x}\|_2$ for brevity. The distance from $\boldsymbol{x} \in \mathbb{R}^n$ to a set $\Omega \subset \mathbb{R}^n$ is defined by $\operatorname{dist}(\boldsymbol{x}, \Omega) := \inf_{\boldsymbol{y} \in \Omega} \|\boldsymbol{x} - \boldsymbol{y}\|$, and the projection onto $\Omega$ is $P_\Omega(\boldsymbol{x}) := \operatorname{argmin}_{y \in \Omega} \|\boldsymbol{x} - \boldsymbol{y}\|$. For the matrix $\boldsymbol{X} \in \mathbb{R}^{m \times n}$, $\boldsymbol{X} \geq 0$ denotes that all elements of matrix are non-negative. $\operatorname{Diag}(\boldsymbol{x})$ denotes the diagonal matrix formed by the vector $\boldsymbol{x}$, while the notation $\operatorname{diag}(\boldsymbol{X})$ refers to the vector obtained by extracting the diagonal elements of the matrix $\boldsymbol{X}$. The spectral norm of a matrix $\boldsymbol{X}$ is defined as $\|\boldsymbol{X}\|_2 = \max_i \sigma_i(\boldsymbol{X})$, where $\sigma_i(\boldsymbol{X})$ is the $i$-th largest singular value of $\boldsymbol{X}$. The nuclear norm of $\boldsymbol{X}$ is defined as $\|\boldsymbol{X}\|_* = \sum_i \sigma_i(\boldsymbol{X})$. For the tensor $\boldsymbol{\mathcal{X}} \in \mathbb{R}^{I_1 \times \cdots \times I_n}$, $\boldsymbol{\mathcal{X}} \geq 0$ denotes that all elements of tensor are non-negative. The mode-$n$ matricization of $\boldsymbol{\mathcal{X}}$ is a matrix denoted by $\boldsymbol{X}_{(n)}$ whose columns are mode-$n$ fibers in lexicographical order. Especially, we denote $\boldsymbol{\mathcal{Y}}(i, j, :)$ as the mode-3 fiber of $\boldsymbol{\mathcal{Y}}$ for $i = 1 \ldots, I_1$ and $j = 1, \ldots, I_2$. The mode-$n$ product of $\boldsymbol{\mathcal{X}}$ with $\boldsymbol{A}$ denoted as $\boldsymbol{\mathcal{Y}} = \boldsymbol{\mathcal{X}} \times_n \boldsymbol{A}$, is defined element-wise as

$$(\boldsymbol{\mathcal{Y}})_{i_1 \cdots i_{n-1} j i_{n+1} \cdots i_N} = \sum_{i_n=1}^{I_n} x_{i_1 i_2 \cdots i_N} a_{j i_n}.$$

The mode-$n$ product can also be represented as matrix multiplication, $\boldsymbol{Y}_{(n)} = \boldsymbol{A}\boldsymbol{X}_{(n)}$. The inner product of $\boldsymbol{\mathcal{X}}$, $\boldsymbol{\mathcal{Y}} \in \mathbb{R}^{I_1 \times \cdots \times I_N}$ is denoted by $\langle \boldsymbol{\mathcal{X}}, \boldsymbol{\mathcal{Y}} \rangle = \sum_{i_1 \ldots i_N} x_{i_1 \cdots i_N} y_{i_1 \cdots i_N}$. The Frobenius norm of $\boldsymbol{\mathcal{X}}$ is denoted by $\|\boldsymbol{\mathcal{X}}\|_F = \sqrt{\langle \boldsymbol{\mathcal{X}}, \boldsymbol{\mathcal{X}} \rangle}$.

Let $\mathbf{\Phi} \in \mathbb{C}^{n_3 \times n_3}$ be a unitary transform matrix with $\mathbf{\Phi}\mathbf{\Phi}^H = \mathbf{\Phi}^H\mathbf{\Phi} = \boldsymbol{I}_{n_3}$, where $\boldsymbol{I}_{n_3}$ is an $n_3 \times n_3$ identity matrix. Denote $\hat{\boldsymbol{\mathcal{Y}}}_{\mathbf{\Phi}}$ as the third-order tensor obtained by multiplying by $\mathbf{\Phi}$ on all slices of $\boldsymbol{\mathcal{Y}}$ as $\hat{\boldsymbol{\mathcal{Y}}}_{\mathbf{\Phi}}(i,j,:) = \mathbf{\Phi}(\boldsymbol{\mathcal{Y}}(i,j,:))$. For brevity, we denote $\hat{\boldsymbol{\mathcal{Y}}}_{\mathbf{\Phi}} = \mathbf{\Phi}[\boldsymbol{\mathcal{Y}}]$. A block diagonal matrix constructed from the frontal slices of $\hat{\boldsymbol{\mathcal{Y}}}_{\mathbf{\Phi}}$ is defined as

$$\overline{\boldsymbol{\mathcal{Y}}}_{\mathbf{\Phi}} = \text{blockdiag}\left(\hat{\boldsymbol{\mathcal{Y}}}_{\mathbf{\Phi}}\right) := \begin{pmatrix} \hat{\boldsymbol{\mathcal{Y}}}_{\mathbf{\Phi}}^{(1)} & & & \\ & \hat{\boldsymbol{\mathcal{Y}}}_{\mathbf{\Phi}}^{(2)} & & \\ & & \ddots & \\ & & & \hat{\boldsymbol{\mathcal{Y}}}_{\mathbf{\Phi}}^{(n_3)} \end{pmatrix}.$$

For a block-diagonal matrix $\boldsymbol{M} \in \mathbb{C}^{n_1 n_3 \times n_2 n_3}$ consisting of $n_3$ blocks $\boldsymbol{M}^{(i)} \in \mathbb{C}^{n_1 \times n_2}$, the operator fold($\boldsymbol{M}$) assembles these blocks into a third-order tensor $\boldsymbol{\mathcal{A}} \in \mathbb{C}^{n_1 \times n_2 \times n_3}$ such that the $i$-th frontal slice of $\boldsymbol{\mathcal{A}}$ corresponds to the $i$-th block of $\boldsymbol{M}$, i.e., $\boldsymbol{\mathcal{A}}(:,:,i) = \boldsymbol{M}^{(i)}$.

The $\mathbf{\Phi}$-product of $\mathcal{B} \in \mathbb{C}^{n_1 \times n_2 \times n_3}$ and $\boldsymbol{\mathcal{C}} \in \mathbb{C}^{n_1 \times n_2 \times n_3}$ is a tensor $\boldsymbol{\mathcal{A}} \in \mathbb{C}^{n_1 \times n_2 \times n_3}$, given by

$$\boldsymbol{\mathcal{A}} = \mathcal{B} *_{\mathbf{\Phi}} \boldsymbol{\mathcal{C}} = \mathbf{\Phi}^H \left[\text{fold}\left(\text{blockdiag}\left(\hat{\mathcal{B}}_{\mathbf{\Phi}}\right) \times \text{blockdiag}\left(\hat{\boldsymbol{\mathcal{C}}}_{\mathbf{\Phi}}\right)\right)\right].$$

The conjugate transpose of $\boldsymbol{\mathcal{Y}} \in \mathbb{C}^{n_1 \times n_2 \times n_3}$ with respect to $\mathbf{\Phi}$, denoted by $\boldsymbol{\mathcal{Y}}^H \in \mathbb{C}^{n_2 \times n_1 \times n_3}$, is defined as

$$\boldsymbol{\mathcal{Y}}^H = \mathbf{\Phi}^H \left[\text{fold}\left(\text{blockdiag}\left(\hat{\boldsymbol{\mathcal{Y}}}_{\mathbf{\Phi}}\right)^H\right)\right].$$

The tensor $\boldsymbol{\mathcal{Y}}$ is unitary with respect to $\mathbf{\Phi}$-product if $\boldsymbol{\mathcal{Y}}^H *_{\mathbf{\Phi}} \boldsymbol{\mathcal{Y}} = \boldsymbol{\mathcal{Y}} *_{\mathbf{\Phi}} \boldsymbol{\mathcal{Y}}^H = \boldsymbol{\mathcal{I}}_{\mathbf{\Phi}}$, where $\boldsymbol{\mathcal{I}}_{\mathbf{\Phi}} \in \mathbb{C}^{n \times n \times n_3}$ is the identity tensor with each slice being identity matrix.

Definition 2.1 (Transformed t-SVD [40]). *The transformed t-SVD of tensor* $\boldsymbol{\mathcal{Y}}$ *is defined as*

$$\boldsymbol{\mathcal{Y}} = \boldsymbol{\mathcal{U}} *_{\mathbf{\Phi}} \boldsymbol{\mathcal{D}} *_{\mathbf{\Phi}} \boldsymbol{\mathcal{V}}^H,$$

*where* $\boldsymbol{\mathcal{U}} \in \mathbb{C}^{n_1 \times n_1 \times n_3}$, $\boldsymbol{\mathcal{V}} \in \mathbb{C}^{n_2 \times n_2 \times n_3}$ *are unitary tensors with respect to the* $\mathbf{\Phi}$*-product, and* $\boldsymbol{\mathcal{D}} \in \mathbb{C}^{n_1 \times n_2 \times n_3}$ *is a diagonal tensor.*

Definition 2.2 (The transformed tubal nuclear norm [40]). *The transformed tubal nuclear norm (TTNN) of* $\boldsymbol{\mathcal{Y}}$ *is the sum of the nuclear norms of all frontal slices of* $\hat{\boldsymbol{\mathcal{Y}}}_{\mathbf{\Phi}}$*, defined as*

$$\|\boldsymbol{\mathcal{Y}}\|_{TTNN} = \sum_{i=1}^{n_3} \left\|\hat{\boldsymbol{\mathcal{Y}}}_{\mathbf{\Phi}}^{(i)}\right\|_*.$$

Let $f : \mathbb{R}^n \to \mathbb{R}_\infty := \mathbb{R} \cup \{+\infty\}$ be a proper closed function. The function $f$ is said to be coercive if $f(\boldsymbol{x}) \to +\infty$ as $\|\boldsymbol{x}\| \to +\infty$. Given $\boldsymbol{x} \in \mathbb{R}^n$ and $t > 0$, the proximity associated to $f$ and its corresponding Moreau proximal envelope are defined respectively by

$$\text{prox}_{tf}(\boldsymbol{x}) := \underset{\boldsymbol{u} \in \mathbb{R}^n}{\text{argmin}} \left\{ f(\boldsymbol{u}) + \frac{1}{2t}\|\boldsymbol{u} - \boldsymbol{x}\|^2 \right\},$$

and

$$f_t(\boldsymbol{x}) := \min_{\boldsymbol{u} \in \mathbb{R}^n} \left\{ f(\boldsymbol{u}) + \frac{1}{2t}\|\boldsymbol{u} - \boldsymbol{x}\|^2 \right\}.$$

Lemma 2.3 (The proximity of TTNN [40]). *For any* $\boldsymbol{\mathcal{Y}} = \boldsymbol{\mathcal{U}} *_{\boldsymbol{\Phi}} \boldsymbol{\mathcal{S}} *_{\boldsymbol{\Phi}} \boldsymbol{\mathcal{V}}^H$*, the proximity of* $\|\cdot\|_{TTNN}$ *is given by*

$$\begin{aligned} prox_{\lambda\|\cdot\|_{TTNN}}(\boldsymbol{\mathcal{Y}}) &= \underset{\boldsymbol{\mathcal{X}}}{argmin} \left\{ \lambda \|\boldsymbol{\mathcal{X}}\|_{TTNN} + \frac{1}{2}\|\boldsymbol{\mathcal{X}} - \boldsymbol{\mathcal{Y}}\|_F^2 \right\} \\ &= \boldsymbol{\mathcal{U}} *_{\boldsymbol{\Phi}} \boldsymbol{\mathcal{D}}_\lambda *_{\boldsymbol{\Phi}} \boldsymbol{\mathcal{V}}^H, \end{aligned}$$

*where* $\boldsymbol{\mathcal{D}}_\lambda = \boldsymbol{\Phi}^H \left[ \hat{\boldsymbol{\mathcal{D}}}_\lambda \right]$ *and* $\hat{\boldsymbol{\mathcal{D}}}_\lambda = \max \left\{ \hat{\boldsymbol{\mathcal{D}}}_{\boldsymbol{\Phi}} - \lambda, 0 \right\}$.

Lemma 2.4 (L-smoothness [7]). *Let* $f$ *be a continuously differentiable function. If its gradient* $\nabla f$ *is Lipschitz continuous with constant* $\ell$*, then the following inequality holds for any* $\boldsymbol{x}, \boldsymbol{y} \in \mathbb{R}^n$*:*

$$f(\boldsymbol{y}) \leq f(\boldsymbol{x}) + \langle \nabla f(\boldsymbol{x}), \boldsymbol{y} - \boldsymbol{x} \rangle + \frac{\ell}{2}\|\boldsymbol{y} - \boldsymbol{x}\|^2.$$

Lemma 2.5 (Differentiability of the Moreau envelope [37]). *Let* $f : \mathbb{R}^n \to \mathbb{R}_\infty$ *be a proper, closed, and convex function. Then, its Moreau envelope* $f_t(\cdot)$ *is continuously differentiable for any* $t > 0$*. Furthermore, its gradient* $\nabla f_t(x)$ *is* $\frac{1}{t}$*-Lipschitz continuous and is given by*

$$\nabla f_t(x) = \frac{x - prox_{tf}(x)}{t}.$$

Lemma 2.6 (Moreau envelope of the indicator function [37]). *The Moreau envelope of the indicator function of a set* $\Omega$ *is the distance function, i.e.,*

$$\mathcal{I}_t(\boldsymbol{x}) = \min_{\boldsymbol{u}} \left\{ \mathcal{I}_\Omega(\boldsymbol{u}) + \frac{1}{2t}\|\boldsymbol{u} - \boldsymbol{x}\|^2 \right\} = \frac{1}{2t} dist^2(\boldsymbol{x}, \Omega).$$

*Moreover, if* $\Omega$ *is a convex set, then*

$$\nabla \mathcal{I}_t(\boldsymbol{x}) = \frac{1}{t}\left(\boldsymbol{x} - P_\Omega(\boldsymbol{x})\right).$$

**3. Unified tensor fusion framework and modeling.** In this section, we formulate blind fusion as a coupled inverse problem, integrating blind deconvolution in the spatial domain with blind unmixing in the spectral domain. Building on this perspective, we develop a unified tensor framework, followed by a corresponding proposed model.

**3.1. Blind deconvolution meets blind unmixing.** The target HR-HSI preserves the spectral resolution of the input HSI while overcoming the reduced spatial resolution. This degradation process can be modeled as an unknown two-dimensional convolution operation applied to each spectral band, followed by spatial downsampling. Notably, our current study focuses specifically on separable two-dimensional blur operators. In the spatial domain, this constitutes a blind deconvolution problem.

Unlike the HSI, the MSI maintains the original spatial resolution while exhibiting reduced spectral resolution. Each MSI band can be mathematically expressed as a linear mixture of corresponding spectral bands from the target HR-HSI. Specifically, the MSI is generated by multiplying the HR-HSI with an unknown spectral response matrix. This formulation naturally leads to a blind unmixing problem in the spectral domain.

Fundamentally, the blind fusion problem comprises two key components: blind deconvolution in the spatial domain and blind unmixing in the spectral domain. This constitutes a coupled inverse problem that necessitates a joint solution across both dimensions.

**3.2. Tensor data-fitting terms.** We denote $\boldsymbol{\mathcal{H}} \in \mathbb{R}^{I' \times J' \times K}$ and $\boldsymbol{\mathcal{M}} \in \mathbb{R}^{I \times J \times K'}$ as the given HSI and MSI, respectively. In practice, an HSI captures a broad range of spectrum, involving hundreds of spectral bands/wavelengths whilst an MSI is measured at less than 20 bands/wavelengths, i.e., $K' \ll K$. On the other hand, from the spatial perspective, an MSI has finer resolution than the HSI, i.e., $I' \ll I$ and $J' \ll J$. The HR-HSI is represented as the tensor $\boldsymbol{\mathcal{S}} \in \mathbb{R}^{I \times J \times K}$. Mathematically, $\boldsymbol{\mathcal{H}}$ and $\boldsymbol{\mathcal{M}}$ are derived from $\boldsymbol{\mathcal{S}}$ by the tensor product as follows:

$$\begin{aligned} \boldsymbol{\mathcal{H}} &= \boldsymbol{\mathcal{S}} \times_1 \boldsymbol{P}_1 \times_2 \boldsymbol{P}_2 + \boldsymbol{\mathcal{N}}_h, \\ \boldsymbol{\mathcal{M}} &= \boldsymbol{\mathcal{S}} \times_3 \boldsymbol{P}_3 + \boldsymbol{\mathcal{N}}_m, \end{aligned} \tag{3.1}$$

where $\boldsymbol{P}_1 \in \mathbb{R}^{I' \times I}$ and $\boldsymbol{P}_2 \in \mathbb{R}^{J' \times J}$ are the spatial degradation operators acting on width and height modes, respectively, and $\boldsymbol{P}_3 \in \mathbb{R}^{K' \times K}$ is the spectral degradation operator. $\boldsymbol{\mathcal{N}}_h$ and $\boldsymbol{\mathcal{N}}_m$ denote the additive noise for $\boldsymbol{\mathcal{H}}$ and $\boldsymbol{\mathcal{M}}$, respectively.

Equation (3.1) reflects the underlying degradation mechanism and constitutes two tensor data-fitting terms of the inverse multi-dimensional imaging problem. Therefore, we adopt the nonlinear least squares optimization to recover $\boldsymbol{\mathcal{S}}$ as follows

$$\min_{\boldsymbol{\mathcal{S}}, \boldsymbol{P}_1, \boldsymbol{P}_2, \boldsymbol{P}_3} \frac{1}{2} \|\boldsymbol{\mathcal{H}} - \boldsymbol{\mathcal{S}} \times_1 \boldsymbol{P}_1 \times_2 \boldsymbol{P}_2\|_F^2 + \frac{\lambda_1}{2} \|\boldsymbol{\mathcal{M}} - \boldsymbol{\mathcal{S}} \times_3 \boldsymbol{P}_3\|_F^2, \tag{3.2}$$

where $\lambda_1 > 0$ is the trade-off parameter. To tackle the ill-posedness of the inverse problem, regularizations and constraints reflecting prior knowledge of target image $\boldsymbol{\mathcal{S}}$ and three degradation operators $\boldsymbol{P}_{1,2,3}$ should be incorporated.

**3.3. Structural constraints for PSF and SRF.** Intrinsically, the spatial degradation includes three primary components: alignment, blurring, and downsampling. This paper assumes that the given HSIs are accurately aligned, thereby allowing us to focus exclusively on the effects of blurring and downsampling. The blurring process, which is assumed to be band-independent, can be mathematically represented as a linear combination of neighboring pixels through a specific convolution kernel operating in both spatial dimensions (width and height). Downsampling is modeled as a uniform subsampling operation applied to the image. Consequently, the spatial degradation operator can be formally expressed as the composition of two sequential operations: the blurring operator followed by the downsampling operator.

Thus, the spatial degradation operator can generally be formulated as

$$\boldsymbol{P} = \boldsymbol{D}\boldsymbol{B}. \tag{3.3}$$

Here, the downsampling operator $\boldsymbol{D} \in \mathbb{R}^{m \times n}$ consists of a subset of the rows of the identity matrix, which is determined by the spatial resolution ratio between the observed HSI-MSI pair. $\boldsymbol{B} \in \mathbb{R}^{n\times n}$ corresponds to the blurring operator associated with the blur kernel $\boldsymbol{b} \in \mathbb{R}^n$. It is a typical circulant matrix of the form

$$\boldsymbol{B} = \begin{pmatrix} b_0 & b_{n-1} & \dots & b_2 & b_1 \\ b_1 & b_0 & b_{n-1} & & b_2 \\ \vdots & b_1 & b_0 & \ddots & \vdots \\ b_{n-2} & & \ddots & \ddots & b_{n-1} \\ b_{n-1} & b_{n-2} & \cdots & b_1 & b_0 \end{pmatrix}. \tag{3.4}$$

Note that the first column of $\boldsymbol{B}$ is given by $\boldsymbol{b}$. As stated in [8], $\boldsymbol{B}$ admits the eigendecomposition as

$$\boldsymbol{B} = \boldsymbol{F}\boldsymbol{S}\boldsymbol{F}^*, \tag{3.5}$$

where $\boldsymbol{F} \in \mathbb{C}^{n\times n}$ is a normalized discrete Fourier transform (DFT) matrix ($\boldsymbol{F}^*\boldsymbol{F} = \boldsymbol{F}\boldsymbol{F}^* = \boldsymbol{I}_n$), $\boldsymbol{S} \in \mathbb{R}^{n\times n}$ is a diagonal matrix and $*$ denotes the conjugate transpose operator. Furthermore, for a circulant matrix, its eigenvalues can be easily computed by taking the Fourier transform of its first column, namely,

$$\operatorname{diag}(\boldsymbol{S}) = \sqrt{n}\boldsymbol{F}\boldsymbol{b}.$$

In summary, the spatial degradation operator has the general form,

$$\boldsymbol{P} = \boldsymbol{D}\boldsymbol{F}\operatorname{Diag}(\sqrt{n}\boldsymbol{F}\boldsymbol{b})\boldsymbol{F}^*, \tag{3.6}$$

where $\boldsymbol{b}$ denotes the unknown blur kernel vector for PSF that needs to be estimated. This formulation thus leads to a blind deconvolution problem.

The spectral degradation operator is commonly represented by a spectral response matrix that operates on the shared wavelength ranges between the targeted HR-HSI and MSI. More specifically, each MSI band can be modeled as a linear combination of its neighboring spectral bands of the HR-HSI. Mathematically, each row of the spectral response matrix $\boldsymbol{P}_3$ can be expressed as

$$\boldsymbol{P}_3(i,:) = (\boldsymbol{D}_3^i \boldsymbol{b}_3^i)^\top, \quad \boldsymbol{b}_3^i \in \mathbb{S}_{k_i} \quad i = 1, \dots, K'. \tag{3.7}$$

Here, $\boldsymbol{D}_3^i$, the spectral response index matrix for the $i$-th MSI band, is formed by selecting specific rows from the identity matrix. It is determined by comparing the wavelength ranges of the HSI and MSI bands. The mixing vector $\boldsymbol{b}_3^i$, subject to a particular constraint, remains an unknown SRF parameter to be estimated.

**3.4. The unified tensor fusion framework.** In summary, we present the unified tensor framework as the following minimization,

$$\begin{aligned} \min_{\boldsymbol{\mathcal{S}},\boldsymbol{P}_1,\boldsymbol{P}_2,\boldsymbol{P}_3} \; & \frac{1}{2}\|\boldsymbol{\mathcal{H}} - \boldsymbol{\mathcal{S}} \times_1 \boldsymbol{P}_1 \times_2 \boldsymbol{P}_2\|_F^2 + \frac{\lambda_1}{2}\|\boldsymbol{\mathcal{M}} - \boldsymbol{\mathcal{S}} \times_3 \boldsymbol{P}_3\|_F^2 + \varphi(\boldsymbol{\mathcal{S}}) \\ \text{s.t.} \quad & \boldsymbol{P}_1 \in \Phi_1, \quad \boldsymbol{P}_2 \in \Phi_2, \quad \boldsymbol{P}_3 \in \Phi_3, \end{aligned} \tag{3.8}$$

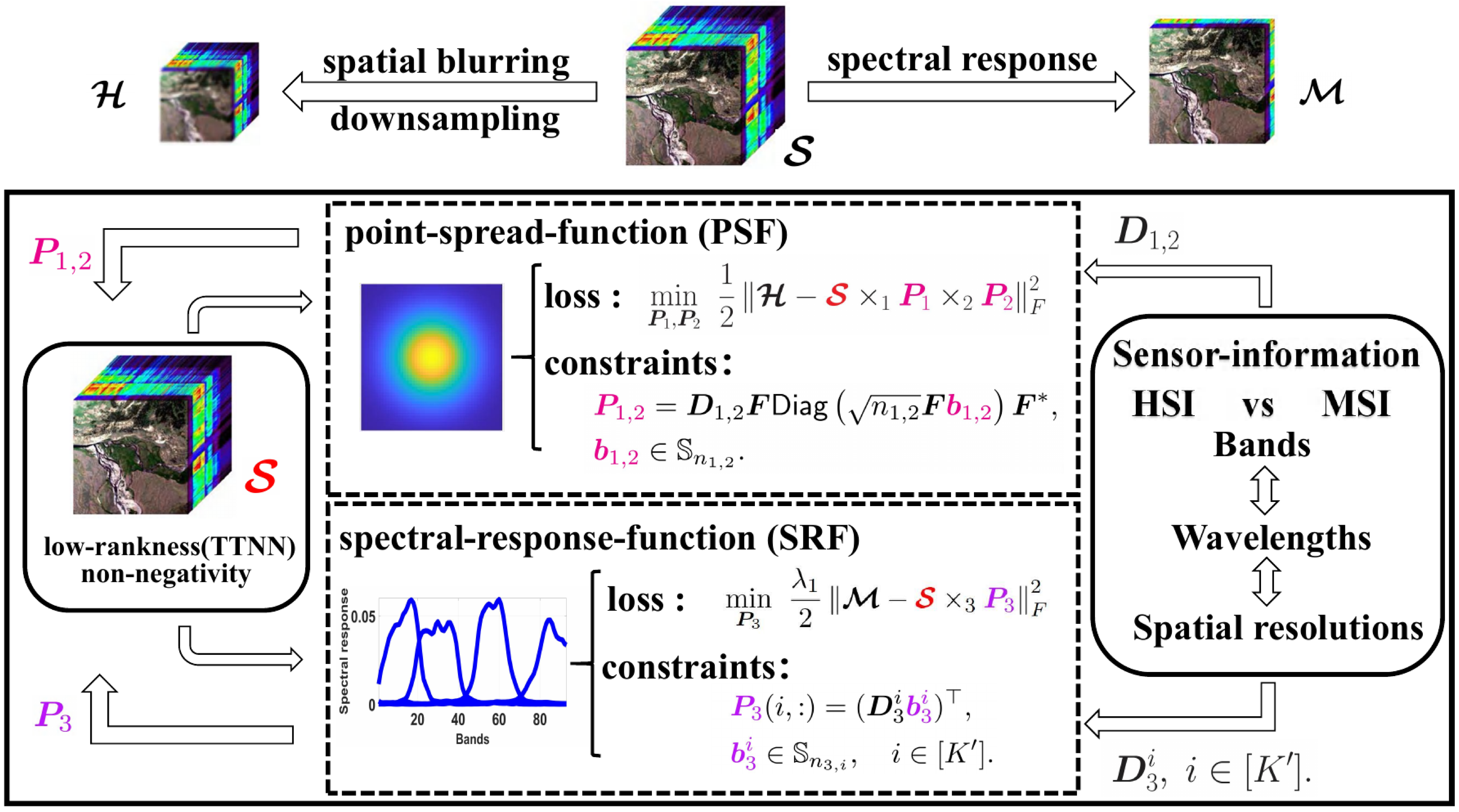

**Figure 3.1.** *The unified tensor fusion framework and proposed model* (3.10).

where

$$\Phi_1 = \left\{ \boldsymbol{P}_1 \mid \boldsymbol{P}_1 = \boldsymbol{D}_1 \boldsymbol{F}_1 \mathrm{Diag}(\sqrt{I}\boldsymbol{F}_1 \boldsymbol{b}_1)\boldsymbol{F}_1^*, \boldsymbol{b}_1 \in \mathbb{S}_1 \right\},$$
$$\Phi_2 = \left\{ \boldsymbol{P}_2 \mid \boldsymbol{P}_2 = \boldsymbol{D}_2 \boldsymbol{F}_2 \mathrm{Diag}(\sqrt{J}\boldsymbol{F}_2 \boldsymbol{b}_2)\boldsymbol{F}_2^*, \boldsymbol{b}_2 \in \mathbb{S}_2 \right\},$$
$$\Phi_3 = \left\{ \boldsymbol{P}_3 \mid \boldsymbol{P}_3(i,:) = (\boldsymbol{D}_3^i \boldsymbol{b}_3^i)^\top, \ \ \boldsymbol{b}_3^i \in \mathbb{S}_3^i, \ i = 1, \dots, K' \right\},$$

where $\lambda_1 > 0$ is the trade-off parameter. The downsampling operators $\boldsymbol{D}_1$, $\boldsymbol{D}_2$, and $\boldsymbol{D}_3^i$ ($i = 1, \dots, K'$) are provided by the sensor information. The regularization term $\varphi(\cdot)$ for $\boldsymbol{\mathcal{S}}$ incorporates the priors for the target HR-HSI, while the constraint sets $\mathbb{S}_1$, $\mathbb{S}_2$, and $\mathbb{S}_3^i$ characterize the specific feasible regions for unknown blur kernels and mixing vectors, thereby imposing constraints on the variables $\boldsymbol{P}_{1,2,3}$. The overall architecture of the proposed unified tensor fusion framework is illustrated in Figure 3.1.

*Remark* 3.1. We elaborate on the limitations, potential extensions, and refinement strategies for the proposed framework as follows.

(a) Boundary effects: The circulant convolution in $\boldsymbol{P}_{1,2}$ assumes periodic boundaries for frequency-domain efficiency. This may induce edge artifacts near high-contrast structures, which can be alleviated via symmetric padding, edge tapering, or patch-based processing.
(b) PSF invariance: The model postulates a band-independent PSF to reduce complexity. While effective for stable systems, this assumption may deviate from reality when the PSF exhibits significant chromatic aberration or wavelength-dependent variation.

(c) SRF truncation: The hard-window local mixing ignores the long-tail characteristics of realistic SRFs, potentially causing spectral bias. To address this, a Tikhonov-type regularization $\lambda \sum_{i=1}^{K'} \|\mathbf{D}_3^i \boldsymbol{b}_3^i\|_2^2$ can be integrated to promote weight decay. Alternatively, incorporating an $\ell_1$ penalty into sparse non-local mixing promotes an adaptive "soft" window, thereby enhancing the precision of spectral alignment.

(d) Noise modeling: While the current i.i.d. Gaussian assumption facilitates theoretical analysis and algorithm design, it may be suboptimal for physical sensors characterized by complex Poisson-Gaussian noise. Consequently, the incorporation of noise-adaptive weighting or band-dependent modeling constitutes a promising extension of the proposed framework to enhance its robustness against non-stationary degradation.

**3.5. Our proposed model.** For the target image $\boldsymbol{\mathcal{S}}$, the inherent low-rankness and non-negativity are crucial structural properties. Owing to the superiority of transformed t-SVD, we incorporate the TTNN regularization [40] of $\boldsymbol{\mathcal{S}}$ to capture its low-rank structure. The regularization term involving the indicator function on the non-negativity set is also included. Moreover, to address the inherent non-uniqueness of blind deconvolution, a widely adopted approach (see, e.g., [26, 42, 2]) is to normalize the kernel with respect to the 1-norm. This normalization, combined with the non-negativity constraint, naturally restricts the kernel to reside within the unit simplex, i.e.,

$$\boldsymbol{b} \in \mathbb{S}_n := \left\{ \boldsymbol{b} \in \mathbb{R}^n \mid b_i \geq 0, \ \sum_{i=1}^{n} b_i = 1 \right\}. \tag{3.9}$$

Hence, the fusion framework in (3.8) is reformulated as

$$\begin{aligned} \min_{\substack{\boldsymbol{\mathcal{S}},\boldsymbol{\mathcal{Y}}, \\ \boldsymbol{b}_1,\boldsymbol{b}_2,\boldsymbol{b}_3^i}} \quad & \frac{1}{2}\|\boldsymbol{\mathcal{H}} - \boldsymbol{\mathcal{S}} \times_1 \boldsymbol{P}_1 \times_2 \boldsymbol{P}_2\|_F^2 + \frac{\lambda_1}{2}\|\boldsymbol{\mathcal{M}} - \boldsymbol{\mathcal{S}} \times_3 \boldsymbol{P}_3\|_F^2 + \lambda_2 \|\boldsymbol{\mathcal{S}}\|_{\text{TTNN}} + \mathcal{I}(\boldsymbol{\mathcal{Y}}) \\ \text{s.t.} \quad & \boldsymbol{P}_1 = \boldsymbol{D}_1 \boldsymbol{F}_1 \mathrm{Diag}(\sqrt{I}\boldsymbol{F}_1 \boldsymbol{b}_1)\boldsymbol{F}_1^*, \qquad \boldsymbol{b}_1 \in \mathbb{S}_I, \\ & \boldsymbol{P}_2 = \boldsymbol{D}_2 \boldsymbol{F}_2 \mathrm{Diag}(\sqrt{J}\boldsymbol{F}_2 \boldsymbol{b}_2)\boldsymbol{F}_2^*, \qquad \boldsymbol{b}_2 \in \mathbb{S}_J, \\ & \boldsymbol{P}_3(i,:) = (\boldsymbol{D}_3^i \boldsymbol{b}_3^i)^\top, \qquad \boldsymbol{b}_3^i \in \mathbb{S}_{k_i}, \ i = 1, \dots, K', \\ & \boldsymbol{\mathcal{S}} = \boldsymbol{\mathcal{Y}}. \end{aligned} \tag{3.10}$$

Herein, the indicator function represents the non-negativity constraint as

$$\mathcal{I}(\boldsymbol{\mathcal{Y}}) := \begin{cases} 0, & \boldsymbol{\mathcal{Y}} \geq 0, \\ \infty, & \text{otherwise.} \end{cases}$$

*Remark* 3.2. The TTNN framework is flexible regarding the choice of the unitary transform matrix $\boldsymbol{\Phi}$ in Definition 2.2. We primarily employ a data-dependent transformation [40]. By performing a singular value decomposition on the tensor unfolding matrix along the spectral mode, we define $\boldsymbol{\Phi}$ as the conjugate transpose of the resulting left singular vector matrix. This approach facilitates a highly significant low-rank approximation tailored to the specific characteristics of the dataset.

Other choices for the transform matrix include: (a) Discrete Fourier Transform (DFT) Matrix [31]: This aligns with the conventional t-SVD framework, where the tensor product is defined via circular convolution in the Fourier domain; and (b) Discrete Wavelet Transform (DWT) Matrix [10]: By employing bases such as the Daubechies 4 (db4) wavelet, the model effectively captures localized, multi-scale features within the tensor.

*Remark* 3.3. The model (3.10) is actually a structured nonconvex optimization problem, which presents three significant theoretical and practical challenges. First, the selection of an appropriate starting point is crucial for nonconvex optimization, particularly for the variable $\boldsymbol{\mathcal{S}}$. Second, the algorithmic design poses substantial difficulties: (3.10) is a multiblock, nonseparable nonlinear problem with a linear equality constraint. This structure naturally suggests adopting a Gauss-Seidel-type framework applied to the augmented Lagrangian function. Nevertheless, the inherent nonlinearity of the tensor data-fitting terms renders the solution of the resulting subproblems computationally challenging. Third, and perhaps most critically, the convergence analysis presents substantial theoretical obstacles. Specifically, the nonsmoothness of the indicator functions poses additional challenges in establishing convergence.

## 4. Algorithm.

This section first presents the partially linearized ADMM with Moreau envelope smoothing to solve model (3.10). Then, a tailored initialization algorithm is devised based on the characteristics of the fusion problem.

### 4.1. A partially linearized ADMM with Moreau envelope smoothing.

ADMM is an effective splitting algorithm designed for solving linearly constrained optimization problems. The canonical ADMM minimizes the augmented Lagrangian function associated with problem (3.10) cyclically over $\boldsymbol{\mathcal{S}}$, $\boldsymbol{\mathcal{Y}}$, $\boldsymbol{\Lambda}$, $\boldsymbol{P}_1$, $\boldsymbol{P}_2$, and $\boldsymbol{P}_3$ while fixing the remaining blocks at the last updated values. However, solving the $\boldsymbol{\mathcal{S}}$, $\boldsymbol{P}_1$, $\boldsymbol{P}_2$, and $\boldsymbol{P}_3$ subproblems presents significant challenges due to the inherent nonlinearity of the two tensor data-fitting terms.

Motivated by the aforementioned considerations, we devise a partially linearized ADMM to solve (3.10). More specifically, we linearize the loss functions $L_1$ and $L_2$, which facilitates an efficient solver for the $\boldsymbol{\mathcal{S}}$-subproblem. We adopt Moreau envelope smoothing to address the nonsmoothness of the indicator function. For $\boldsymbol{P}_{1,2,3}$-subproblems, we reformulate the model with respect to variables $\boldsymbol{b}_{1,2,3}$ and adopt the linearization. We now present our partially proximal linearized ADMM.

For denotation convenience, variables $\boldsymbol{b}_1$, $\boldsymbol{b}_2$, and $\boldsymbol{b}_3^i$ are transformed into $\boldsymbol{P}_1$, $\boldsymbol{P}_2$, and $\boldsymbol{P}_3$, respectively. The Moreau envelope augmented Lagrangian function associated with problem (3.10) is presented as

$$\begin{aligned}\mathcal{L}_{\beta,\mu}(\boldsymbol{\mathcal{S}},\boldsymbol{\mathcal{Z}},\boldsymbol{\Lambda},\boldsymbol{P}_1,\boldsymbol{P}_2,\boldsymbol{P}_3) :=& L_1(\boldsymbol{\mathcal{S}},\boldsymbol{P}_1,\boldsymbol{P}_2)+L_2(\boldsymbol{\mathcal{S}},\boldsymbol{P}_3)+\lambda_2\|\boldsymbol{\mathcal{S}}\|_{\mathrm{TTNN}}+\mathcal{I}_\mu(\boldsymbol{\mathcal{Z}})\\ &+\langle\boldsymbol{\Lambda},\boldsymbol{\mathcal{S}}-\boldsymbol{\mathcal{Z}}\rangle+\frac{\beta}{2}\|\boldsymbol{\mathcal{S}}-\boldsymbol{\mathcal{Z}}\|_F^2,\end{aligned}$$

where

$$L_1(\boldsymbol{\mathcal{S}},\boldsymbol{P}_1,\boldsymbol{P}_2):=\frac{1}{2}\left\|\boldsymbol{\mathcal{H}}-\boldsymbol{\mathcal{S}}\times_1\boldsymbol{P}_1\times_2\boldsymbol{P}_2\right\|_F^2,\quad L_2(\boldsymbol{\mathcal{S}},\boldsymbol{P}_3):=\frac{\lambda_1}{2}\left\|\boldsymbol{\mathcal{M}}-\boldsymbol{\mathcal{S}}\times_3\boldsymbol{P}_3\right\|_F^2, \tag{4.1}$$

are the data-fitting terms that measure the approximation in (3.1), and

$$\mathcal{I}_\mu(\boldsymbol{\mathcal{Z}}) = \underset{\boldsymbol{\mathcal{Y}}}{\operatorname{argmin}}\, \mathcal{I}(\boldsymbol{\mathcal{Y}}) + \frac{1}{2\mu}\|\boldsymbol{\mathcal{Y}} - \boldsymbol{\mathcal{Z}}\|_F^2,$$

is the Moreau envelope of the indicator function for the nonnegative orthant $\Omega$.

For notational convenience, define

$$L(\boldsymbol{\mathcal{S}}, \boldsymbol{P}_{1,2,3}) := L_1(\boldsymbol{\mathcal{S}}, \boldsymbol{P}_1, \boldsymbol{P}_2) + L_2(\boldsymbol{\mathcal{S}}, \boldsymbol{P}_3).$$

**Updating variable $\boldsymbol{\mathcal{S}}$:**

$$\boldsymbol{\mathcal{S}}^{t+1} = \underset{\boldsymbol{\mathcal{S}}}{\operatorname{argmin}}\Big\{\lambda_2\|\boldsymbol{\mathcal{S}}\|_{\text{TTNN}} + \langle \boldsymbol{\Lambda}^t, \boldsymbol{\mathcal{S}} - \boldsymbol{\mathcal{Z}}^t\rangle + \langle \nabla_{\boldsymbol{\mathcal{S}}} L(\boldsymbol{\mathcal{S}}^t, \boldsymbol{P}_{1,2,3}^t), \boldsymbol{\mathcal{S}} - \boldsymbol{\mathcal{S}}^t\rangle + \frac{\alpha^t}{2}\|\boldsymbol{\mathcal{S}} - \boldsymbol{\mathcal{S}}^t\|_F^2 + \frac{\beta}{2}\|\boldsymbol{\mathcal{S}} - \boldsymbol{\mathcal{Z}}^t\|_F^2\Big\},$$

where $\alpha^t$ is the varying stepsize determined by the Lipschitz constant of $\nabla_{\boldsymbol{\mathcal{S}}} L(\cdot, \boldsymbol{P}_{1,2,3}^t)$. Accordingly, by the definition of proximity in Theorem 2.3, the closed-form solution of $\boldsymbol{\mathcal{S}}$-subproblem can be formalized as

$$\begin{cases} \tilde{\boldsymbol{\mathcal{S}}}^t = \dfrac{\alpha^t}{\alpha^t + \beta}\left(\boldsymbol{\mathcal{S}}^t - \dfrac{1}{\alpha^t}\nabla_{\boldsymbol{\mathcal{S}}} L(\boldsymbol{\mathcal{S}}^t, \boldsymbol{P}_{1,2,3}^t)\right), \\ \tilde{\boldsymbol{\mathcal{Z}}}^t = \dfrac{\beta}{\alpha^t + \beta}\left(\boldsymbol{\mathcal{Z}}^t - \dfrac{1}{\beta}\boldsymbol{\Lambda}^t\right), \\ \boldsymbol{\mathcal{S}}^{t+1} = \operatorname{prox}_{\frac{\lambda_2}{\alpha^t+\beta}\|\cdot\|_{\text{TTNN}}}\left(\tilde{\boldsymbol{\mathcal{S}}}^t + \tilde{\boldsymbol{\mathcal{Z}}}^t\right). \end{cases} \tag{4.2}$$

**Updating variable $\boldsymbol{\mathcal{Z}}$:**

$$\boldsymbol{\mathcal{Z}}^{t+1} = \underset{\boldsymbol{\mathcal{Z}}}{\operatorname{argmin}}\Big\{\mathcal{I}_\mu(\boldsymbol{\mathcal{Z}}) + \langle \boldsymbol{\Lambda}^t, \boldsymbol{\mathcal{S}}^{t+1} - \boldsymbol{\mathcal{Z}}\rangle + \frac{\beta}{2}\|\boldsymbol{\mathcal{S}}^{t+1} - \boldsymbol{\mathcal{Z}}\|_F^2\Big\}. \tag{4.3}$$

Lemma 4.1. *The solution of subproblem (4.3) is given by*

$$\boldsymbol{\mathcal{Y}}^{t+1} = P_\Omega(\boldsymbol{\mathcal{S}}^{t+1} + \frac{1}{\beta}\boldsymbol{\Lambda}^t), \tag{4.4a}$$

$$\boldsymbol{\mathcal{Z}}^{t+1} = \frac{\mu}{1+\mu\beta}\left(\frac{1}{\mu}\boldsymbol{\mathcal{Y}}^{t+1} + \boldsymbol{\Lambda}^t + \beta\boldsymbol{\mathcal{S}}^{t+1}\right). \tag{4.4b}$$

*Proof.* The $\boldsymbol{\mathcal{Z}}$-subproblem in (4.3) is equivalent to the following minimization,

$$(\boldsymbol{\mathcal{Y}}^{t+1}, \boldsymbol{\mathcal{Z}}^{t+1}) = \underset{\boldsymbol{\mathcal{Y}},\, \boldsymbol{\mathcal{Z}}}{\operatorname{argmin}}\Big\{\mathcal{I}(\boldsymbol{\mathcal{Y}}) + \frac{1}{2\mu}\|\boldsymbol{\mathcal{Y}} - \boldsymbol{\mathcal{Z}}\|_F^2 + \langle \boldsymbol{\Lambda}^t, \boldsymbol{\mathcal{S}}^{t+1} - \boldsymbol{\mathcal{Z}}\rangle + \frac{\beta}{2}\|\boldsymbol{\mathcal{S}}^{t+1} - \boldsymbol{\mathcal{Z}}\|_F^2\Big\}.$$

The optimality condition of the above problem is

$$0 = \frac{1}{\mu}\left(\boldsymbol{\mathcal{Z}}^{t+1} - \boldsymbol{\mathcal{Y}}^{t+1}\right) - \boldsymbol{\Lambda}^t + \beta\left(\boldsymbol{\mathcal{S}}^{t+1} - \boldsymbol{\mathcal{Z}}^{t+1}\right), \tag{4.5a}$$

$$0 \in \partial\mathcal{I}(\boldsymbol{\mathcal{Y}}^{t+1}) + \frac{1}{\mu}(\boldsymbol{\mathcal{Y}}^{t+1} - \boldsymbol{\mathcal{Z}}^{t+1}), \tag{4.5b}$$

(4.5a) immediately yields (4.4b). Plugging (4.4b) into (4.5b) yields

$$0 \in \frac{1+\mu\beta}{\beta}\partial\mathcal{I}(\boldsymbol{\mathcal{Y}}^{t+1}) + \boldsymbol{\mathcal{Y}}^{t+1} - \left(\boldsymbol{\mathcal{S}}^{t+1} - \frac{1}{\beta}\boldsymbol{\Lambda}^t\right),$$

which implies that

$$\boldsymbol{\mathcal{Y}}^{t+1} = \text{prox}_{\frac{1+\mu\beta}{\beta}\mathcal{I}}\left(\boldsymbol{\mathcal{S}}^{t+1} + \frac{1}{\beta}\boldsymbol{\Lambda}^t\right). \qquad \blacksquare$$

**Updating multiplier variable $\boldsymbol{\Lambda}$:**

$$\boldsymbol{\Lambda}^{t+1} = \boldsymbol{\Lambda}^t + \beta(\boldsymbol{\mathcal{S}}^{t+1} - \boldsymbol{\mathcal{Y}}^{t+1}). \tag{4.6}$$

**Updating variables $\boldsymbol{P}_{1,2}$:** to facilitate the upcoming subproblems, we first formulate the generalized form of $\boldsymbol{P}_{1,2,3}$ subproblems as follows

$$\min_{\boldsymbol{P}\in\Phi} \frac{1}{2}\|\boldsymbol{P}\boldsymbol{A} - \boldsymbol{W}\|_F^2,$$

where $\Phi := \{\boldsymbol{P} \in \mathbb{R}^{m\times n} \mid \boldsymbol{P} = \boldsymbol{D}\boldsymbol{F}\text{Diag}(\sqrt{n}\boldsymbol{F}\boldsymbol{b})\boldsymbol{F}^*, \boldsymbol{b} \in \mathbb{S}_n\}$. Furthermore, we reformulate the above model through substituting the variable $\boldsymbol{P}$ with the variable $\boldsymbol{b}$ as follows,

$$\min_{\boldsymbol{b}\in\mathbb{S}_n} f(\boldsymbol{b}) = \frac{1}{2}\|\boldsymbol{X}(\boldsymbol{b})\|_F^2, \tag{4.7}$$

where $\boldsymbol{X}(\boldsymbol{b}) := \boldsymbol{D}\boldsymbol{F}\text{Diag}(\sqrt{n}\boldsymbol{F}\boldsymbol{b})\boldsymbol{F}^*\boldsymbol{A} - \boldsymbol{W}$ and $\mathbb{S}_n$ is the unit simplex constraint defined in (3.9). The minimization model (4.7) is a least squares problem with the unit simplex constraint, which inherently exhibits convexity. The projected gradient method is adopted to solve (4.7) as follows,

$$\boldsymbol{b}^{t+1} = P_{\mathbb{S}_n}\left(\boldsymbol{b}^t - \tau^t\nabla f(\boldsymbol{b}^t)\right), \tag{4.8}$$

where $P_{\mathbb{S}_n}(\cdot)$ denotes the projection on the unit simplex $\mathbb{S}_n$, $\tau^t$ is the varying stepsize determined by the Lipschitz constant of $\nabla f(\boldsymbol{b}^t)$ and

$$\begin{aligned} \boldsymbol{Q}^t &= \boldsymbol{F}^*\boldsymbol{D}^\top\boldsymbol{X}(\boldsymbol{b}^t)\boldsymbol{A}^\top\boldsymbol{F}, \\ \nabla f(\boldsymbol{b}^t) &= \boldsymbol{F}^*\text{diag}(\boldsymbol{Q}^t). \end{aligned} \tag{4.9}$$

Furthermore, we derive

$$\boldsymbol{P}^{t+1} = \boldsymbol{D}\boldsymbol{F}\text{Diag}(\sqrt{n}\boldsymbol{F}\boldsymbol{b}^{t+1})\boldsymbol{F}^*. \tag{4.10}$$

More specifically, we present our algorithm for updating the two spatial degradation operators $\boldsymbol{P}_{1,2}$.

$$\min_{\boldsymbol{P}_1 \in \Phi} \frac{1}{2} \left\| \boldsymbol{P}_1 \boldsymbol{A}_{(1)}^t - \boldsymbol{H}_{(1)} \right\|_F^2,$$

where $\boldsymbol{H}_{(1)} \in \mathbb{R}^{I' \times J'K}$ and $\boldsymbol{A}_{(1)}^t \in \mathbb{R}^{I \times J'K}$ are the mode-1 unfolding matrices of the tensors $\boldsymbol{\mathcal{H}}$ and $\boldsymbol{\mathcal{S}}^{t+1} \times_2 \boldsymbol{P}_2^t$, respectively. The reformulation of the $\boldsymbol{P}_1$-subproblem is

$$\min_{\boldsymbol{b}_1 \in \mathbb{S}_1} f_1(\boldsymbol{b}_1) = \frac{1}{2} \|\boldsymbol{X}_1(\boldsymbol{b}_1)\|_F^2$$

where $\boldsymbol{X}_1(\boldsymbol{b}_1) = \boldsymbol{D}_1 \boldsymbol{F}_1 \mathrm{Diag}(\sqrt{I} \boldsymbol{F}_1 \boldsymbol{b}_1) \boldsymbol{F}_1^* \boldsymbol{A}_{(1)}^t - \boldsymbol{H}_{(1)}$. Then the iterative point $\boldsymbol{P}_1^{t+1}$ is updated by

$$\begin{cases} \boldsymbol{Q}_1^t = \boldsymbol{F}_1^* \boldsymbol{D}_1^\top \boldsymbol{X}_1(\boldsymbol{b}_1^t) \left( \boldsymbol{A}_{(1)}^t \right)^\top \boldsymbol{F}_1, \\ \nabla f_1(\boldsymbol{b}_1^t) = \boldsymbol{F}_1^* \mathrm{diag}(\boldsymbol{Q}_1^t), \\ \boldsymbol{b}_1^{t+1} = P_{\mathbb{S}_I} \left( \boldsymbol{b}_1^t - \tau_1^t \nabla f_1(\boldsymbol{b}_1^t) \right), \\ \boldsymbol{P}_1^{t+1} = \boldsymbol{D}_1 \boldsymbol{F}_1 \mathrm{Diag}(\sqrt{I} \boldsymbol{F}_1 \boldsymbol{b}_1^{t+1}) \boldsymbol{F}_1^*. \end{cases} \tag{4.11}$$

Analogously, $\boldsymbol{P}_2^{t+1}$ is updated by

$$\begin{cases} \boldsymbol{Q}_2^t = \boldsymbol{F}_2^* \boldsymbol{D}_2^\top \boldsymbol{X}_2(\boldsymbol{b}_2^t) \left( \boldsymbol{A}_{(2)}^t \right)^\top \boldsymbol{F}_2, \\ \nabla f_2(\boldsymbol{b}_2^t) = \boldsymbol{F}_2^* \mathrm{diag}(\boldsymbol{Q}_2^t), \\ \boldsymbol{b}_2^{t+1} = P_{\mathbb{S}_J} \left( \boldsymbol{b}_2^t - \tau_2^t \nabla f_2(\boldsymbol{b}_2^t) \right), \\ \boldsymbol{P}_2^{t+1} = \boldsymbol{D}_2 \boldsymbol{F}_2 \mathrm{Diag}(\sqrt{J} \boldsymbol{F}_2 \boldsymbol{b}_2^{t+1}) \boldsymbol{F}_2^*. \end{cases} \tag{4.12}$$

Here, $\boldsymbol{X}_2(\boldsymbol{b}) = \boldsymbol{D}_2 \boldsymbol{F}_2 \mathrm{Diag}(\sqrt{J} \boldsymbol{F}_2 \boldsymbol{b}_2) \boldsymbol{F}_2^* \boldsymbol{A}_{(2)}^t - \boldsymbol{H}_{(2)}$, where $\boldsymbol{A}_{(2)}^t \in \mathbb{R}^{J \times I'K}$ and $\boldsymbol{H}_{(2)} \in \mathbb{R}^{J' \times I'K}$ are the mode-2 unfolding matrices of the tensors $\boldsymbol{\mathcal{S}}^{t+1} \times_1 \boldsymbol{P}_1^{t+1}$ and $\boldsymbol{\mathcal{H}}$, respectively.

**Updating variable $\boldsymbol{P}_3$:** from the optimization problem in (3.10), this subproblem can be further simplified as

$$\min_{\boldsymbol{b}_3^i \in \mathbb{S}_{k_i}} g(\boldsymbol{b}_3^i) := \sum_{i=1}^{K'} \|\boldsymbol{M}^i - \boldsymbol{\mathcal{S}} \times_3 \boldsymbol{D}_3^i \times_3 \boldsymbol{b}_3^i\|_F^2,$$

where $\boldsymbol{M}_i$ denotes the $i$-th spectral band of the MSI, and $\mathbb{S}_{k_i}$ represents the feasible set for the weight vector. This subproblem is likewise addressed using the projected gradient descent algorithm, and $\boldsymbol{P}_3^{t+1}$ is updated as follows,

$$\begin{cases} \boldsymbol{\mathcal{T}} = \boldsymbol{\mathcal{S}} \times_3 \boldsymbol{D}_3^i, \\ \left[ \nabla g((\boldsymbol{b}_3^i)^t) \right]_{k=1}^{k_i} = 2 \langle \boldsymbol{\mathcal{T}} \times_3 (\boldsymbol{b}_3^i)^t - \boldsymbol{M}^i, \boldsymbol{\mathcal{T}}(:,:,k) \rangle, \\ (\boldsymbol{b}_3^i)^{t+1} = P_{\mathbb{S}_{k_i}} \left( (\boldsymbol{b}_3^i)^t - \tau_{3,i}^t \nabla g((\boldsymbol{b}_3^i)^t) \right), \\ \boldsymbol{P}_3^{t+1}(i,:) = (\boldsymbol{D}_3^i (\boldsymbol{b}_3^i)^{t+1})^\top, \quad i = 1, \dots, K'. \end{cases} \tag{4.13}$$

The pseudo-code of our partially linearized ADMM for solving (3.10) is presented in Algorithm 4.1.

**Algorithm 4.1** Partially Linearized ADMM for solving model (3.10).

**Input:** $\boldsymbol{\mathcal{H}}$, $\boldsymbol{\mathcal{M}}$, initial points $(\boldsymbol{\mathcal{S}}^0, \boldsymbol{\mathcal{Z}}^0, \boldsymbol{\Lambda}^0, \boldsymbol{P}_1^0, \boldsymbol{P}_2^0, \boldsymbol{P}_3^0)$, and parameters $\lambda_1$, $\lambda_2$, $\beta$ and $\mu$.

```
repeat
  for t = 1, ..., T do
    Compute S^{t+1} by (4.2).
    Compute Z^{t+1} by (4.4).
    Compute Λ^{t+1} by (4.6).
    Compute P_1^{t+1} by (4.11).
    Compute P_2^{t+1} by (4.12).
    Compute P_3^{t+1} by (4.13).
  end for
until The stopping criteria are satisfied.
```

**Output:** $\boldsymbol{\mathcal{S}}$.

**4.2. Initialization.** The selection of initial points plays a crucial role in determining the convergence behavior of algorithms for nonconvex optimization problems. An effective initialization strategy should incorporate both the inherent characteristics of the specific problem and relevant prior knowledge.

The most challenging aspect is determining the initial point for the HR-HSI. Recently, Alparone *et al.* [4] proposed an efficient nested hypersharpening approach based on multivariate linear regression. This approach first estimates a weight vector by linking the degraded MSI bands to each HSI band. Subsequently, it applies the weighted combination of all high spatial-resolution MSI bands to enhance the spatial resolution of each HSI band. Inspired by this work, we employ this computationally efficient algorithm to initialize $\boldsymbol{\mathcal{S}}^0$ for our subsequent iterative algorithm.

Specifically, the weight vectors are combined into a matrix, and the optimal solution is obtained by solving the following regularized least-squares problem:

$$\boldsymbol{W}^* = \underset{W \in \mathbb{R}^{(K'+1)\times K}}{\operatorname{argmin}} \left\| \boldsymbol{\mathcal{H}} - \boldsymbol{\mathcal{M}}_d^+ \times_3 \boldsymbol{W} \right\|_F^2 + \frac{\gamma}{2} \|\boldsymbol{W}\|_F^2 , \tag{4.14}$$

where the tensor $\boldsymbol{\mathcal{M}}_d^+ \in \mathbb{R}^{I' \times J' \times (K'+1)}$ is constructed by applying bicubic convolution with downsampling to the MSI and concatenating it with an all-ones matrix along the spectral dimension. The corresponding closed-form solution is

$$\boldsymbol{W}^* = \left( \boldsymbol{M}_{(3)}^+ (\boldsymbol{M}_{(3)}^+)^\top + \gamma \boldsymbol{I} \right)^{-1} \boldsymbol{M}_{(3)}^+ \boldsymbol{H}_{(3)}.$$

The initial HR-HSI can then be reconstructed via

$$\boldsymbol{\mathcal{S}}^0 = \boldsymbol{\mathcal{M}}^+ \times_3 \boldsymbol{W}^*,$$

where $\boldsymbol{\mathcal{M}}^+$ is formed by concatenating the MSI with an all-ones matrix along the spectral dimension.

**4.3. Complexity Analysis.** For the observed tensors $\boldsymbol{\mathcal{H}} \in \mathbb{R}^{I' \times J' \times K}$ and $\boldsymbol{\mathcal{M}} \in \mathbb{R}^{I \times J \times K'}$, the complexity of the initialization is $\mathcal{O}(I'J'(K'+1)^2 + I'J'K(K'+1))$. Within the proximal partially linearized ADMM framework, the per-iteration complexity is analyzed as follows:

(i) $\boldsymbol{\mathcal{S}}$-subproblem: the gradient computation incurs a cost of $\mathcal{O}(KI'J(I+J') + KIJ'(I'+J) + IJKK')$, while the TTNN proximal operator requires a cost of $\mathcal{O}(I^2JK + IJK^2)$. More specifically, the TTNN proximal cost is dominated by two steps.

  (a) The unitary transformation incurs a cost of $\mathcal{O}(IJK^2)$, involving the application of the transform $\mathbf{\Phi} \in \mathbb{C}^{K \times K}$ to each of the $I \times J$ mode-3 fibers of the target tensor $\boldsymbol{\mathcal{S}}$.

  (b) The SVD operations require $\mathcal{O}(I^2JK)$, performed on the transformed slices.

(ii) $\boldsymbol{\mathcal{Z}}$-subproblem: the dominant cost is $\mathcal{O}(IJK)$.

(iii) $\boldsymbol{P}_{1,2}$-subproblems: each subproblem exhibits a dominant complexity of $\mathcal{O}(I^3 + I^2J'K)$.

(iv) $\boldsymbol{P}_3$-subproblem: this requires a per-band update cost of $\mathcal{O}(k_i IJK)$.

Notably, while the $\mathcal{O}(IJK^2)$ term scales quadratically with the number of spectral bands $K$, this overhead remains computationally tractable for applications where $K \ll I, J$. For scenarios involving an exceptionally large $K$, the complexity can be further mitigated through several practical strategies: (a) Spectral subspace projection: The HSI can be projected onto a lower-dimensional subspace via techniques such as principal component analysis. This dimensionality reduction decreases the effective $K$ before the application of TTNN regularization; (b) Fast orthogonal transforms: Computational efficiency can be further enhanced by employing fast orthogonal transforms (e.g., fast Fourier transform, empirical wavelet transform, and fast discrete cosine transform). These operators possess inherent symmetry and periodicity, which enable the reduction of mode-3 fiber transformation complexity to $\mathcal{O}(K \log K)$; (c) Randomized approximations: By shifting computations from the full $K$-dimensional space to a smaller $r$-dimensional subspace ($r \ll K$), randomized methods such as randomized SVD circumvent the quadratic complexity bottleneck, reducing the total cost to $\mathcal{O}(IJKr)$.

## 5. Convergence.

We now establish the convergence analysis of the proposed algorithm.

### 5.1. General optimization problem.

For the sequel convergence analysis of Algorithm 4.1, let us cast the reformulated model (3.10) into a general optimization problem as follows:

$$\begin{aligned} \min_{\boldsymbol{x},\boldsymbol{s},\boldsymbol{y}_1,\boldsymbol{y}_2,\boldsymbol{y}_3} \quad & f(\boldsymbol{x}, \boldsymbol{y}_1, \boldsymbol{y}_2) + g(\boldsymbol{x}, \boldsymbol{y}_3) + \psi(\boldsymbol{x}) + \varphi(\boldsymbol{s}) \\ \text{s.t.} \quad & \boldsymbol{x} = \boldsymbol{s},\ \boldsymbol{y}_1 \in \mathcal{B}_1,\ \boldsymbol{y}_2 \in \mathcal{B}_2,\ \boldsymbol{y}_3 \in \mathcal{B}_3, \end{aligned} \tag{5.1}$$

where $\mathcal{B}_1 \subseteq \mathbb{R}^{n_1}$, $\mathcal{B}_2 \subseteq \mathbb{R}^{n_2}$, and $\mathcal{B}_3 \subseteq \mathbb{R}^{n_3}$ denote three bounded convex sets, $f : \mathbb{R}^m \times \mathbb{R}^{n_1} \times \mathbb{R}^{n_2} \to \mathbb{R}$ and $g : \mathbb{R}^m \times \mathbb{R}^{n_3} \to \mathbb{R}$ are continuously differentiable functions, $\psi : \mathbb{R}^m \to \mathbb{R}_\infty$ and $\varphi : \mathbb{R}^m \to \mathbb{R}_\infty$ are proper, closed, and convex functions.

The Lagrangian function of (5.1) is

$$\mathcal{L}(\boldsymbol{x}, \boldsymbol{y}, \boldsymbol{s}, \lambda) := f(\boldsymbol{x}, \boldsymbol{y}_1, \boldsymbol{y}_2) + g(\boldsymbol{x}, \boldsymbol{y}_3) + \psi(\boldsymbol{x}) + \varphi(\boldsymbol{s}) + \langle \boldsymbol{\lambda}, \boldsymbol{x} - \boldsymbol{s} \rangle.$$

The corresponding subdifferential mapping $\partial \mathcal{L}$ is defined as

$$\partial \mathcal{L}(\boldsymbol{x}, \boldsymbol{y}, \boldsymbol{s}, \boldsymbol{\lambda}) := \begin{pmatrix} \partial \psi(\boldsymbol{x}) + \nabla_{\boldsymbol{x}} f(\boldsymbol{x}, \boldsymbol{y}_1, \boldsymbol{y}_2) + \nabla_{\boldsymbol{x}} g(\boldsymbol{x}, \boldsymbol{y}_3) + \boldsymbol{\lambda} \\ \partial \varphi(\boldsymbol{s}) - \boldsymbol{\lambda} \\ \partial \mathcal{I}_{\mathcal{B}_1}(\boldsymbol{y}_1) + \nabla_{\boldsymbol{y}_1} f(\boldsymbol{x}, \boldsymbol{y}_1, \boldsymbol{y}_2) \\ \partial \mathcal{I}_{\mathcal{B}_2}(\boldsymbol{y}_2) + \nabla_{\boldsymbol{y}_2} f(\boldsymbol{x}, \boldsymbol{y}_1, \boldsymbol{y}_2) \\ \partial \mathcal{I}_{\mathcal{B}_3}(\boldsymbol{y}_3) + \nabla_{\boldsymbol{y}_3} g(\boldsymbol{x}, \boldsymbol{y}_3) \\ \boldsymbol{x} - \boldsymbol{s} \end{pmatrix}. \tag{5.2}$$

Definition 5.1. *$(\boldsymbol{x}^*, \boldsymbol{y}^*, \boldsymbol{s}^*, \boldsymbol{\lambda}^*)$ is called an $\epsilon$-stationary point of problem (5.1) if there exists $G \in \partial \mathcal{L}(\boldsymbol{x}^*, \boldsymbol{y}^*, \boldsymbol{s}^*, \boldsymbol{\lambda}^*)$ such that $\|G\|_2 \le \epsilon$.*

**5.2. Algorithmic recursion.** Before outlining the recursive steps, we establish the Moreau smoothed problem that motivates our algorithm.

$$\begin{aligned} \min_{\boldsymbol{x}, \boldsymbol{z}, \boldsymbol{y}_1, \boldsymbol{y}_2, \boldsymbol{y}_3} \quad & f(\boldsymbol{x}, \boldsymbol{y}_1, \boldsymbol{y}_2) + g(\boldsymbol{x}, \boldsymbol{y}_3) + \psi(\boldsymbol{x}) + \varphi_\mu(\boldsymbol{z}) \\ \text{s.t.} \quad & \boldsymbol{x} = \boldsymbol{z}, \ \boldsymbol{y}_1 \in \mathcal{B}_1, \ \boldsymbol{y}_2 \in \mathcal{B}_2, \ \boldsymbol{y}_3 \in \mathcal{B}_3, \end{aligned} \tag{5.3}$$

where $\varphi_\mu$ is

$$\varphi_\mu(\boldsymbol{z}) = \min_{\boldsymbol{s}} \left\{ \varphi(\boldsymbol{s}) + \frac{1}{2\mu} \|\boldsymbol{s} - \boldsymbol{z}\|^2 \right\}.$$

Definition 5.2. *$(\boldsymbol{x}^*, \boldsymbol{y}^*, \boldsymbol{z}^*, \boldsymbol{\lambda}^*)$ is called a stationary point of the Moreau smoothed problem (5.3) if*

$$0 \in \partial \mathcal{L}_\mu(\boldsymbol{x}^*, \boldsymbol{y}^*, \boldsymbol{z}^*, \boldsymbol{\lambda}^*),$$

*where $\mathcal{L}_\mu$ denotes the Lagrangian function associated with (5.3), defined as*

$$\mathcal{L}_\mu(\boldsymbol{x}, \boldsymbol{y}, \boldsymbol{z}, \lambda) := f(\boldsymbol{x}, \boldsymbol{y}_1, \boldsymbol{y}_2) + g(\boldsymbol{x}, \boldsymbol{y}_3) + \psi(\boldsymbol{x}) + \varphi_\mu(\boldsymbol{z}) + \langle \boldsymbol{\lambda}, \boldsymbol{x} - \boldsymbol{z} \rangle. \tag{5.4}$$

*Furthermore, $(\boldsymbol{x}^*, \boldsymbol{y}^*, \boldsymbol{z}^*, \boldsymbol{\lambda}^*)$ is called an $\epsilon$-stationary point of problem (5.3) if there exists $\tilde{G} \in \partial \mathcal{L}_\mu(\boldsymbol{x}^*, \boldsymbol{y}^*, \boldsymbol{z}^*, \boldsymbol{\lambda}^*)$ such that $\|\tilde{G}\|_2 \le \epsilon$.*

In the following lemma, we establish the relationship between a stationary point of the smoothed problem (5.3) and an $\epsilon$-stationary point of the original problem (5.1).

Lemma 5.3. *Suppose $\mu = \mathcal{O}(\epsilon)$ and the multiplier $\boldsymbol{\lambda}^*$ is bounded. Then, $(\boldsymbol{x}^*, \boldsymbol{y}^*, \boldsymbol{s}^*, \boldsymbol{\lambda}^*)$ associated with the stationary point $(\boldsymbol{x}^*, \boldsymbol{y}^*, \boldsymbol{z}^*, \boldsymbol{\lambda}^*)$ of the Moreau smoothed problem (5.3) constitutes an $\epsilon$-stationary point of the original problem (5.1). Here, $(\boldsymbol{z}^*, \boldsymbol{s}^*)$ is the solution to the following subproblem involving the Lagrangian function $\mathcal{L}_\mu$:*

$$(\boldsymbol{z}^*, \boldsymbol{s}^*) = \arg\min_{\boldsymbol{s}, \boldsymbol{z}} \left\{ \varphi(\boldsymbol{s}) + \frac{1}{2\mu} \|\boldsymbol{s} - \boldsymbol{z}\|^2 + \langle \boldsymbol{\lambda}^*, \boldsymbol{x}^* - \boldsymbol{z} \rangle \right\}. \tag{5.5}$$

*Proof.* If $(\boldsymbol{x}^*, \boldsymbol{y}^*, \boldsymbol{z}^*, \boldsymbol{\lambda}^*)$ is a stationary point of the Moreau smoothed problem (5.3), then we have $\boldsymbol{x}^* = \boldsymbol{z}^*$, $\boldsymbol{\lambda}^* = \nabla\varphi_\mu(\boldsymbol{z}^*)$, and

$$0 \in \begin{pmatrix} \partial\psi(\boldsymbol{x}^*) + \nabla_{\boldsymbol{x}} f(\boldsymbol{x}^*, \boldsymbol{y}_1^*, \boldsymbol{y}_2^*) + \nabla_{\boldsymbol{x}} g(\boldsymbol{x}^*, \boldsymbol{y}_3^*) + \boldsymbol{\lambda}^* \\ \partial\mathcal{I}_{\mathcal{B}_1}(\boldsymbol{y}_1^*) + \nabla_{\boldsymbol{y}_1} f(\boldsymbol{x}^*, \boldsymbol{y}_1^*, \boldsymbol{y}_2^*) \\ \partial\mathcal{I}_{\mathcal{B}_2}(\boldsymbol{y}_2^*) + \nabla_{\boldsymbol{y}_2} f(\boldsymbol{x}^*, \boldsymbol{y}_1^*, \boldsymbol{y}_2^*) \\ \partial\mathcal{I}_{\mathcal{B}_3}(\boldsymbol{y}_3^*) + \nabla_{\boldsymbol{y}_3} g(\boldsymbol{x}^*, \boldsymbol{y}_3^*) \end{pmatrix}.$$

It follows from the optimality condition of subproblem (5.5) that

$$\boldsymbol{\lambda}^* = -\frac{1}{\mu}(\boldsymbol{s}^* - \boldsymbol{z}^*) \in \partial\varphi(\boldsymbol{s}^*).$$

Consequently, $0 \in \partial\varphi(\boldsymbol{s}^*) - \boldsymbol{\lambda}^*$ and $\boldsymbol{z}^* - \boldsymbol{s}^* = \mu\boldsymbol{\lambda}^*$. Since $\boldsymbol{x}^* = \boldsymbol{z}^*$, we obtain

$$\boldsymbol{x}^* - \boldsymbol{s}^* = \boldsymbol{z}^* - \boldsymbol{s}^* = \mu\boldsymbol{\lambda}^*.$$

According to Definition 5.1, there exists a vector $G = (\mathbf{0}; \mathbf{0}; \mathbf{0}; \mathbf{0}; \mathbf{0}; \mu\boldsymbol{\lambda}^*)^\top \in \partial\mathcal{L}(\boldsymbol{x}^*, \boldsymbol{y}^*, \boldsymbol{s}^*, \boldsymbol{\lambda}^*)$. Given that $\boldsymbol{\lambda}^*$ is bounded, we can choose a sufficiently small $\mu = k\epsilon$ such that $\|G\|_2 \le \epsilon$. Consequently, the point $(\boldsymbol{x}^*, \boldsymbol{y}^*, \boldsymbol{s}^*, \boldsymbol{\lambda}^*)$ is verified as an $\epsilon$-stationary point of the original problem (5.1). ■

Furthermore, we introduce the following shorthand notation for the partial gradients involved in the algorithm,

$$\begin{aligned} \boldsymbol{g}_{\boldsymbol{x}}^t &= \nabla_{\boldsymbol{x}} f(\boldsymbol{x}^t, \boldsymbol{y}_1^t, \boldsymbol{y}_2^t) + \nabla_{\boldsymbol{x}} g(\boldsymbol{x}^t, \boldsymbol{y}_3^t), & \boldsymbol{g}_{\boldsymbol{y}_3}^t &= \nabla_{\boldsymbol{y}_3} g(\boldsymbol{x}^{t+1}, \boldsymbol{y}_3^t), \\ \boldsymbol{g}_{\boldsymbol{y}_1}^t &= \nabla_{\boldsymbol{y}_1} f(\boldsymbol{x}^{t+1}, \boldsymbol{y}_1^t, \boldsymbol{y}_2^t), & \boldsymbol{g}_{\boldsymbol{y}_2}^t &= \nabla_{\boldsymbol{y}_2} f(\boldsymbol{x}^{t+1}, \boldsymbol{y}_1^{t+1}, \boldsymbol{y}_2^t). \end{aligned}$$

The iterative recursion associated with Algorithm 4.1 is given by

$$(5.6)\quad \begin{cases} \boldsymbol{x}^{t+1} = \operatorname*{argmin}_{\boldsymbol{x}} \Big\{ \psi(\boldsymbol{x}) + \langle \boldsymbol{g}_{\boldsymbol{x}}^t, \boldsymbol{x} - \boldsymbol{x}^t \rangle + \frac{\tau_{\boldsymbol{x}}^t}{2}\|\boldsymbol{x} - \boldsymbol{x}^t\|^2 + \langle \boldsymbol{\lambda}^t, \boldsymbol{x} - \boldsymbol{z}^t \rangle + \frac{\beta}{2}\|\boldsymbol{x} - \boldsymbol{z}^t\|^2 \Big\}, \\ \boldsymbol{z}^{t+1} = \operatorname*{argmin}_{\boldsymbol{z}} \Big\{ \varphi_\mu(\boldsymbol{z}) + \langle \boldsymbol{\lambda}^t, \boldsymbol{x}^{t+1} - \boldsymbol{z} \rangle + \frac{\beta}{2}\|\boldsymbol{x}^{t+1} - \boldsymbol{z}\|^2 \Big\}, \\ \boldsymbol{\lambda}^{t+1} = \boldsymbol{\lambda}^t + \beta(\boldsymbol{x}^{t+1} - \boldsymbol{z}^{t+1}), \\ \boldsymbol{y}_1^{t+1} = P_{\mathcal{B}_1}\left(\boldsymbol{y}_1^t - \dfrac{\boldsymbol{g}_{\boldsymbol{y}_1}^t}{\tau_{\boldsymbol{y}_1}^t}\right), \quad \boldsymbol{y}_2^{t+1} = P_{\mathcal{B}_2}\left(\boldsymbol{y}_2^t - \dfrac{\boldsymbol{g}_{\boldsymbol{y}_2}^t}{\tau_{\boldsymbol{y}_2}^t}\right), \quad \boldsymbol{y}_3^{t+1} = P_{\mathcal{B}_3}\left(\boldsymbol{y}_3^t - \dfrac{\boldsymbol{g}_{\boldsymbol{y}_3}^t}{\tau_{\boldsymbol{y}_3}^t}\right). \end{cases}$$

**5.3. Lipschitz continuity of the partial gradients and gradients.** The following lemma establishes global Lipschitz continuity of the partial gradients, $\nabla_{\boldsymbol{\mathcal{S}}} L_1$ and $\nabla_{\boldsymbol{\mathcal{S}}} L_2$ associated with (4.1), whose properties correspond to those of $f(\cdot)$ and $g(\cdot)$ in the general minimization problem (5.1).

Lemma 5.4 ([49]). *For any fixed $\boldsymbol{P}_1$, $\boldsymbol{P}_2$, and $\boldsymbol{P}_3$, there exist $\ell(\boldsymbol{P}_1, \boldsymbol{P}_2) > 0$ and $\tilde{\ell}(\boldsymbol{P}_3) > 0$ such that*

$$\begin{aligned} &\|\nabla_{\boldsymbol{\mathcal{S}}} L_1(\boldsymbol{\mathcal{S}}_1, \boldsymbol{P}_1, \boldsymbol{P}_2) - \nabla_{\boldsymbol{\mathcal{S}}} L_1(\boldsymbol{\mathcal{S}}_2, \boldsymbol{P}_1, \boldsymbol{P}_2)\|_F \le \ell(\boldsymbol{P}_1, \boldsymbol{P}_2)\,\|\boldsymbol{\mathcal{S}}_1 - \boldsymbol{\mathcal{S}}_2\|_F, \quad \forall\, \boldsymbol{\mathcal{S}}_1, \boldsymbol{\mathcal{S}}_2. \\ &\|\nabla_{\boldsymbol{\mathcal{S}}} L_2(\boldsymbol{\mathcal{S}}_1, \boldsymbol{P}_3) - \nabla_{\boldsymbol{\mathcal{S}}} L_2(\boldsymbol{\mathcal{S}}_2, \boldsymbol{P}_3)\|_F \le \tilde{\ell}(\boldsymbol{P}_3)\|\boldsymbol{\mathcal{S}}_1 - \boldsymbol{\mathcal{S}}_2\|_F, \end{aligned}$$

*where*

$$\ell(\boldsymbol{P}_1,\boldsymbol{P}_2) = \max\left\{1, \left\|(\boldsymbol{P}_1)^\top \boldsymbol{P}_1\right\|_2 \left\|(\boldsymbol{P}_2)^\top \boldsymbol{P}_2\right\|_2\right\},$$
$$\tilde{\ell}(\boldsymbol{P}_3) = \max\left\{1, \left\|(\boldsymbol{P}_3)^\top \boldsymbol{P}_3\right\|_2\right\}.$$

The following lemma establishes an important property regarding the $\boldsymbol{P}$-subproblems.

Lemma 5.5. *$\nabla f$ in (4.7) is globally Lipschitz continuous. Namely, there exists $\ell_f$ such that*

$$\|\nabla f(\boldsymbol{b}_1) - \nabla f(\boldsymbol{b}_2)\| \le \ell_f \|\boldsymbol{b}_1 - \boldsymbol{b}_2\|, \quad \forall\, \boldsymbol{b}_1, \boldsymbol{b}_2 \in \mathbb{S}_n,$$

*where $\ell_f := \|\boldsymbol{F}^* \boldsymbol{D}^\top \boldsymbol{D}\boldsymbol{A}\boldsymbol{A}^\top \boldsymbol{F}\|_2$.*

*Proof.* It follows from the definition of $\nabla f$ in (4.9) that, for any $\boldsymbol{b}_1$, $\boldsymbol{b}_2$ we have

$$\begin{aligned}
\|\nabla f(\boldsymbol{b}_1) - \nabla f(\boldsymbol{b}_2)\| \le & \|\boldsymbol{F}^*\mathrm{diag}(\boldsymbol{Q}_1) - \boldsymbol{F}^*\mathrm{diag}(\boldsymbol{Q}_2)\| \\
= & \|\mathrm{diag}(\boldsymbol{Q}_1 - \boldsymbol{Q}_2)\| \\
= & \|\mathrm{diag}(\boldsymbol{F}^*\boldsymbol{D}^\top \boldsymbol{X}(\boldsymbol{b}_1 - \boldsymbol{b}_2)\boldsymbol{A}^\top \boldsymbol{A}^\top \boldsymbol{F})\| \\
\le & \|\boldsymbol{F}^*\boldsymbol{D}^\top \boldsymbol{D}\boldsymbol{A}\boldsymbol{A}^\top \boldsymbol{F}\|_2 \|\boldsymbol{b}_1 - \boldsymbol{b}_2\| \\
= & \ell_f \|\boldsymbol{b}_1 - \boldsymbol{b}_2\|,
\end{aligned}$$

where the final inequality follows from the operator norm property $\|\boldsymbol{A}\boldsymbol{b}\| \le \|\boldsymbol{A}\|_2\|\boldsymbol{b}\|$, while the concluding equality is a consequence of the norm-preserving property of unitary transformations. This completes the proof. ∎

Lemma 5.6 ([49]). *$\nabla L_1$ and $\nabla L_2$ associated with (4.1) are Lipschitz continuous on bounded subsets of $\mathbb{R}^{I\times J\times K} \times \mathcal{B}_1 \times \mathcal{B}_2$ and $\mathbb{R}^{I\times J\times K} \times \mathcal{B}_3$, respectively.*

**5.4. Convergence analysis.** Throughout this subsection, let $B_1 \times B_2 \times B_3$ and $B_1 \times B_4$ be any bounded subsets of $\mathbb{R}^m \times \mathcal{B}_1 \times \mathcal{B}_2$ and $\mathbb{R}^m \times \mathcal{B}_3$, respectively.

*Assumption* 5.7.
- (i) $\nabla f$ and $\nabla g$ are Lipschitz continuous on $B_1 \times B_2 \times B_3$ and $B_1 \times B_4$, respectively, with constants $M_1 > 0$ and $M_2 > 0$.
- (ii) The partial gradients of $f$ and $g$ are Lipschitz continuous on $\mathbb{R}^m \times \mathcal{B}_1 \times \mathcal{B}_2$ and $\mathbb{R}^m \times \mathcal{B}_3$, respectively.
- (iii) There exist lower and upper bounds for positive Lipschitz constants on $B_1 \times B_2 \times B_3$ and $B_1 \times B_4$; that is, $L_{\boldsymbol{x}}^t \in [L_{\boldsymbol{x}}, L'_{\boldsymbol{x}}]$, $\tilde{L}_{\boldsymbol{x}}^t \in [\tilde{L}_{\boldsymbol{x}}, \tilde{L}'_{\boldsymbol{x}}]$, $L_{\boldsymbol{y}_1}^t \in [L_{\boldsymbol{y}_1}, L'_{\boldsymbol{y}_1}]$, $L_{\boldsymbol{y}_2}^t \in [L_{\boldsymbol{y}_2}, L'_{\boldsymbol{y}_2}]$, and $L_{\boldsymbol{y}_3}^t \in [L_{\boldsymbol{y}_3}, L'_{\boldsymbol{y}_3}]$.

To facilitate the subsequent analysis, we define the augmented Lagrangian function of problem (5.3) as

$$\mathcal{L}_{\beta,\mu}(\boldsymbol{x},\boldsymbol{y},\boldsymbol{z},\boldsymbol{\lambda}) := f(\boldsymbol{x},\boldsymbol{y}_1,\boldsymbol{y}_2) + g(\boldsymbol{x},\boldsymbol{y}_3) + \psi(\boldsymbol{x}) + \varphi_\mu(\boldsymbol{z}) + \langle \boldsymbol{\lambda}, \boldsymbol{x} - \boldsymbol{z}\rangle + \frac{\beta}{2}\|\boldsymbol{x} - \boldsymbol{z}\|^2. \tag{5.7}$$

Now we establish a lower bound for the augmented Lagrangian function $\mathcal{L}_{\beta,\mu}$ in the following lemma. Let $F^*$ denote the optimal value of (3.9). Note that the model (3.9) is equivalent to

$$\begin{aligned}
\min_{\boldsymbol{x},\boldsymbol{y}_1,\boldsymbol{y}_2,\boldsymbol{y}_3} \quad & f(\boldsymbol{x},\boldsymbol{y}_1,\boldsymbol{y}_2) + g(\boldsymbol{x},\boldsymbol{y}_3) + \psi(\boldsymbol{x}) + \varphi(\boldsymbol{x}) \\
\text{s.t.} \quad & \boldsymbol{y}_1 \in \mathcal{B}_1,\ \boldsymbol{y}_2 \in \mathcal{B}_2,\ \boldsymbol{y}_3 \in \mathcal{B}_3,
\end{aligned} \tag{5.8}$$

$F^*$ is also the optimal value of (5.8).

Lemma 5.8. *Suppose that $\beta\mu > \sqrt{2}$, and let $\varphi(\cdot)$ be the indicator function of a closed convex set. Let $\{\boldsymbol{u}^t := (\boldsymbol{x}^t, \boldsymbol{y}^t, \boldsymbol{z}^t, \boldsymbol{\lambda}^t)\}_{t=0}^{\infty}$ be the sequence generated by recursion (5.6), then $\{\mathcal{L}_{\beta,\mu}(\boldsymbol{u}^t)\}_{t=0}^{\infty}$ is uniformly lower bounded by $F^*$.*

*Proof.* By the L-smoothness of $\varphi_\mu$ with constant $1/\mu$ and $\nabla\varphi_\mu(\boldsymbol{z}^t) = \boldsymbol{\lambda}^t$, we get

$$\varphi_\mu(\boldsymbol{x}^t) \le \varphi_\mu(\boldsymbol{z}^t) + \langle \nabla\varphi_\mu(\boldsymbol{z}^t), \boldsymbol{x}^t - \boldsymbol{z}^t\rangle + \frac{1}{2\mu}\|\boldsymbol{x}^t - \boldsymbol{z}^t\|^2.$$

Hence,

$$\begin{aligned}
\mathcal{L}_{\beta,\mu}(\boldsymbol{u}^t) =& f(\boldsymbol{x}^t, \boldsymbol{y}_1^t, \boldsymbol{y}_2^t) + g(\boldsymbol{x}^t, \boldsymbol{y}_3^t) + \psi(\boldsymbol{x}^t) + \varphi_\mu(\boldsymbol{z}^t) + \langle \boldsymbol{\lambda}^t, \boldsymbol{x}^t - \boldsymbol{z}^t\rangle + \frac{\beta}{2}\|\boldsymbol{x}^t - \boldsymbol{z}^t\|^2 \\
\ge& f(\boldsymbol{x}^t, \boldsymbol{y}_1^t, \boldsymbol{y}_2^t) + g(\boldsymbol{x}^t, \boldsymbol{y}_3^t) + \psi(\boldsymbol{x}^t) + \varphi_\mu(\boldsymbol{x}^t) + \left(\frac{\beta}{2} - \frac{1}{2\mu}\right)\|\boldsymbol{x}^t - \boldsymbol{z}^t\|^2 \\
\ge& f(\boldsymbol{x}^*, \boldsymbol{y}_1^*, \boldsymbol{y}_2^*) + g(\boldsymbol{x}^*, \boldsymbol{y}_3^*) + \psi(\boldsymbol{x}^*) + \varphi(\boldsymbol{x}^*) \\
=& F^*.
\end{aligned}$$

The penultimate inequality follows from the condition $\beta\mu > \sqrt{2}$ and the fact that $\varphi(\cdot)$ is the indicator function of a closed convex set, where $F^*$ denotes the optimal value of the minimization problem (5.8). ■

In the following lemma, we establish the nonincreasing property of the augmented Lagrangian function $\mathcal{L}_{\beta,\mu}$. The parameters in recursion (5.6) are set as

$$\begin{aligned}
&\tau_{\boldsymbol{x}}^t := \left(L_{\boldsymbol{x}}^t + \tilde{L}_{\boldsymbol{x}}^t\right)/2, \quad \tau_{\boldsymbol{y}_1}^t := \gamma_{\boldsymbol{y}_1} L_{\boldsymbol{y}_1}^t/2, \quad \tau_{\boldsymbol{y}_2}^t := \gamma_{\boldsymbol{y}_2} L_{\boldsymbol{y}_2}^t/2, \quad \tau_{\boldsymbol{y}_3}^t := \gamma_{\boldsymbol{y}_3} L_{\boldsymbol{y}_3}^t/2, \\
&\gamma_{\boldsymbol{y}_1}, \gamma_{\boldsymbol{y}_2}, \gamma_{\boldsymbol{y}_3} > 1, \quad \beta > \frac{\sqrt{2}}{\mu}.
\end{aligned} \tag{5.9}$$

Lemma 5.9. *Suppose that Assumption 5.7 holds and let $\varphi(\cdot)$ be the indicator function of a closed convex set. Let the parameters in recursion (5.6) satisfy the conditions in (5.9). Then the following statements hold.*

*(i) The sequence $\{\mathcal{L}_{\beta,\mu}(\boldsymbol{u}^t)\}_{t=0}^{\infty}$ is nonincreasing and there exists $\delta > 0$ such that*

$$\mathcal{L}_{\beta,\mu}(\boldsymbol{u}^{t+1}) \le \mathcal{L}_{\beta,\mu}(\boldsymbol{u}^t) - \delta\|\boldsymbol{u}^{t+1} - \boldsymbol{u}^t\|^2. \tag{5.10}$$

*(ii) $\sum_{t=0}^{\infty}\|\boldsymbol{u}^{t+1} - \boldsymbol{u}^t\|^2 < +\infty$, hence*

$$\lim_{t\to+\infty}\|\boldsymbol{u}^{t+1} - \boldsymbol{u}^t\| = 0.$$

*Proof.* It follows from the definition of $\mathcal{L}_{\beta,\mu}$ in (5.7) that

$$\begin{aligned}
&\mathcal{L}_{\beta,\mu}(\boldsymbol{x}^t,\boldsymbol{z}^t,\boldsymbol{\lambda}^t,\boldsymbol{y}^t_{1,2,3})-\mathcal{L}_{\beta,\mu}(\boldsymbol{x}^{t+1},\boldsymbol{z}^t,\boldsymbol{\lambda}^t,\boldsymbol{y}^t_{1,2,3})\\
=&\psi(\boldsymbol{x}^t)-\psi(\boldsymbol{x}^{t+1})+f(\boldsymbol{x}^t,\boldsymbol{y}^t_{1,2})-f(\boldsymbol{x}^{t+1},\boldsymbol{y}^t_{1,2})+g(\boldsymbol{x}^t,\boldsymbol{y}^t_3)-g(\boldsymbol{x}^{t+1},\boldsymbol{y}^t_3)\\
&+\langle\boldsymbol{\lambda}^t,\boldsymbol{x}^t-\boldsymbol{x}^{t+1}\rangle+\frac{\beta}{2}\|\boldsymbol{x}^t-\boldsymbol{z}^t\|^2-\frac{\beta}{2}\|\boldsymbol{x}^{t+1}-\boldsymbol{z}^t\|^2\\
\geq&\langle\partial\psi(\boldsymbol{x}^{t+1})+\boldsymbol{g}^t_{\boldsymbol{x}},\boldsymbol{x}^t-\boldsymbol{x}^{t+1}\rangle+\frac{\beta}{2}\|\boldsymbol{x}^{t+1}-\boldsymbol{x}^t\|^2+\langle\boldsymbol{x}^t-\boldsymbol{x}^{t+1},\boldsymbol{\lambda}^t+\beta(\boldsymbol{x}^{t+1}-\boldsymbol{z}^t)\rangle\\
&-\left(\frac{L^t_{\boldsymbol{x}}}{2}+\frac{\tilde{L}^t_{\boldsymbol{x}}}{2}\right)\|\boldsymbol{x}^t-\boldsymbol{x}^{t+1}\|^2.
\end{aligned}\tag{5.11}$$

By the optimality conditions of the $\boldsymbol{x}$-subproblem in recursion (5.6), we have

$$-\boldsymbol{g}^t_{\boldsymbol{x}}-\tau^t_{\boldsymbol{x}}(\boldsymbol{x}^{t+1}-\boldsymbol{x}^t)-\boldsymbol{\lambda}^t-\beta(\boldsymbol{x}^{t+1}-\boldsymbol{z}^t)\in\partial\psi(\boldsymbol{x}^{t+1}).$$

Substituting $\partial\psi(\boldsymbol{x}^{t+1})$ into (5.11), we obtain

$$\mathcal{L}_{\beta,\mu}(\boldsymbol{x}^{t+1},\boldsymbol{z}^t,\boldsymbol{\lambda}^t,\boldsymbol{y}^t)-\mathcal{L}_{\beta,\mu}(\boldsymbol{x}^t,\boldsymbol{z}^t,\boldsymbol{\lambda}^t,\boldsymbol{y}^t)\leq-\frac{2\tau^t_{\boldsymbol{x}}+\beta-(L^t_{\boldsymbol{x}}+\tilde{L}^t_{\boldsymbol{x}})}{2}\|\boldsymbol{x}^{t+1}-\boldsymbol{x}^t\|^2.\tag{5.12}$$

On the other hand, the $\boldsymbol{z}$-subproblem implies that

$$\begin{aligned}
&\mathcal{L}_{\beta,\mu}(\boldsymbol{x}^{t+1},\boldsymbol{z}^t,\boldsymbol{\lambda}^t,\boldsymbol{y}^t)-\mathcal{L}_{\beta,\mu}(\boldsymbol{x}^{t+1},\boldsymbol{z}^{t+1},\boldsymbol{\lambda}^t,\boldsymbol{y}^t)\\
=&\varphi_\mu(\boldsymbol{z}^t)-\varphi_\mu(\boldsymbol{z}^{t+1})+\langle\boldsymbol{\lambda}^t,\boldsymbol{z}^{t+1}-\boldsymbol{z}^t\rangle+\frac{\beta}{2}(\|\boldsymbol{x}^{t+1}-\boldsymbol{z}^t\|^2-\|\boldsymbol{x}^{t+1}-\boldsymbol{z}^{t+1}\|^2)\\
=&\varphi_\mu(\boldsymbol{z}^t)-\varphi_\mu(\boldsymbol{z}^{t+1})+\frac{\beta}{2}\|\boldsymbol{z}^{t+1}-\boldsymbol{z}^t\|^2+\langle\boldsymbol{z}^{t+1}-\boldsymbol{z}^t,\boldsymbol{\lambda}^{t+1}\rangle\\
\geq&\langle\nabla\varphi_\mu(\boldsymbol{z}^{t+1}),\boldsymbol{z}^t-\boldsymbol{z}^{t+1}\rangle+\frac{\beta}{2}\|\boldsymbol{z}^{t+1}-\boldsymbol{z}^t\|^2+\langle\boldsymbol{z}^{t+1}-\boldsymbol{z}^t,\boldsymbol{\lambda}^{t+1}\rangle\\
=&\frac{\beta}{2}\|\boldsymbol{z}^{t+1}-\boldsymbol{z}^t\|^2,
\end{aligned}$$

where the inequality follows from the convexity of $\varphi_\mu$, and the last equality holds because $\nabla\varphi_\mu(\boldsymbol{z}^{t+1})=\boldsymbol{\lambda}^{t+1}$. Hence,

$$\mathcal{L}_{\beta,\mu}(\boldsymbol{x}^{t+1},\boldsymbol{z}^{t+1},\boldsymbol{\lambda}^t,\boldsymbol{y}^t)-\mathcal{L}_{\beta,\mu}(\boldsymbol{x}^{t+1},\boldsymbol{z}^t,\boldsymbol{\lambda}^t,\boldsymbol{y}^t)\leq-\frac{\beta}{2}\|\boldsymbol{z}^{t+1}-\boldsymbol{z}^t\|^2.\tag{5.13}$$

Similarly, we have

$$\begin{aligned}
&\mathcal{L}_{\beta,\mu}(\boldsymbol{x}^{t+1},\boldsymbol{z}^{t+1},\boldsymbol{\lambda}^{t+1},\boldsymbol{y}^t)-\mathcal{L}_{\beta,\mu}(\boldsymbol{x}^{t+1},\boldsymbol{z}^{t+1},\boldsymbol{\lambda}^t,\boldsymbol{y}^t)\\
&=\langle\boldsymbol{\lambda}^{t+1}-\boldsymbol{\lambda}^t,\boldsymbol{x}^{t+1}-\boldsymbol{z}^{t+1}\rangle=\frac{1}{\beta}\|\boldsymbol{\lambda}^{t+1}-\boldsymbol{\lambda}^t\|^2,
\end{aligned}$$

where the second equality follows from the recursion $\boldsymbol{\lambda}^{t+1}=\boldsymbol{\lambda}^t+\beta(\boldsymbol{x}^{t+1}-\boldsymbol{z}^{t+1})$. Based on Lemma 2.5, we have

$$\|\boldsymbol{\lambda}^{t+1}-\boldsymbol{\lambda}^t\|=\|\nabla\varphi_\mu(\boldsymbol{z}^{t+1})-\nabla\varphi_\mu(\boldsymbol{z}^t)\|\leq\frac{1}{\mu}\|\boldsymbol{z}^{t+1}-\boldsymbol{z}^t\|.\tag{5.14}$$

Thus,

$$\mathcal{L}_{\beta,\mu}(\boldsymbol{x}^{t+1},\boldsymbol{z}^{t+1},\boldsymbol{\lambda}^{t+1},\boldsymbol{y}^t)-\mathcal{L}_{\beta,\mu}(\boldsymbol{x}^{t+1},\boldsymbol{z}^{t+1},\boldsymbol{\lambda}^{t},\boldsymbol{y}^t)\le\frac{1}{\beta\mu^2}\|\boldsymbol{z}^{t+1}-\boldsymbol{z}^t\|^2. \tag{5.15}$$

Furthermore, denoting $\boldsymbol{w}:=(\boldsymbol{x},\boldsymbol{z},\boldsymbol{\lambda})$ and using the Lipshchitz continuous of $\nabla f(\boldsymbol{x},\cdot,\boldsymbol{y}_2)$, we have

$$\begin{aligned}\mathcal{L}_{\beta,\mu}(\boldsymbol{w}^{t+1},\boldsymbol{y}^t)-\mathcal{L}_{\beta,\mu}(\boldsymbol{w}^{t+1},\boldsymbol{y}_1^{t+1},\boldsymbol{y}_{2,3}^t)&=f(\boldsymbol{x}^{t+1},\boldsymbol{y}_1^t,\boldsymbol{y}_2^t)-f(\boldsymbol{x}^{t+1},\boldsymbol{y}_1^{t+1},\boldsymbol{y}_2^t)\\&\ge\langle\boldsymbol{g}_{\boldsymbol{y}_1}^t,\boldsymbol{y}_1^t-\boldsymbol{y}_1^{t+1}\rangle-\frac{L_{\boldsymbol{y}_1}^t}{2}\|\boldsymbol{y}_1^{t+1}-\boldsymbol{y}_1^t\|^2.\end{aligned}$$

Actually, $\boldsymbol{y}_1^{t+1}=P_{\mathcal{B}_1}(\boldsymbol{y}_1^t-\frac{1}{\tau_{\boldsymbol{y}_1}^t}\boldsymbol{g}_{\boldsymbol{y}_1}^t)$. Hence, it follows from the projection inequality $\langle\boldsymbol{y}_1^t-\frac{1}{\tau_{\boldsymbol{y}_1}^t}\boldsymbol{g}_{\boldsymbol{y}_1}^t-\boldsymbol{y}_1^{t+1},\boldsymbol{y}_1^t-\boldsymbol{y}_1^{t+1}\rangle\le 0$ that

$$\mathcal{L}_{\beta,\mu}(\boldsymbol{w}^{t+1},\boldsymbol{y}_1^{t+1},\boldsymbol{y}_{2,3}^t)-\mathcal{L}_{\beta,\mu}(\boldsymbol{w}^{t+1},\boldsymbol{y}^t)\le-\frac{2\tau_{\boldsymbol{y}_1}^t-L_{\boldsymbol{y}_1}^t}{2}\|\boldsymbol{y}_1^{t+1}-\boldsymbol{y}_1^t\|^2. \tag{5.16}$$

Similarly, we have

$$\mathcal{L}_{\beta,\mu}(\boldsymbol{w}^{t+1},\boldsymbol{y}_{1,2}^{t+1},\boldsymbol{y}_3^t)-\mathcal{L}_{\beta,\mu}(\boldsymbol{w}^{t+1},\boldsymbol{y}_1^{t+1},\boldsymbol{y}_{2,3}^t)\le-\frac{2\tau_{\boldsymbol{y}_2}^t-L_{\boldsymbol{y}_2}^t}{2}\|\boldsymbol{y}_2^{t+1}-\boldsymbol{y}_2^t\|^2, \tag{5.17}$$

$$\mathcal{L}_{\beta,\mu}(\boldsymbol{w}^{t+1},\boldsymbol{y}^{t+1})-\mathcal{L}_{\beta,\mu}(\boldsymbol{w}^{t+1},\boldsymbol{y}_{1,2}^{t+1},\boldsymbol{y}_3^t)\le-\frac{2\tau_{\boldsymbol{y}_3}^t-L_{\boldsymbol{y}_3}^t}{2}\|\boldsymbol{y}_3^{t+1}-\boldsymbol{y}_3^t\|^2. \tag{5.18}$$

By summing (5.12)-(5.17), we get

$$\begin{aligned}&\mathcal{L}_{\beta,\mu}(\boldsymbol{u}^{t+1})-\mathcal{L}_{\beta,\mu}(\boldsymbol{u}^t)\\&\le-\frac{2\tau_{\boldsymbol{x}}^t+\beta-(L_{\boldsymbol{x}}^t+\tilde{L}_{\boldsymbol{x}}^t)}{2}\|\boldsymbol{x}^{t+1}-\boldsymbol{x}^t\|^2-\frac{2\tau_{\boldsymbol{y}_1}^t-L_{\boldsymbol{y}_1}^t}{2}\|\boldsymbol{y}_1^{t+1}-\boldsymbol{y}_1^t\|^2\\&\quad-\frac{2\tau_{\boldsymbol{y}_2}^t-L_{\boldsymbol{y}_2}^t}{2}\|\boldsymbol{y}_2^{t+1}-\boldsymbol{y}_2^t\|^2-\frac{2\tau_{\boldsymbol{y}_3}^t-L_{\boldsymbol{y}_3}^t}{2}\|\boldsymbol{y}_3^{t+1}-\boldsymbol{y}_3^t\|^2-\left(\frac{\beta}{2}-\frac{1}{\beta\mu^2}\right)\|\boldsymbol{z}^{t+1}-\boldsymbol{z}^t\|^2.\end{aligned}$$

Based on the definitions in (5.9) and Assumption 5.7 (iii), we define

$$\delta:=\min\left\{\frac{\beta}{2},\frac{(\gamma_{\boldsymbol{y}_1}-1)L_{\boldsymbol{y}_1}}{2},\frac{(\gamma_{\boldsymbol{y}_2}-1)L_{\boldsymbol{y}_2}}{2},\frac{(\gamma_{\boldsymbol{y}_3}-1)L_{\boldsymbol{y}_3}}{2},\frac{\beta}{2}-\frac{1}{\beta\mu^2}\right\}. \tag{5.19}$$

It follows from the conditions in (5.9) that $\delta>0$. Then we derive (5.10), which implies the statement (i). By summing (5.10) over $t=0,1,\ldots$, and applying the lower bound from Lemma 5.8, we have

$$\sum_{t=0}^{\infty}\|\boldsymbol{u}^{t+1}-\boldsymbol{u}^t\|^2\le\frac{1}{\delta}(\mathcal{L}_{\beta,\mu}(\boldsymbol{u}^0)-F^*)<+\infty,$$

which indicates the statement (ii). The proof is completed. ■

We now present a lemma concerning the boundedness of the subgradients of the Lagrangian functions $\mathcal{L}$ and $\mathcal{L}_\mu$ for problems (5.1) and (5.3), respectively. To this end, we first reformulate the $\boldsymbol{z}$-subproblem in recursion (5.6) into the following form,

$$(\boldsymbol{s}^{t+1}, \boldsymbol{z}^{t+1}) := \arg\min_{\boldsymbol{s},\boldsymbol{z}} \left\{ \varphi(\boldsymbol{s}) + \frac{1}{2\mu}\|\boldsymbol{s}-\boldsymbol{z}\|^2 + \langle \boldsymbol{\lambda}^t, \boldsymbol{x}^{t+1} - \boldsymbol{z}\rangle + \frac{\beta}{2}\|\boldsymbol{x}^{t+1} - \boldsymbol{z}\|^2 \right\}. \tag{5.20}$$

Based on the optimality conditions of (5.20) and $\boldsymbol{z}$-subproblem in (5.6), we obtain

$$\boldsymbol{z}^{t+1} - \boldsymbol{s}^{t+1} = \mu\boldsymbol{\lambda}^{t+1}, \quad \boldsymbol{\lambda}^{t+1} \in \partial\varphi(\boldsymbol{s}^{t+1}), \quad \boldsymbol{\lambda}^{t+1} = \nabla\varphi_\mu(\boldsymbol{z}^{t+1}). \tag{5.21}$$

To facilitate the following analysis, we introduce $\boldsymbol{d}^{t+1} := (\boldsymbol{d}_{\boldsymbol{x}}^{t+1}; \boldsymbol{d}_{\boldsymbol{s}}^{t+1}; \boldsymbol{d}_{\lambda}^{t+1}; \boldsymbol{d}_{\boldsymbol{y}_1}^{t+1}; \boldsymbol{d}_{\boldsymbol{y}_2}^{t+1}; \boldsymbol{d}_{\boldsymbol{y}_3}^{t+1})^\top$, and $\tilde{\boldsymbol{d}}^{t+1} := (\boldsymbol{d}_{\boldsymbol{x}}^{t+1}; \tilde{\boldsymbol{d}}_{\boldsymbol{z}}^{t+1}; \tilde{\boldsymbol{d}}_{\lambda}^{t+1}; \boldsymbol{d}_{\boldsymbol{y}_1}^{t+1}; \boldsymbol{d}_{\boldsymbol{y}_2}^{t+1}; \boldsymbol{d}_{\boldsymbol{y}_3}^{t+1})^\top$. We will establish that these two vectors satisfy $\boldsymbol{d}^{t+1} \in \partial\mathcal{L}(\boldsymbol{x}^{t+1}, \boldsymbol{s}^{t+1}, \boldsymbol{y}^{t+1}, \boldsymbol{\lambda}^{t+1})$ and $\tilde{\boldsymbol{d}}^{t+1} \in \partial\mathcal{L}_\mu(\boldsymbol{x}^{t+1}, \boldsymbol{z}^{t+1}, \boldsymbol{y}^{t+1}, \boldsymbol{\lambda}^{t+1})$ with their block components defined as

$$\begin{cases} \boldsymbol{d}_{\boldsymbol{x}}^{t+1} = \boldsymbol{g}_{\boldsymbol{x}}^{t+1} - \boldsymbol{g}_{\boldsymbol{x}}^{t} - \tau_{\boldsymbol{x}}^t(\boldsymbol{x}^{t+1} - \boldsymbol{x}^t) + \beta(\boldsymbol{z}^t - \boldsymbol{z}^{t+1}), \\ \boldsymbol{d}_{\boldsymbol{s}}^{t+1} = \tilde{\boldsymbol{d}}_{\boldsymbol{z}}^{t+1} = 0, \\ \boldsymbol{d}_{\boldsymbol{\lambda}}^{t+1} = \frac{1}{\beta}(\boldsymbol{\lambda}^{t+1} - \boldsymbol{\lambda}^t) + \mu\boldsymbol{\lambda}^{t+1}, \quad \tilde{\boldsymbol{d}}_{\boldsymbol{\lambda}}^{t+1} = \frac{1}{\beta}(\boldsymbol{\lambda}^{t+1} - \boldsymbol{\lambda}^t), \\ \boldsymbol{d}_{\boldsymbol{y}_1}^{t+1} = \nabla_{\boldsymbol{y}_1} f(\boldsymbol{x}^{t+1}, \boldsymbol{y}_1^{t+1}, \boldsymbol{y}_2^{t+1}) - \boldsymbol{g}_{\boldsymbol{x}}^t - \tau_{\boldsymbol{y}_1}^t(\boldsymbol{y}_1^{t+1} - \boldsymbol{y}_1^t), \\ \boldsymbol{d}_{\boldsymbol{y}_2}^{t+1} = \nabla_{\boldsymbol{y}_2} f(\boldsymbol{x}^{t+1}, \boldsymbol{y}_1^{t+1}, \boldsymbol{y}_2^{t+1}) - \boldsymbol{g}_{\boldsymbol{y}_2}^t - \tau_{\boldsymbol{y}_2}^t(\boldsymbol{y}_2^{t+1} - \boldsymbol{y}_2^t), \\ \boldsymbol{d}_{\boldsymbol{y}_3}^{t+1} = \nabla_{\boldsymbol{y}_3} g(\boldsymbol{x}^{t+1}, \boldsymbol{y}_3^{t+1}) - \boldsymbol{g}_{\boldsymbol{y}_3}^t - \tau_{\boldsymbol{y}_3}^t(\boldsymbol{y}_3^{t+1} - \boldsymbol{y}_3^t). \end{cases} \tag{5.22}$$

Lemma 5.10. *Suppose that Assumption 5.7 holds and let $\varphi(\cdot)$ be a coercive function. Let $\{\boldsymbol{v}^t := (\boldsymbol{x}^t, \boldsymbol{s}^t, \boldsymbol{z}^t, \boldsymbol{y}^t, \boldsymbol{\lambda}^t)\}_{t=0}^\infty$ be the sequence generated by recursion (5.6), and $\{\boldsymbol{u}^t\}_{t=0}^\infty$ be defined in Lemma 5.8. Then, the vectors $\boldsymbol{d}^{t+1}$ and $\tilde{\boldsymbol{d}}^{t+1}$ defined in (5.22) satisfy*

$$\boldsymbol{d}^{t+1} \in \partial\mathcal{L}(\boldsymbol{x}^{t+1}, \boldsymbol{s}^{t+1}, \boldsymbol{y}^{t+1}, \boldsymbol{\lambda}^{t+1}), \qquad \tilde{\boldsymbol{d}}^{t+1} \in \partial\mathcal{L}_\mu(\boldsymbol{x}^{t+1}, \boldsymbol{z}^{t+1}, \boldsymbol{y}^{t+1}, \boldsymbol{\lambda}^{t+1}).$$

*Furthermore, there exists a constant $L_\varphi > 0$ such that*

$$\left\|\boldsymbol{d}^{t+1}\right\| \le \Delta\|\boldsymbol{u}^{t+1} - \boldsymbol{u}^t\| + \mu L_\varphi, \quad \left\|\tilde{\boldsymbol{d}}^{t+1}\right\| \le \Delta\|\boldsymbol{u}^{t+1} - \boldsymbol{u}^t\|, \tag{5.23}$$

*where*

$$\Delta := \max\left\{ M + \frac{\gamma L}{2}, \ \frac{1+\mu\beta^2}{\mu\beta}, \ \frac{L'_{\boldsymbol{x}} + \tilde{L}'_{\boldsymbol{x}}}{2} \right\},$$

*with $M = \max\{M_1, M_2\}$, $\gamma = \max\{\gamma_{\boldsymbol{y}_1}, \gamma_{\boldsymbol{y}_2}, \gamma_{\boldsymbol{y}_3}\}$, $L = \max\{L'_{\boldsymbol{y}_1}, L'_{\boldsymbol{y}_2}, L'_{\boldsymbol{y}_3}\}$. Particularly,*

$$\lim_{t\to+\infty} dist(\mathbf{0}, \partial\mathcal{L}_\mu(\boldsymbol{x}^{t+1}, \boldsymbol{y}^{t+1}, \boldsymbol{z}^{t+1}, \boldsymbol{\lambda}^{t+1})) = 0.$$

*Proof.* The optimality of the $x$-subproblem in (5.6) implies

$$-\boldsymbol{g}_{\boldsymbol{x}}^t - \tau_{\boldsymbol{x}}^t(\boldsymbol{x}^{t+1} - \boldsymbol{x}^t) - \boldsymbol{\lambda}^t - \beta(\boldsymbol{x}^{t+1} - \boldsymbol{z}^t) \in \partial\psi(\boldsymbol{x}^{t+1}).$$

By (5.2), it holds that

$$\boldsymbol{d}_{\boldsymbol{x}}^{t+1} = \boldsymbol{g}_{\boldsymbol{x}}^{t+1} - \boldsymbol{g}_{\boldsymbol{x}}^t - \tau_{\boldsymbol{x}}^t(\boldsymbol{x}^{t+1} - \boldsymbol{x}^t) + \beta(\boldsymbol{z}^t - \boldsymbol{z}^{t+1}) \in \partial_{\boldsymbol{x}}\mathcal{L}(\boldsymbol{x}^{t+1}, \boldsymbol{s}^{t+1}, \boldsymbol{y}^{t+1}, \boldsymbol{\lambda}^{t+1}).$$

Then, under Assumption 5.7 (i) and the definition of $\tau_{\boldsymbol{x}}^t$ in (5.9), the norm of $\boldsymbol{d}_{\boldsymbol{x}}^{t+1}$ is bounded as follows

$$\|\boldsymbol{d}_{\boldsymbol{x}}^{t+1}\| \le M\|\boldsymbol{y}^{t+1} - \boldsymbol{y}^t\| + \frac{L'_{\boldsymbol{x}} + \tilde{L}'_{\boldsymbol{x}}}{2}\|\boldsymbol{x}^{t+1} - \boldsymbol{x}^t\| + \beta\|\boldsymbol{z}^{t+1} - \boldsymbol{z}^t\|, \tag{5.24}$$

where $M = \max\{M_1, M_2\}$. By (5.21), we have $0 \in \partial_{\boldsymbol{s}}\mathcal{L}(\boldsymbol{x}^{t+1}, \boldsymbol{s}^{t+1}, \boldsymbol{y}^{t+1}, \boldsymbol{\lambda}^{t+1})$ and $0 = \nabla_{\boldsymbol{z}}\mathcal{L}_\mu(\boldsymbol{x}^{t+1}, \boldsymbol{z}^{t+1}, \boldsymbol{y}^{t+1}, \boldsymbol{\lambda}^{t+1})$. Furthermore, based on the update scheme for $\boldsymbol{\lambda}^{t+1}$ and the definitions of $\boldsymbol{d}_{\boldsymbol{\lambda}}^{t+1}$ and $\tilde{\boldsymbol{d}}_{\boldsymbol{\lambda}}^{t+1}$ in (5.22), we obtain

$$\nabla_{\boldsymbol{\lambda}}\mathcal{L}_\mu(\boldsymbol{x}^{t+1}, \boldsymbol{z}^{t+1}, \boldsymbol{y}^{t+1}, \boldsymbol{\lambda}^{t+1}) = \boldsymbol{x}^{t+1} - \boldsymbol{z}^{t+1} = \frac{1}{\beta}(\boldsymbol{\lambda}^{t+1} - \boldsymbol{\lambda}^t) = \tilde{\boldsymbol{d}}_{\boldsymbol{\lambda}}^{t+1}, \tag{5.25}$$

and

$$\nabla_{\boldsymbol{\lambda}}\mathcal{L}(\boldsymbol{x}^{t+1}, \boldsymbol{s}^{t+1}, \boldsymbol{y}^{t+1}, \boldsymbol{\lambda}^{t+1}) = \boldsymbol{x}^{t+1} - \boldsymbol{s}^{t+1} = \boldsymbol{x}^{t+1} - \boldsymbol{z}^{t+1} + \boldsymbol{z}^{t+1} - \boldsymbol{s}^{t+1}. \tag{5.26}$$

Hence, (5.21) implies $\boldsymbol{d}_{\boldsymbol{\lambda}}^{t+1} = (\boldsymbol{\lambda}^{t+1} - \boldsymbol{\lambda}^t)/\beta + \mu\boldsymbol{\lambda}^{t+1} = \nabla_{\boldsymbol{\lambda}}\mathcal{L}(\boldsymbol{x}^{t+1}, \boldsymbol{s}^{t+1}, \boldsymbol{y}^{t+1}, \boldsymbol{\lambda}^{t+1})$. By (5.14), we have

$$\left\|\tilde{\boldsymbol{d}}_{\boldsymbol{\lambda}}^{t+1}\right\| \le \frac{1}{\mu\beta}\|\boldsymbol{z}^{t+1} - \boldsymbol{z}^t\|.$$

Due to the assumption that $\psi(\cdot)$ is coercive, we have the sequence $\{\boldsymbol{x}^t\}_{t=0}^\infty$ is bounded [7]. Furthermore, owing to the presence of the penalty term $\frac{\beta}{2}\|\boldsymbol{x} - \boldsymbol{z}\|^2$, the sequence $\{\boldsymbol{z}^t\}_{t=0}^\infty$ is also bounded. Hence, there exists $L_\varphi > 0$ such that

$$\|\boldsymbol{\lambda}^{t+1}\| = \|\nabla\varphi_\mu(\boldsymbol{z}^{t+1})\| \le L_\varphi.$$

Combined with the boundedness of $\mathcal{B}_1$, $\mathcal{B}_2$, and $\mathcal{B}_3$, this implies that the sequence $\{\boldsymbol{v}^t\}_{t=0}^\infty$ is bounded. Consequently,

$$\|\boldsymbol{d}_{\boldsymbol{\lambda}}^{t+1}\| \le \frac{1}{\beta}\|\boldsymbol{\lambda}^{t+1} - \boldsymbol{\lambda}^t\| + \mu\|\boldsymbol{\lambda}^{t+1}\| \le \frac{1}{\mu\beta}\|\boldsymbol{z}^{t+1} - \boldsymbol{z}^t\| + \mu L_\varphi,$$

where the last equality is obtained from (5.14). Similarly, it follows from the $\boldsymbol{y}_1$-subproblem that

$$-\boldsymbol{g}_{\boldsymbol{y}_1}^t - \tau_{\boldsymbol{y}_1}^t(\boldsymbol{y}_1^{t+1} - \boldsymbol{y}_1^t) \in \partial\mathcal{I}(\boldsymbol{y}_1^{t+1}).$$

Then

$$\begin{aligned} &\boldsymbol{d}_{\boldsymbol{y}_1}^{t+1} = \nabla_{\boldsymbol{y}_1} f(\boldsymbol{x}^{t+1}, \boldsymbol{y}_1^{t+1}, \boldsymbol{y}_2^{t+1}) - \boldsymbol{g}_{\boldsymbol{y}_1}^t - \tau_{\boldsymbol{y}_1}^t(\boldsymbol{y}_1^{t+1} - \boldsymbol{y}_1^t) \in \partial_{\boldsymbol{y}_1}\mathcal{L}(\boldsymbol{x}^{t+1}, \boldsymbol{s}^{t+1}, \boldsymbol{y}^{t+1}, \boldsymbol{\lambda}^{t+1}), \\ &\|\boldsymbol{d}_{\boldsymbol{y}_1}^{t+1}\| \le M_1\|\boldsymbol{y}_1^{t+1} - \boldsymbol{y}_1^t\| + M_1\|\boldsymbol{y}_2^{t+1} - \boldsymbol{y}_2^t\| + \frac{\gamma_{\boldsymbol{y}_1} L'_{\boldsymbol{y}_1}}{2}\|\boldsymbol{y}_1^{t+1} - \boldsymbol{y}_1^t\|, \end{aligned} \tag{5.27}$$

where the second inequality is based on Assumption (5.7) (ii) and (5.9). Analogously, we get

$$\begin{aligned} &\boldsymbol{d}_{\boldsymbol{y}_2}^{t+1} = \nabla_{\boldsymbol{y}_2} f(\boldsymbol{x}^{t+1}, \boldsymbol{y}_1^{t+1}, \boldsymbol{y}_2^{t+1}) - \boldsymbol{g}_{\boldsymbol{y}_2}^t - \tau_{\boldsymbol{y}_2}^t(\boldsymbol{y}_2^{t+1} - \boldsymbol{y}_2^t) \in \partial_{\boldsymbol{y}_2}\mathcal{L}(\boldsymbol{x}^{t+1}, \boldsymbol{s}^{t+1}, \boldsymbol{y}^{t+1}, \boldsymbol{\lambda}^{t+1}), \\ &\|\boldsymbol{d}_{\boldsymbol{y}_2}^{t+1}\| \le M_1\|\boldsymbol{y}_2^{t+1} - \boldsymbol{y}_2^t\| + \frac{\gamma_{\boldsymbol{y}_2} L'_{\boldsymbol{y}_2}}{2}\|\boldsymbol{y}_2^{t+1} - \boldsymbol{y}_2^t\|, \end{aligned} \tag{5.28}$$

and

$$\begin{aligned} &\boldsymbol{d}_{\boldsymbol{y}_3}^{t+1} = \nabla_{\boldsymbol{y}_3} g(\boldsymbol{x}^{t+1}, \boldsymbol{y}_3^{t+1}) - \boldsymbol{g}_{\boldsymbol{y}_3}^t - \tau_{\boldsymbol{y}_3}^t(\boldsymbol{y}_3^{t+1} - \boldsymbol{y}_3^t) \in \partial_{\boldsymbol{y}_3}\mathcal{L}(\boldsymbol{x}^{t+1}, \boldsymbol{s}^{t+1}, \boldsymbol{y}^{t+1}, \boldsymbol{\lambda}^{t+1}), \\ &\|\boldsymbol{d}_{\boldsymbol{y}_3}^{t+1}\| \le M_2\|\boldsymbol{y}_3^{t+1} - \boldsymbol{y}_3^t\| + \frac{\gamma_{\boldsymbol{y}_3} L'_{\boldsymbol{y}_3}}{2}\|\boldsymbol{y}_3^{t+1} - \boldsymbol{y}_3^t\|. \end{aligned} \tag{5.29}$$

Combining all equations in (5.24)-(5.29), and the definitions of $\mathcal{L}$ and $\mathcal{L}_\mu$ defined in (5.2) and (5.4), we get

$$\begin{aligned} \boldsymbol{d}^{t+1} &= (\boldsymbol{d}_{\boldsymbol{x}}^{t+1}; \boldsymbol{d}_{\boldsymbol{s}}^{t+1}; \boldsymbol{d}_{\boldsymbol{\lambda}}^{t+1}; \boldsymbol{d}_{\boldsymbol{y}_1}^{t+1}; \boldsymbol{d}_{\boldsymbol{y}_2}^{t+1}; \boldsymbol{d}_{\boldsymbol{y}_3}^{t+1})^\top \in \partial\mathcal{L}(\boldsymbol{x}^{t+1}, \boldsymbol{s}^{t+1}, \boldsymbol{y}^{t+1}, \boldsymbol{\lambda}^{t+1}), \\ \tilde{\boldsymbol{d}}^{t+1} &= (\boldsymbol{d}_{\boldsymbol{x}}^{t+1}; \tilde{\boldsymbol{d}}_{\boldsymbol{z}}^{t+1}; \tilde{\boldsymbol{d}}_{\boldsymbol{\lambda}}^{t+1}; \boldsymbol{d}_{\boldsymbol{y}_1}^{t+1}; \boldsymbol{d}_{\boldsymbol{y}_2}^{t+1}; \boldsymbol{d}_{\boldsymbol{y}_3}^{t+1})^\top \in \partial\mathcal{L}_\mu(\boldsymbol{x}^{t+1}, \boldsymbol{z}^{t+1}, \boldsymbol{y}^{t+1}, \boldsymbol{\lambda}^{t+1}). \end{aligned}$$

By the triangle inequality, we obtain

$$\begin{aligned} \left\|\tilde{\boldsymbol{d}}^{t+1}\right\| &\le (M + \frac{\gamma L}{2})\|\boldsymbol{y}^{t+1} - \boldsymbol{y}^t\| + \frac{(1+\mu\beta^2)}{\mu\beta}\|\boldsymbol{z}^{t+1} - \boldsymbol{z}^t\| + \frac{L'_{\boldsymbol{x}} + \tilde{L}'_{\boldsymbol{x}}}{2}\|\boldsymbol{x}^{t+1} - \boldsymbol{x}^t\|, \\ \|\boldsymbol{d}^{t+1}\| &\le (M + \frac{\gamma L}{2})\|\boldsymbol{y}^{t+1} - \boldsymbol{y}^t\| + \frac{(1+\mu\beta^2)}{\mu\beta}\|\boldsymbol{z}^{t+1} - \boldsymbol{z}^t\| + \frac{L'_{\boldsymbol{x}} + \tilde{L}'_{\boldsymbol{x}}}{2}\|\boldsymbol{x}^{t+1} - \boldsymbol{x}^t\| + \mu L_\varphi, \end{aligned}$$

where $\gamma = \max\{\gamma_{\boldsymbol{y}_1}, \gamma_{\boldsymbol{y}_2}, \gamma_{\boldsymbol{y}_3}\}$, and $L = \max\{L'_{\boldsymbol{y}_1}, L'_{\boldsymbol{y}_2}, L'_{\boldsymbol{y}_3}\}$. Denote

$$\Delta := \max\left\{M + \frac{\gamma L}{2},\ \frac{1+\mu\beta^2}{\mu\beta},\ \frac{L'_{\boldsymbol{x}} + \tilde{L}'_{\boldsymbol{x}}}{2}\right\}, \tag{5.30}$$

then the inequality (5.23) holds. Since the sequence $\{\boldsymbol{v}^t\}_{t=0}^\infty$ is bounded, it follows that $\{\boldsymbol{u}^t\}_{t=0}^\infty$ is also bounded. Consequently, there exists a subsequence $\{\boldsymbol{u}^{t_l}\}_{l=1}^\infty$ such that $\lim_{l\to\infty} \boldsymbol{u}^{t_l} = \boldsymbol{u}^*$. Hence, $\lim_{l\to+\infty}\|\boldsymbol{u}^{t_l+1} - \boldsymbol{u}^{t_l}\| = 0$. Set $t \ge t_l$, it follows from all inequalities in (5.24)-(5.29), we have $\lim_{l\to+\infty}\|\tilde{\boldsymbol{d}}^{t_l+1}\| = 0$, i.e., $\lim_{k\to+\infty}\operatorname{dist}(\mathbf{0}, \mathcal{L}_\mu(\boldsymbol{x}^{t+1}, \boldsymbol{y}^{t+1}, \boldsymbol{z}^{t+1}, \boldsymbol{\lambda}^{t+1})) = 0$, which completes the proof. ∎

Lemma 5.11. *The Lagrangian function $\mathcal{L}_\mu(\boldsymbol{x}, \boldsymbol{y}, \boldsymbol{z}, \lambda)$ defined in* (5.4) *is a KŁ function.*

*Proof.* As stated in [5], the indicator function on the non-negative set is semi-algebraic. Besides, according to [40], it has been established that the transformed tensor nuclear norm $\|\cdot\|_{\text{TTNN}}$ is also semi-algebraic. Therefore, the function $\mathcal{L}_\mu$ defined in (5.4) is semi-algebraic, thereby satisfying the KŁ property. The proof is completed. ∎

Following a similar approach to the proof of Theorem 1 in [6], we establish the global convergence of the considered sequence $\{\boldsymbol{u}^t\}_{t=0}^{\infty}$ under KŁ assumption of $\mathcal{L}_\mu$. To this end, we summarize the required conditions in the following assumption.

*Assumption* 5.12. The functions in problem (5.1) associated with model (3.10) satisfy:
(i) $\varphi : \mathbb{R}^m \to \mathbb{R}_\infty$ is an indicator function of a closed convex set.
(ii) $\psi : \mathbb{R}^m \to \mathbb{R}_\infty$ is a coercive function.

Theorem 5.13. *Suppose that Assumptions 5.7 and 5.12 hold. Let $\{\boldsymbol{u}^t\}_{t=0}^{\infty}$ be the sequence generated by recursion (5.6) with the parameters satisfying (5.9), then any limit point of the sequence is a critical point of Lagrangian function $\mathcal{L}_\mu$. Moreover, $\{\boldsymbol{u}^t\}_{t=0}^{\infty}$ converges globally to the critical point $\boldsymbol{u}^*$ of problem (5.3).*

*Remark* 5.14. As established in the proof of Lemma 5.10, the sequence $\{\boldsymbol{v}^t\}_{t=0}^{\infty}$ is bounded. Hence, Assumption 5.7 is naturally satisfied. However, for the sake of completeness and clarity, we retain Assumption 5.7 in Theorem 5.13, explicitly including the associated constants.

Next, we establish the iteration complexity of the proposed partially linearized ADMM for achieving an $\epsilon$-stationary point for both the Moreau smoothed problem (5.3) and the original problem (5.1).

Theorem 5.15. *Suppose that Assumptions 5.7 and 5.12 hold. Let $\{\boldsymbol{v}^t\}_{t=1}^{T}$ be the sequence generated by recursion (5.6) with the parameters satisfying (5.9). For any tolerance $\epsilon > 0$, there exists $\tilde{G}^t \in \partial\mathcal{L}_\mu(\boldsymbol{x}^t, \boldsymbol{y}^t, \boldsymbol{z}^t, \boldsymbol{\lambda}^t)$, for $t = 1, \ldots, T$ such that $\min_{t=1,\ldots,T} \|\tilde{G}^t\|^2 \leq \epsilon^2$, provided that $T = \mathcal{O}(\epsilon^{-2})$.*

*Proof.* According to Lemma 5.10, there exists $\tilde{G}^t \in \partial\mathcal{L}_\mu(\boldsymbol{x}^t, \boldsymbol{y}^t, \boldsymbol{z}^t, \boldsymbol{\lambda}^t)$ such that $\|\tilde{G}^t\|^2 \leq \Delta^2\|\boldsymbol{u}^t - \boldsymbol{u}^{t-1}\|^2$. Since the parameters $\delta$ and $\Delta$ (defined in (5.19) and (5.30), respectively) are constants, the descent condition in Lemma 5.9 implies that

$$\left\|\tilde{G}^t\right\|^2 \leq \Delta^2\|\boldsymbol{u}^t - \boldsymbol{u}^{t-1}\|^2 \leq \frac{\Delta^2}{\delta}\left(\mathcal{L}_{\beta,\mu}(\boldsymbol{u}^{t-1}) - \mathcal{L}_{\beta,\mu}(\boldsymbol{u}^t)\right).$$

Now, by summing the inequality over $t = 1, \ldots, T$, it follows from Lemma 5.8 that

$$\frac{1}{T}\sum_{t=1}^{T}\left\|\tilde{G}^t\right\|^2 \leq \frac{\Delta^2}{T\delta}\left(\mathcal{L}_{\beta,\mu}(\boldsymbol{u}^0) - F^*\right).$$

Consequently, to ensure that $\min_{t=1,\ldots,T} \|\tilde{G}^t\| \leq \epsilon^2$, the number of required iterations is $T = \mathcal{O}(\epsilon^{-2})$. This completes the proof. ■

Theorem 5.16. *Suppose that Assumptions 5.7 and 5.12 hold. Let $\{\boldsymbol{v}^t\}_{t=1}^{T}$ be the sequence generated by recursion (5.6) with the parameters satisfying (5.9). For any tolerance $\epsilon > 0$, the parameters $\mu = \mathcal{O}(\epsilon)$ and $\beta = 1/\epsilon$ are chosen to ensure $\beta\mu > \sqrt{2}$. Then, there exists $G^t \in \partial\mathcal{L}(\boldsymbol{x}^t, \boldsymbol{s}^t, \boldsymbol{y}^t, \boldsymbol{\lambda}^t)$, for $t = 1, \ldots, T$ such that $\min_{t=1,\ldots,T} \|G^t\|^2 \leq \epsilon^2$, provided that $T = \mathcal{O}(\epsilon^{-4})$.*

*Proof.* According to Lemma 5.10, there exists $G^t \in \partial\mathcal{L}(\boldsymbol{x}^t, \boldsymbol{s}^t, \boldsymbol{y}^t, \boldsymbol{\lambda}^t)$ such that $\|G^t\|^2 \leq 2\Delta^2\|\boldsymbol{u}^t - \boldsymbol{u}^{t-1}\|^2 + 2\mu^2 L_\varphi^2$. By setting $\beta = \frac{1}{\epsilon}$ and $\mu = \mathcal{O}(\epsilon)$ such that $\beta\mu > \sqrt{2}$, the parameter $\delta$ defined in (5.19) becomes a constant independent of $\beta$ and $\mu$. Moreover, $\Delta$ as defined in (5.30) is given by $\Delta = (1+\mu\beta^2)/\mu\beta$. This, together with the decrease condition in Lemma 5.9, implies that

$$\begin{aligned}\|G^t\|^2 &\le 2\Delta^2\|\boldsymbol{u}^t - \boldsymbol{u}^{t-1}\|^2 + 2\mu^2 L_\varphi^2\\ &= \frac{2(1+\mu\beta^2)^2}{\mu^2\beta^2}\|\boldsymbol{u}^t - \boldsymbol{u}^{t-1}\|^2 + 2\mu^2 L_\varphi^2\\ &\le \frac{2(1+\mu\beta^2)^2}{\delta\mu^2\beta^2}\left(\mathcal{L}_{\beta,\mu}(\boldsymbol{u}^{t-1}) - \mathcal{L}_{\beta,\mu}(\boldsymbol{u}^t)\right) + 2\mu^2 L_\varphi^2.\end{aligned}$$

Now, by summing the inequality over $t = 1, \ldots, T$, it follows from Lemma 5.8 that

$$\frac{1}{T}\sum_{t=1}^{T}\|G^t\|^2 \le \frac{2(1+\mu\beta^2)^2}{T\delta\mu^2\beta^2}\left(\mathcal{L}_{\beta,\mu}(\boldsymbol{u}^0) - F^*\right) + 2\mu^2 L_\varphi^2.$$

Consequently, to ensure that $\min_{t=1,\ldots,T}\|G^t\| \le \epsilon^2$, the number of required iterations is $T = \mathcal{O}(\epsilon^{-4})$. This completes the proof. ■

*Remark* 5.17. It is worth noting that the iteration complexity of $\mathcal{O}(\epsilon^{-4})$ is consistent with existing results [28, 53] for solving nonsmooth, nonconvex optimization problems by ADMM via smoothing frameworks. Since the smoothing parameter $\mu$ must be coupled with the optimality tolerance $\epsilon$ (specifically $\mu = \mathcal{O}(\epsilon)$), the resulting complexity naturally increases from the standard $\mathcal{O}(\epsilon^{-2})$ in Theorem 5.15 (typical of smooth nonconvex problem) to $\mathcal{O}(\epsilon^{-4})$ in Theorem 5.16.

**6. Numerical Experiments.** In this section, we first introduce the evaluation datasets and quality assessment metrics, followed by extensive numerical simulations to validate the proposed model and algorithm. Furthermore, we compare the fused results across different methods through both visual qualitative and quantitative analyses, demonstrating the competitive performance of our framework. This framework is hereafter referred to as Tenfuse.

**6.1. Datasets.** In this section, we evaluate and compare the fusion performance of representative methods across four distinct datasets, including two synthetic datasets (Washington DC Mall[1] and Chikusei[2]) and two real-world datasets (Hyperion[3]/Sentinel-2[4] and EnMAP[5]/Sentinel-2).

**6.1.1. Washington DC mall.** The first dataset is the Washington DC Mall, captured by the HYDICE sensor. The original HSI has a spatial resolution of 2.8 m, dimensions of $1028 \times 307$ pixels, and 210 bands. After removing bands affected by atmospheric absorption, 191 bands remain. A $300 \times 300 \times 191$ subset serves as the ground-truth HR-HSI, as illustrated in Figure 6.1 (a). The LR-HSI is generated by applying a $9 \times 9$ Gaussian filter ($\sigma = 1$) and downsampling by 4. The MSI is synthesized using the Landsat 8 SRF [6]. Additionally, zero-mean Gaussian noise is added to achieve SNRs of 30 dB for the HSI and 35 dB for the MSI. The final input pair is constructed based on a $9 \times 9$ Gaussian PSF and the Landsat 8 SRF.

[1] https://engineering.purdue.edu/~biehl/MultiSpec/hyperspectral.html
[2] https://naotoyokoya.com/Download.html
[3] https://earthexplorer.usgs.gov/
[4] https://dataspace.copernicus.eu/
[5] https://www.enmap.org/data_access/

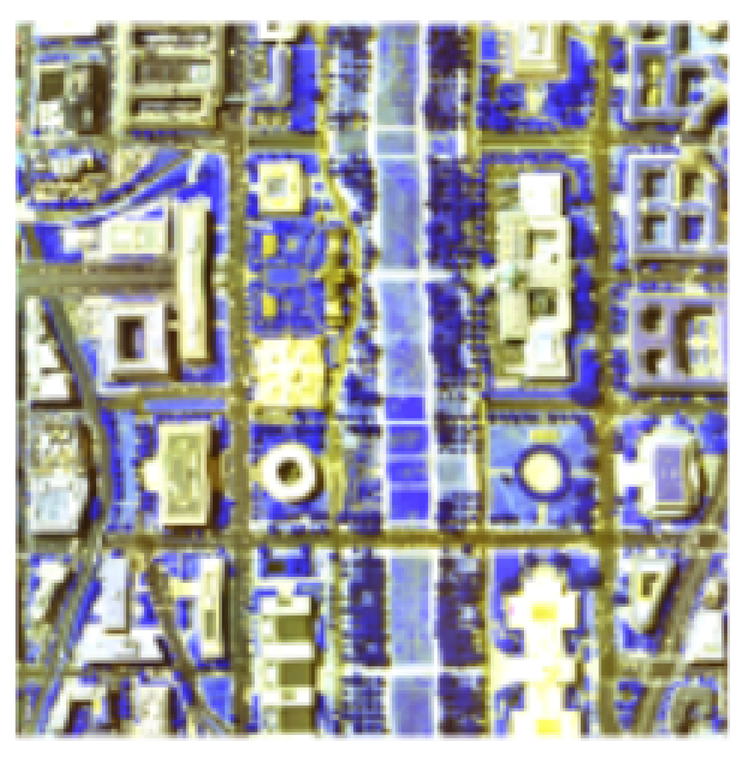

(a) Washington DC mall

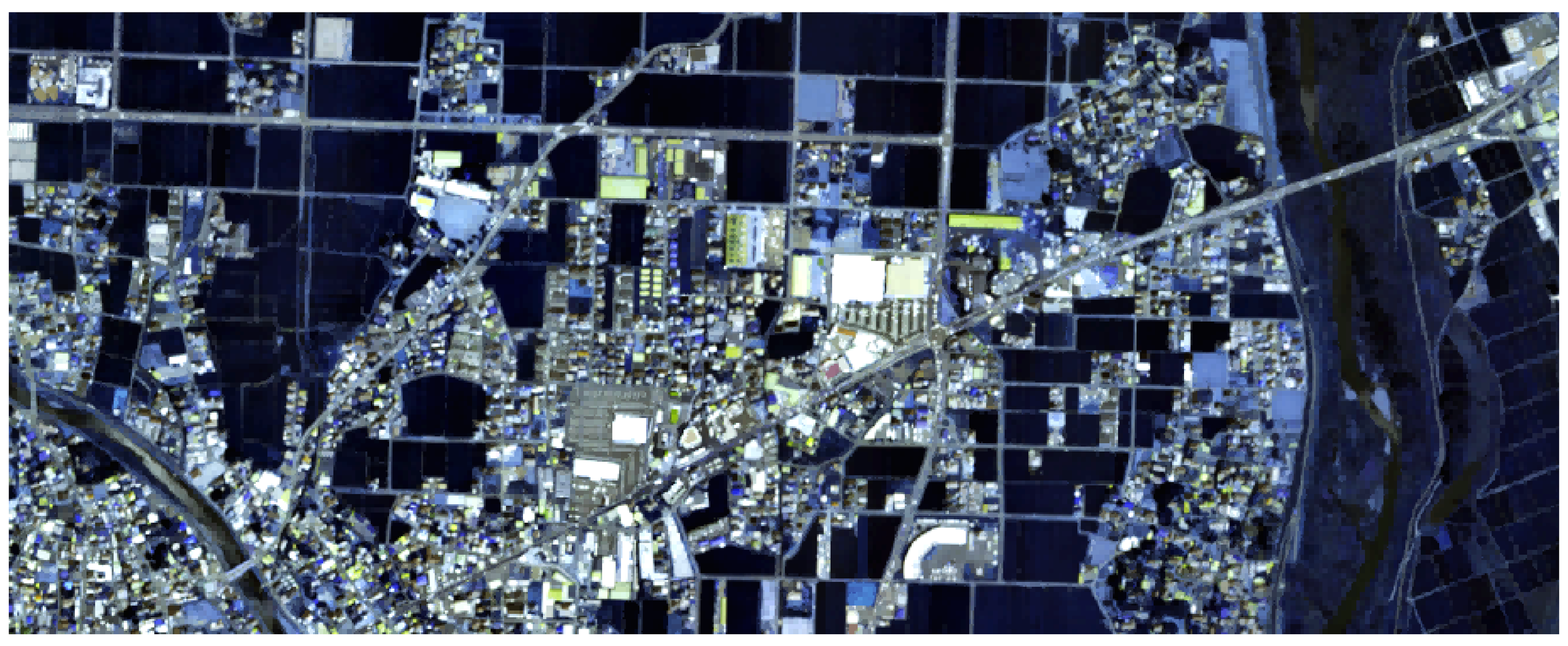

(b) Chikusei

**Figure 6.1.** *Two ground-truth HR-HSIs from synthetic datasets.*

**6.1.2. Chikusei.** The second benchmark dataset is the Chikusei hyperspectral image, acquired over Chikusei, Japan, using a Headwall Hyperspec-VNIR-C sensor. This dataset was publicly released by Naoto Yokoya and Akira Iwasaki at the University of Tokyo. The original hyperspectral dataset comprises 128 spectral bands covering the wavelength range from 360–1020 nm and the scene consists of $2517 \times 2335$ pixels with a spatial resolution of 2.5 m. We extract a $400 \times 1000$ pixels as the ground-truth HR-HSI, as illustrated in Figure 6.1 (b). While the procedure for generating the input LR-HSI and HR-MSI pair mirrors that of the Washington DC Mall dataset, the SRF construction diverges by utilizing the visible-to-near-infrared bands from the Landsat 8 sensor.

**6.1.3. Hyperion/Sentinel-2.** Sentinel-2 is an Earth Observation (EO) satellite launched by the European Space Agency (ESA). The MSI sensor onboard the Sentinel-2 satellite acquires 13 bands from 442 nm to 2185 nm. For the fusion experiments conducted in this study, only bands with the spatial resolution of 10 meters (bands 2, 3, 4, and 8) for Sentinel-2 HR-MSI are selected.

Hyperion is an imaging spectrometer aboard NASA's Earth Observing-1 (EO-1) satellite. The Hyperion sensor acquires hyperspectral imagery (HSI) across 242 contiguous spectral bands covering 400–2500 nm, with a nominal bandwidth of 10 nm per band and a spatial resolution of 30 meters. The sensor employs two separate detector arrays to cover different spectral regions: visible-near infrared (VNIR, 356–1085 nm) and short-wave infrared (SWIR, 852–2577 nm). The VNIR and SWIR spectral ranges exhibit a 20-channel overlap region. After retaining VNIR bands 8–57 and SWIR channels 79–224 while eliminating redundant overlapping bands, the final dataset comprises 196 unique spectral bands.

The captured scene is located in the Freiburg region, Baden-Württemberg, Germany. It contains various landforms such as urban areas, rural villages, forests, rivers, and farmland. The two images were acquired on September 24, 2016 (for Hyperion) and August 23, 2016 (for Sentinel-2). The image registration was performed by aligning the geographic coordinates (47.76-47.79°N, 7.49-7.55°E) between the Hyperion LR-HSI and Sentinel-2 HR-MSI. Subsets with the following dimensions were then obtained for this experiment: $100 \times 160 \times 196$ (Hyperion) and $300 \times 480 \times 4$ (Sentinel-2), as shown in Figure 6.2. This area represents an actual ground coverage of 3 km × 4.8 km.

[6] https://nwp-saf.eumetsat.int/site/

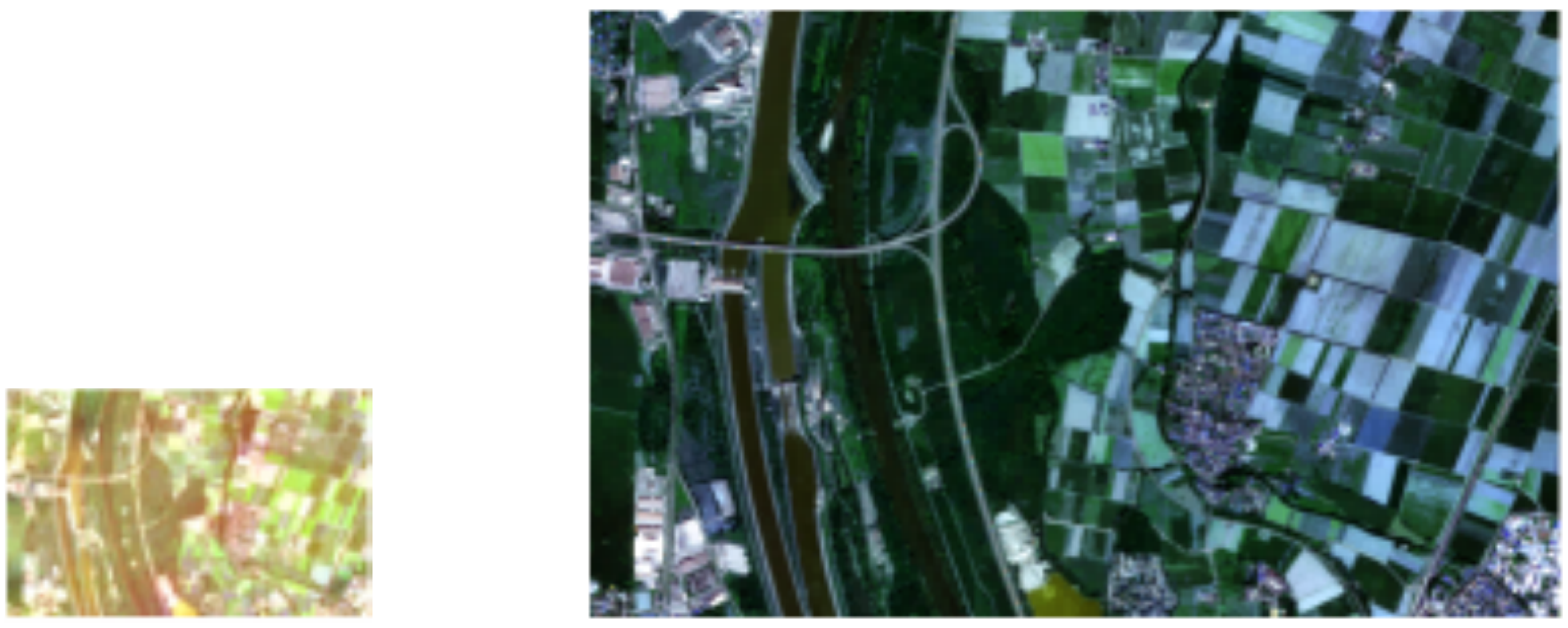

(a) Hyperion LR-HSI (b) Sentinel-2 HR-MSI

**Figure 6.2.** *Hyperion/Sentinel-2 dataset covering a real-world geographic region of 3 km × 4.8 km.*

**6.1.4. EnMAP/Sentinel-2.** EnMAP, developed by DLR and Kayser-Threde, provides hyperspectral imagery with 224 spectral bands (91 VNIR and 133 SWIR) at a 30 m spatial resolution, covering the 420–2450 nm range. The study area covers the Nuremberg region in Germany, featuring a diverse landscape of urban, rural, forest, and farmland areas. Both EnMAP and Sentinel-2A captured the region within a one-hour interval on April 30, 2025 (10:00–11:00 UTC). Following geographic alignment (49.18–49.24°N, 10.07–10.28°E), co-registered subsets with dimensions of $200 \times 500 \times 224$ (EnMAP) and $600 \times 1500 \times 4$ (Sentinel-2) were extracted. The subsets, as illustrated in Figure 6.3, cover a ground area of 6 km × 15 km for the fusion experiments.

## 6.2. Quality assessment.

**6.2.1. Full-reference quality assessment for synthetic datasets.** In blind HSI-MSI fusion, we treat the available HSI as reference, from which synthetic HR-MSI and LR-HSI are generated via SRF and PSF (with downsampling), respectively. The fused result is then quantitatively compared against this reference using standard metrics [11]: Peak Signal-to-Noise Ratio (PSNR), Spectral Angle Mapper (SAM), Universal Image Quality Index (UIQI), and Erreur Relative Globale Adimensionnelle de Synthèse (ERGAS).

For the fused image $\boldsymbol{X} \in \mathbb{R}^{m \times n}$ and its reference image $\boldsymbol{X}^* \in \mathbb{R}^{m \times n}$,

$$\mathrm{PSNR}(\boldsymbol{X}, \boldsymbol{X}^*) = 10 \log \frac{255^2}{\mathrm{MSE}(\boldsymbol{X}, \boldsymbol{X}^*)},$$

where the Mean Squared Error (MSE) measures the $\ell_2$ error as follows,

$$\mathrm{MSE}(\boldsymbol{X}, \boldsymbol{X}^*) = \frac{\|\boldsymbol{X} - \boldsymbol{X}^*\|_F^2}{mn}.$$

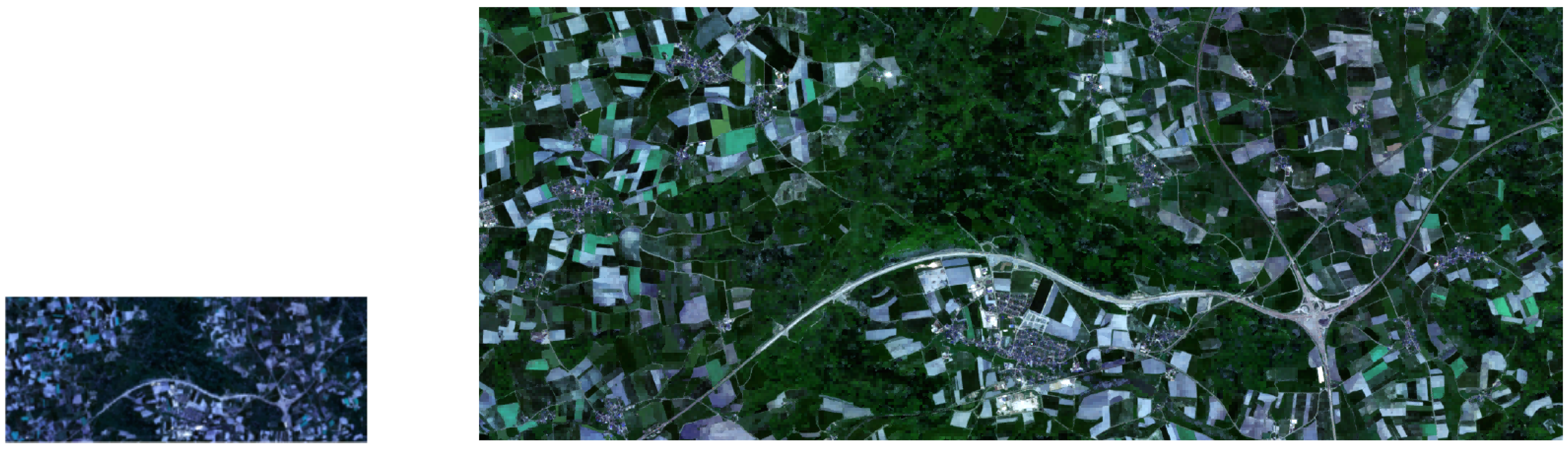

(a) EnMAP LR-HSI (b) Sentinel-2 HR-MSI

**Figure 6.3.** *EnMAP/Sentinel-2 dataset covering a real-world geographic region of 6 km × 15 km.*

The PSNR for HSI is defined as the average value of all bands, and a higher value of PSNR indicates better fusion quality.

$$\mathrm{SAM}(\boldsymbol{X}, \boldsymbol{X}^*) = \frac{1}{n}\sum_{j=1}^{n} \arccos\left(\frac{\langle \boldsymbol{x}_j, \boldsymbol{x}_j^* \rangle}{\|\boldsymbol{x}_j\|\|\boldsymbol{x}_j^*\|}\right).$$

The SAM value for HSI is calculated by averaging the pixel-wise SAM values across the entire image. The UIQI between two image patches is defined as

$$\mathrm{UIQI}(\boldsymbol{x}, \boldsymbol{x}^*) = \left(\frac{2\mu_{\boldsymbol{x}}\mu_{\boldsymbol{x}^*}}{\mu_{\boldsymbol{x}}^2 + \mu_{\boldsymbol{x}^*}^2}\right)\left(\frac{2\sigma_{\boldsymbol{x}}\sigma_{\boldsymbol{x}^*}}{\sigma_{\boldsymbol{x}}^2 + \sigma_{\boldsymbol{x}^*}^2}\right)\left(\frac{\sigma_{\boldsymbol{x},\boldsymbol{x}^*}}{\sigma_{\boldsymbol{x}}\sigma_{\boldsymbol{x}^*}}\right),$$

where $\sigma$ and $\mu$ represent the variance and mean, respectively. $\sigma_{\boldsymbol{x},\boldsymbol{x}^*}$ is the covariance between $\boldsymbol{x}$ and $\boldsymbol{x}^*$, capturing their structural similarity. The average value of all spectral bands is the UIQI for HSI, and a larger value of UIQI means the better fusion quality.

$$\mathrm{ERGAS}(\boldsymbol{X}, \boldsymbol{X}^*) = d\sqrt{\frac{1}{K}\sum_{k=1}^{K}\left(\frac{\mathrm{RMSE}_k}{\mu_k}\right)^2},$$

where $d$ denotes the spatial resolution ratio between HSI and MSI. For the $k$-th spectral band, $\mathrm{RMSE}_k = \sqrt{\mathrm{MSE}_k}$ and $\mu_k$ is the mean radiance value of the $k$-th band in the reference image. The optimal value of ERGAS is 0.

**6.2.2. No-reference quality assessment for real-world datasets.** In real-world HSI-MSI fusion tasks, acquiring a corresponding full-resolution reference image for quantitative evaluation is often impractical or impossible. When a PAN image is available, we adopt three no-reference metrics from [3]: the Spectral Distortion Index ($D_\lambda$), the Spatial Distortion Index ($D_s$), and the Quality with No Reference (QNR) index. The QNR index is defined as $\mathrm{QNR} = (1 - \boldsymbol{D}_\lambda)(1 - \boldsymbol{D}_s)$, where the spectral distortion index $\boldsymbol{D}_\lambda$ quantifies spectral consistency, while the spatial distortion index $\boldsymbol{D}_s$ measures spatial fidelity via multivariate regression between the fused HR-HSI bands and the available PAN image. More specifically,

$$\boldsymbol{D}_\lambda = \sqrt{\frac{1}{K(K-1)} \sum_{i=1}^{K} \sum_{\substack{j=1 \\ j\neq i}}^{K} |Q(\boldsymbol{H}_i, \boldsymbol{H}_j) - Q(\boldsymbol{S}_i, \boldsymbol{S}_j)|},$$

where $\boldsymbol{H}_i$ and $\boldsymbol{S}_i$ denote the $i$-th bands of the LR-HSI and fused HR-HSI, respectively, while $Q$ represents the UIQI index from full-reference quality assessment.

$$\boldsymbol{D}_s = \sqrt{\frac{1}{K} \sum_{k=1}^{K} \left| Q(\boldsymbol{S}_k, \boldsymbol{P}) - Q(\boldsymbol{H}_k, \tilde{\boldsymbol{P}}) \right|},$$

where $\boldsymbol{H}_k$ is the $k$-th band of LR-HSI, while its counterpart $\tilde{\boldsymbol{P}}$ is generated by applying bicubic convolution to the panchromatic image $\boldsymbol{P}$.

**6.3. Test approaches.** We now compare the numerical performance of our proposed method, Tenfuse, with nine state-of-the-art approaches involving three categories, including Hypersharpening-based methods: AWLP [41], SM-SaBC [32], and Nested-GSA [4]; Learning-based methods: UDALN [29], SURE [34], and ZSL [12]; Model-based methods: CNMF [52], Hysure [39], and IR-TenSR [48]. For the implementation of our model (3.10), we set $\gamma = 0.1$ for the initialization of the HR-HSI. Subsequently, we choose $\sigma_s = 1$ for the spatial initialization of $\boldsymbol{P}_{1,2}$ and $\sigma_\lambda = 100$ for the spectral initialization of $\boldsymbol{P}_3$. In model (3.10), the regularization parameters $\boldsymbol{\lambda}_1$ and $\boldsymbol{\lambda}_2$ are tuned from $\{0.1, 1, 10, 50\}$ and $\{0.01, 0.1, 1\}$, respectively, and we report the best recovery results obtained. To facilitate the iterative updates of the Algorithm 4.1, we set $\alpha^t = 1/\max\{1, \tilde{\alpha}^t\}$, and $\beta \in \{1, 10, 20, 100\}$ for $\boldsymbol{\mathcal{S}}$-subproblem solver, where

$$\tilde{\alpha}^t := \left\|(\boldsymbol{P}_1^t)^\top \boldsymbol{P}_1^t\right\|_2 \left\|(\boldsymbol{P}_2^t)^\top \boldsymbol{P}_2^t\right\|_2 + \lambda_1 \left\|(\boldsymbol{P}_3^t)^\top \boldsymbol{P}_3^t\right\|_2.$$

The parameter $\mu$ for the $\boldsymbol{\mathcal{Z}}$-subproblem solver is selected from the set $\{20^{-1}, 20^{-2}, 20^{-3}, 20^{-4}\}$. The step sizes are defined as $\tau_1^t = 1/\left\|\nabla^2 f(\boldsymbol{b}_1^t)\right\|_2$ for the $\boldsymbol{P}_1$-subproblem, $\tau_2^t = 1/\left\|\nabla^2 f(\boldsymbol{b}_2^t)\right\|_2$ for the $\boldsymbol{P}_2$-subproblem, and $\tau_{3,i}^t = 1/\|\nabla^2 g(\boldsymbol{b}_3^i)^t\|_2$ for the $\boldsymbol{P}_3$-subproblem. The results reported below correspond to the optimal performance achieved by each method.

**6.4. Performance on synthetic datasets.** The quantitative results in Table 6.1 demonstrate that Tenfuse, a model-based method, consistently outperforms the competing approaches across both synthetic datasets (Washington DC Mall and Chikusei) in key metrics such as PSNR, SAM, UIQI, and ERGAS. It achieves the highest reconstruction accuracy, while maintaining a reasonable computational time. More notably, learning-based methods such as UDALN and SURE achieve comparable fusion quality in metrics like PSNR and SAM, but suffer from prohibitively long runtimes. These advantages make Tenfuse a preferred choice for practical applications where both precision and speed are critical. The visual results in Figures 6.4 and 6.5 provide further evidence for the effectiveness of Tenfuse.

**6.5. Performance on real-world datasets.**

**Table 6.1**

*Performance of ten test methods on two synthetic datasets. For Washing DC Mall, MSI:* $300 \times 300 \times 8$ *and HSI:* $75 \times 75 \times 191$*; For Chikusei, MSI:* $400 \times 1000 \times 4$ *and HSI:* $100 \times 250 \times 128$.

| Datasets | Metrics | **Hypersharpening-based methods** | | | **Learning-based methods** | | | **Model-based methods** | | | |
|---|---|---|---|---|---|---|---|---|---|---|---|
| | | AWLP | SM-SaBC | Nested-GSA | UDALN | SURE | ZSL | CNMF | Hysure | IR-TenSR | Tenfuse |
| Washington DC Mall | PSNR (∞) | 35.55 | 33.07 | 46.52 | 49.04 | 48.00 | 46.12 | 11.49 | 45.58 | 31.39 | **49.33** |
| | SAM (0) | 0.41 | 0.10 | 0.05 | 0.05 | 0.07 | 0.06 | 0.33 | 0.11 | 0.67 | **0.04** |
| | UIQI (1) | 0.26 | 0.38 | 0.86 | 0.88 | 0.67 | 0.77 | 0.02 | 0.86 | 0.29 | **0.90** |
| | ERGAS (0) | 3.82 | 1.59 | 1.80 | 2.85 | 7.39 | 6.09 | 65.54 | 2.85 | 5.61 | **1.77** |
| | Time (s) | **0.56** | 1.83 | 1.00 | 3216.36 | 5140.25 | 4342.73 | 5.13 | 68.23 | 122.70 | 67.08 |
| Chikusei | PSNR (∞) | 31.61 | 33.08 | 41.35 | 44.29 | 44.20 | 43.27 | 7.95 | 40.53 | 23.10 | **44.54** |
| | SAM (0) | 0.15 | 0.10 | 0.06 | 0.05 | 0.05 | 0.07 | 0.27 | 0.08 | 0.13 | **0.03** |
| | UIQI (1) | 0.37 | 0.38 | 0.87 | 0.89 | 0.88 | 0.77 | 0.06 | 0.82 | 0.44 | **0.92** |
| | ERGAS (0) | 1.78 | 1.5 | 0.72 | 0.58 | 0.56 | 0.65 | 32.05 | 0.77 | 4.58 | **0.46** |
| | Time (s) | **1.29** | 5.62 | 2.65 | 6034.28 | 22323.20 | 14026.88 | 1.46 | 97.32 | 117.49 | 183.96 |

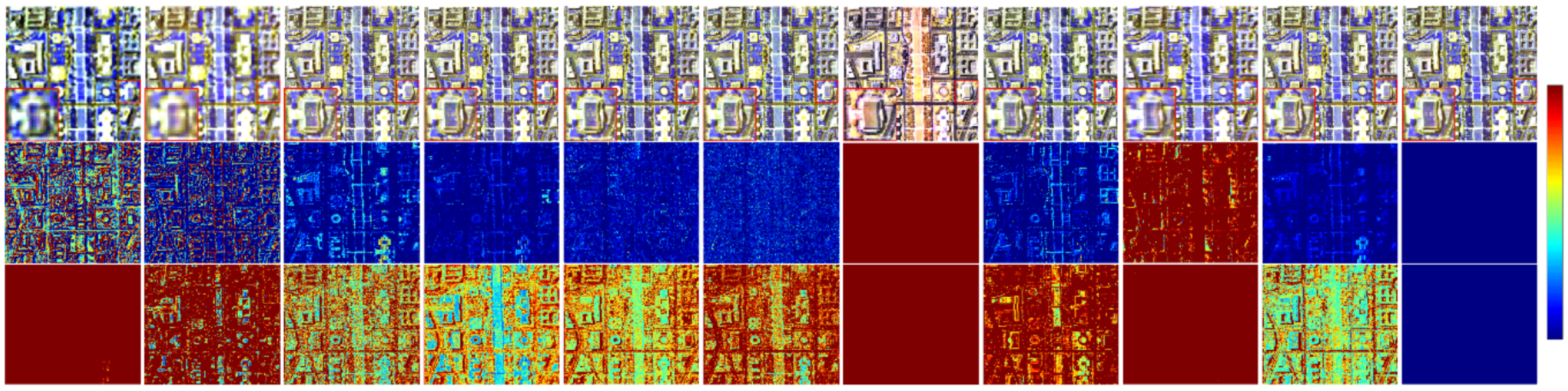

**Figure 6.4.** *Fusion results of synthetic Washing DC mall dataset. Row 1: Pseudo-color representation of reconstructed HR-HSI (R:14/B:66/G:24). Row 2: Residual heatmap at band 53. Row 3: Spectral angle mapper (SAM) spatial distribution.*

**6.5.1. Hyperion/Sentinel-2 fusion.** The QNR metric (ranging from 0 to 1, where higher is better) is employed to assess performance, utilizing the panchromatic band from NASA's ALI for the fusion of real Hyperion/Sentinel-2 data, illustrated in Figure 6.7 (e). As shown in Table 6.2, SM-SaBC exhibits the lowest spectral distortion, while Tenfuse achieves the smallest spatial distortion. Notably, our proposed Tenfuse method attains the highest overall QNR score, demonstrating its balance between spectral fidelity and spatial detail preservation.

Figure 6.6 displays pseudo-color images of all tested methods, generated from the fused HR-HSI using bands 8, 15, and 23 (498–650 nm, visible range) for the tested methods. A representative region is magnified to facilitate visual comparison. The magnified view reveals that Nested-GSA, UDALN, SURE, and Tenfuse yield sharper spatial details compared to other methods. This observation further validates Tenfuse's effectiveness in preserving fine image structures. Furthermore, as illustrated in Figure 6.7, Tenfuse exhibits performance comparable to Nested-GSA and UDALN and superior to SURE in the SWIR spectral range.

**6.5.2. Full-scale fusion of EnMAP/Sentinel-2.** In full-scale fusion of EnMAP/Sentinel-2, no panchromatic image is available. Consequently, only the spectral distortion $D_\lambda$ metric can be derived from the reconstructed HR-HSI and observed LR-HSI. To complement the

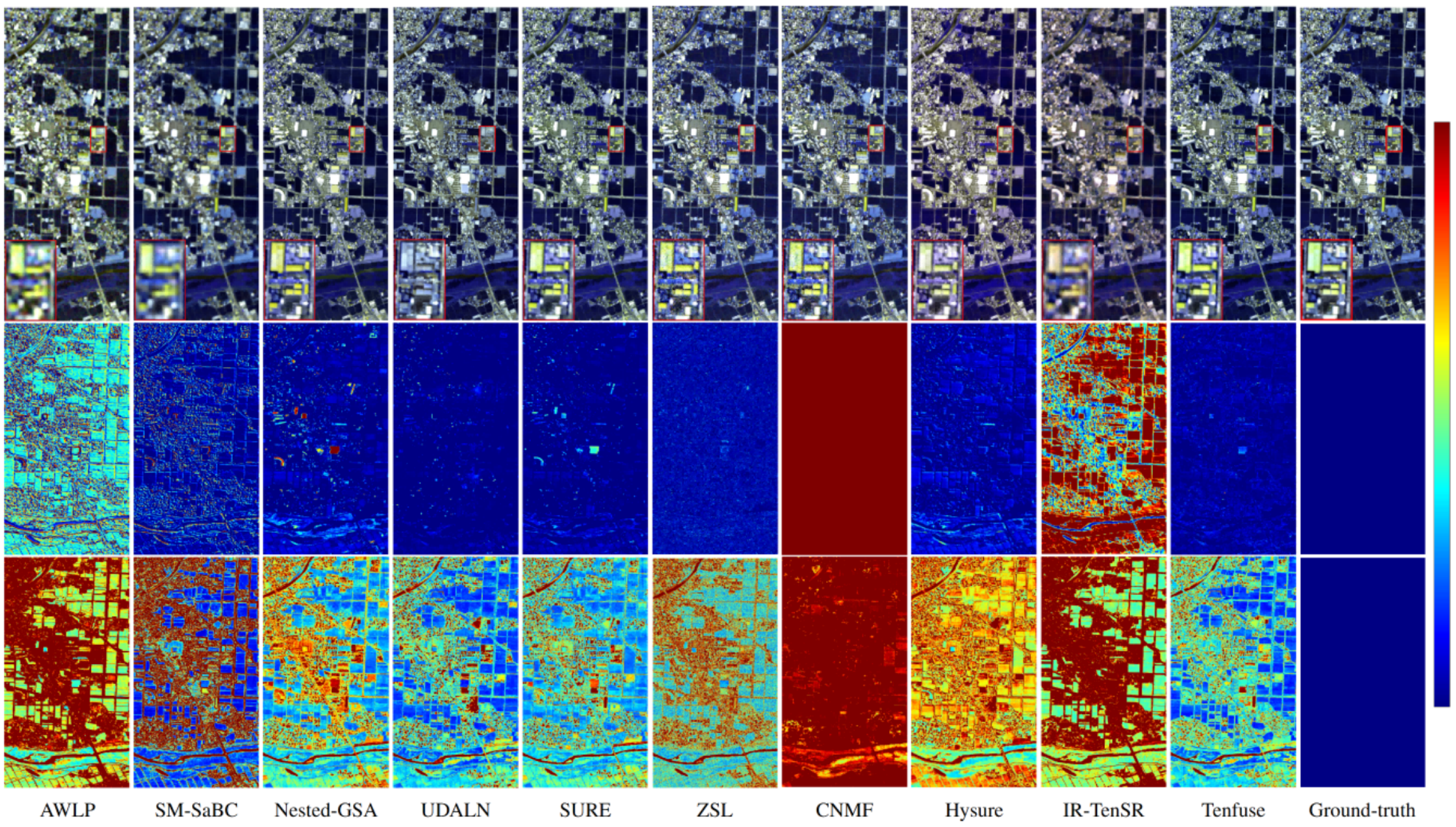

**Figure 6.5.** *Fusion results of synthetic Chikusei dataset. Row 1: Pseudo-color representation of reconstructed HR-HSI (R:17, G:24, B:60). Row 2: Residual heatmap at band 60. Row 3: Spectral angle mapper (SAM) error distribution.*

**Table 6.2**

*Performance of real-world dataset Hyperion/Sentinel-2. MSI:* $300 \times 480 \times 4$ *and HSI:* $100 \times 160 \times 196$.

| Metrics | Hypersharpening-based methods | | | Learning-based methods | | | Model-based methods | | | |
|---|---|---|---|---|---|---|---|---|---|---|
| | **AWLP** | **SM-SaBC** | **Nested-GSA** | **UDALN** | **SURE** | **ZSL** | **CNMF** | **Hysure** | **IR-TenSR** | **Tenfuse** |
| $D_\lambda$ (0) | 0.0278 | **0.0213** | 0.0448 | 0.0462 | 0.0619 | 0.0306 | 0.1104 | 0.0464 | 0.0223 | 0.0511 |
| $D_s$ (0) | 0.0908 | 0.0895 | 0.0475 | 0.0431 | 0.0278 | 0.0674 | 0.2559 | 0.0708 | 0.1032 | **0.0222** |
| QNR (1) | 0.8839 | 0.8911 | 0.9098 | 0.9129 | 0.9120 | 0.9041 | 0.6620 | 0.8861 | 0.8768 | **0.9278** |
| Time (s) | **0.85** | 2.56 | 1.79 | 4947.64 | 20124.88 | 11364.03 | 0.4736 | 68.23 | 39.96 | 30.66 |

spatial consistency analysis, we employed the coefficient of determination ($R^2$ in [1]) from the multivariate linear regression of the HR-HSI onto the combined bands of the HR-MSI.

Table 6.3 presents the $D_\lambda$ and $R^2$ results of various comparative methods for real data fusion. As shown in the table, SM-SaBC demonstrates the lowest spectral distortion, while Tenfuse achieves the highest spatial consistency. Figure 6.8 displays pseudo-color images generated from the fused HR-HSI using bands 15, 30, and 48 (487–653 nm, visible range) for the tested methods. A representative region is magnified to facilitate visual comparison. The magnified view reveals that Nested-GSA, UDALN, and Tenfuse yield sharper spatial details compared to other methods. Furthermore, as demonstrated in Figure 6.9, Tenfuse achieves performance comparable to Nested-GSA and outperforms SURE in the SWIR spectral range, confirming its consistent effectiveness across multiple spectral ranges.

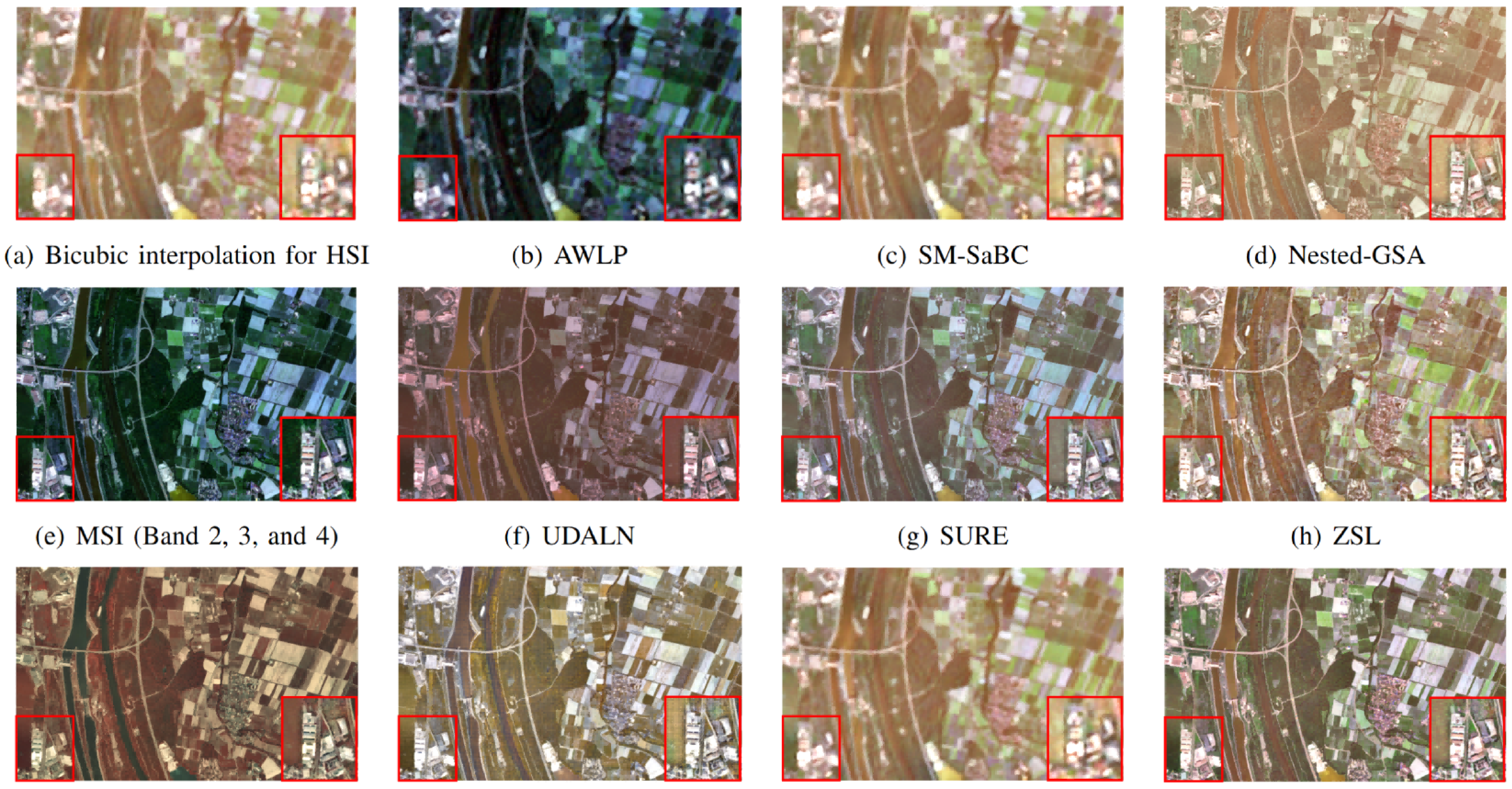

(a) Bicubic interpolation for HSI (b) AWLP (c) SM-SaBC (d) Nested-GSA
(e) MSI (Band 2, 3, and 4) (f) UDALN (g) SURE (h) ZSL
(i) CNMF (j) Hysure (k) IR-TenSR (l) Tenfuse

**Figure 6.6.** *Fusion results of Hyperion/Sentinel-2 in VNIR spectral range. Pseudo-color visualization of reconstructed HR-HSIs using band composition (R:8/B:23/G:15) corresponding to wavelengths (498, 650, and 569 nm).*

**Table 6.3**

*Performance on real-world dataset EnMAP/Sentinel-2. MSI:* $600 \times 1500 \times 4$ *and HSI:* $200 \times 500 \times 224$.

| Metrics | Hypersharpening-based methods | | | Learning-based methods | | | Model-based methods | | | |
|---|---|---|---|---|---|---|---|---|---|---|
| | **AWLP** | **SM-SaBC** | **Nested-GSA** | **UDALN** | **SURE** | **ZSL** | **CNMF** | **Hysure** | **IR-TenSR** | **Tenfuse** |
| $D_\lambda$ (0) | 0.0283 | **0.0246** | 0.1857 | 0.1810 | 0.2836 | 0.0853 | 0.2569 | 0.2718 | 0.0318 | 0.1249 |
| $R^2$ (1) | 0.0236 | 0.0259 | 0.6849 | 0.5728 | 0.6360 | 0.0552 | 0.5345 | 0.5574 | 0.0344 | **0.6952** |
| Time (s) | **5.87** | 18.54 | 9.54 | 39668.05 | 83904.88 | 23224.09 | 11.96 | 37.52 | 298.62 | 241.89 |

**6.6. Robustness analysis.** This subsection presents a sensitivity analysis to evaluate the impact of key parameters on the performance of the Tenfuse model. The analysis focuses on four primary aspects: (a) the trade-off parameter $\lambda_1$; (b) the spectral mixing window size $k_i$; (c) the Moreau envelope parameter $\mu$; and (d) model's robustness against mismatched, band-dependent noise. By systematically adjusting these variables, we assess the robustness of the model and identify the optimal configurations for both spectral and spatial reconstruction.

**6.6.1. Balance between spatial and spectral quality.** The balance between spatial and spectral quality is pivotal for image fusion and is closely governed by the regularization parameters. To investigate this relationship, we conducted an ablation study using the real-world Hyperion/Sentinel-2 dataset by systematically varying the weight of the spectral data-fidelity term, $\lambda_1$, and the spectral mixing window size, $k_i$.

As illustrated in Figure 6.10, increasing $\lambda_1$ significantly reduces spectral distortion $D_\lambda$, allowing Tenfuse to achieve performance comparable to or even exceeding that of SM-SaBC in spectral fidelity. However, this improvement often entails a marginal sacrifice in spatial sharpness $D_s$. This confirms the flexibility of the proposed model, which can be fine-tuned to prioritize either spectral integrity or spatial detail depending on the application requirements.

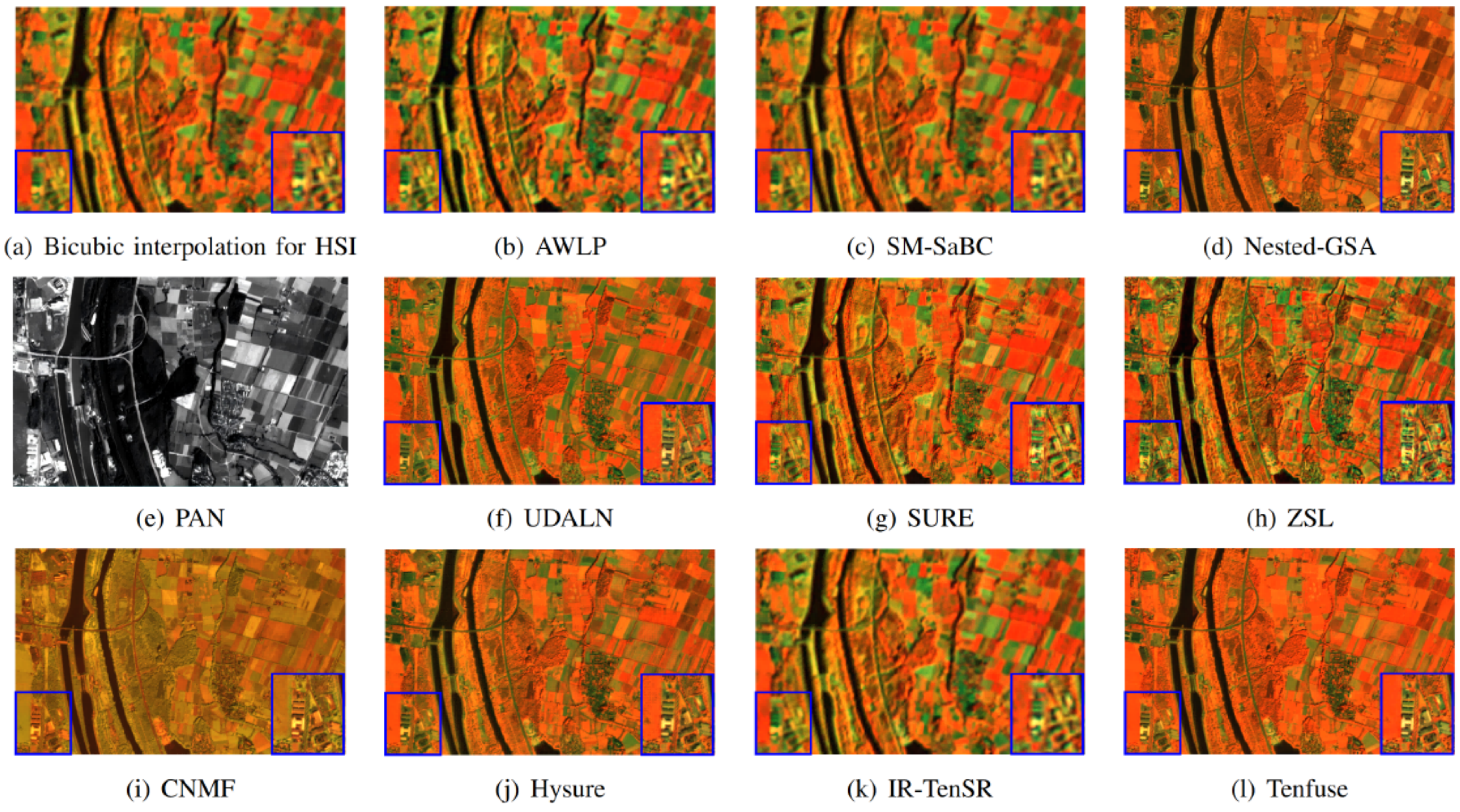

**Figure 6.7.** *Fusion results of Hyperion/Sentinel-2 in SWIR spectral range. Pseudo-color visualization of reconstructed HR-HSIs using band composition (R:46/B:177/G:119) corresponding to wavelengths (864, 2203, 1618 nm).*

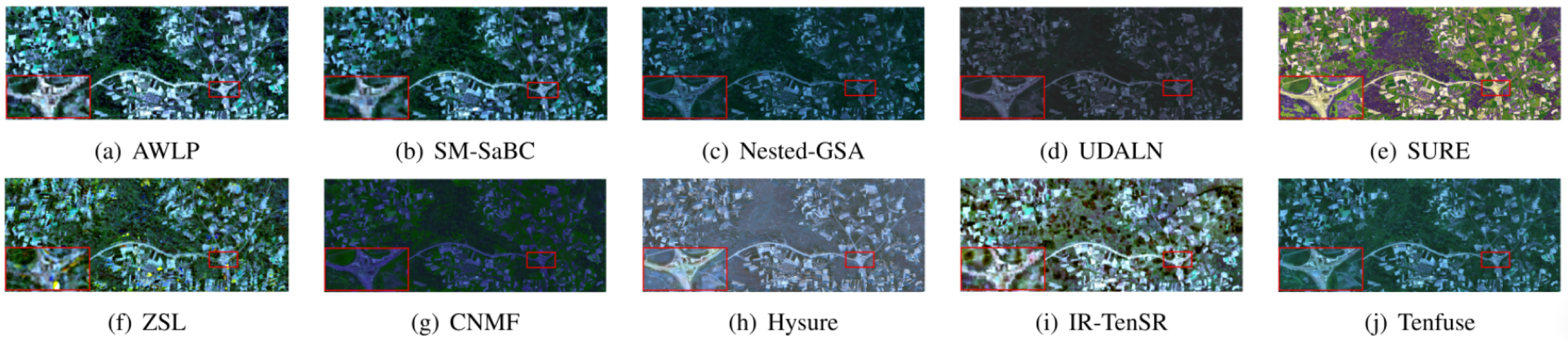

**Figure 6.8.** *Fusion results of real EnMAP/Sentinel-2 data. Pseudo-color visualization of reconstructed HR-HSIs using band composition (R:15/B:48/G:30) corresponding to wavelengths (487, 653, 561 nm).*

Regarding the spectral mixing window size $k_i$, we evaluated the model's sensitivity to the window size scale factor $c$, defined as the multiplier of the baseline window dimensions in the real-world Hyperion/Sentinel-2 experiment. As shown in Figure 6.11, the spectral distortion $D_\lambda$ reaches its minimum at $c = 1.5$, whereas larger scale factors lead to a significant increase in $D_\lambda$ due to excessive spectral mixing. Simultaneously, the spatial distortion $D_s$ rises steadily with $c$, with a sharp spike at $c = 10$ indicating a loss of fine spatial structures. Consequently, $c = 1.5$ is identified as the optimal configuration for achieving a superior balance between spectral fidelity and spatial detail preservation.

**6.6.2. The Moreau envelope parameter $\mu$.** The sensitivity of the model to the parameter $\mu$ was evaluated using the Washington DC Mall synthetic dataset. As illustrated in Figure 6.12, we plotted the PSNR performance against various values of $\mu$ to identify the

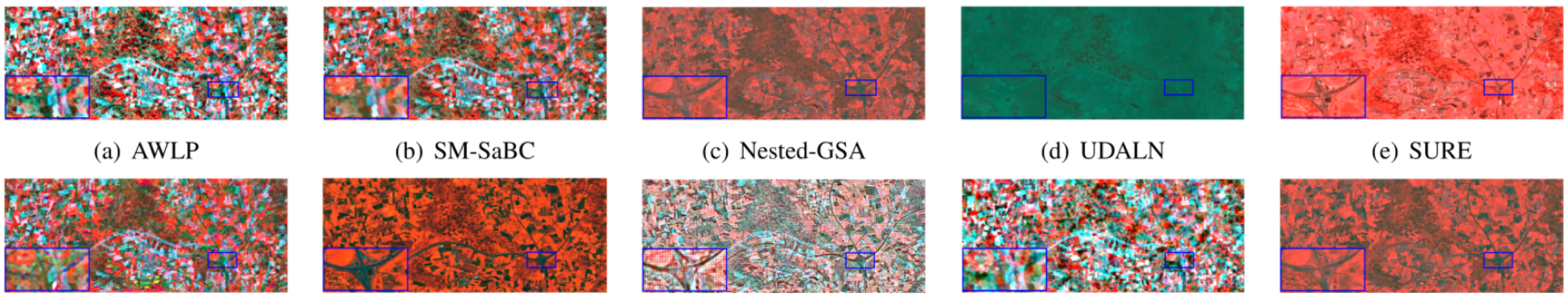

(a) AWLP (b) SM-SaBC (c) Nested-GSA (d) UDALN (e) SURE
(f) ZSL (g) CNMF (h) Hysure (i) IR-TenSR (j) Tenfuse

**Figure 6.9.** *Fusion results of real EnMAP/Sentinel-2 data. Pseudo-color visualization of reconstructed HR-HSIs using band composition (R:72/B:192/G:150) corresponding to wavelengths (823, 2204, 1619 nm).*

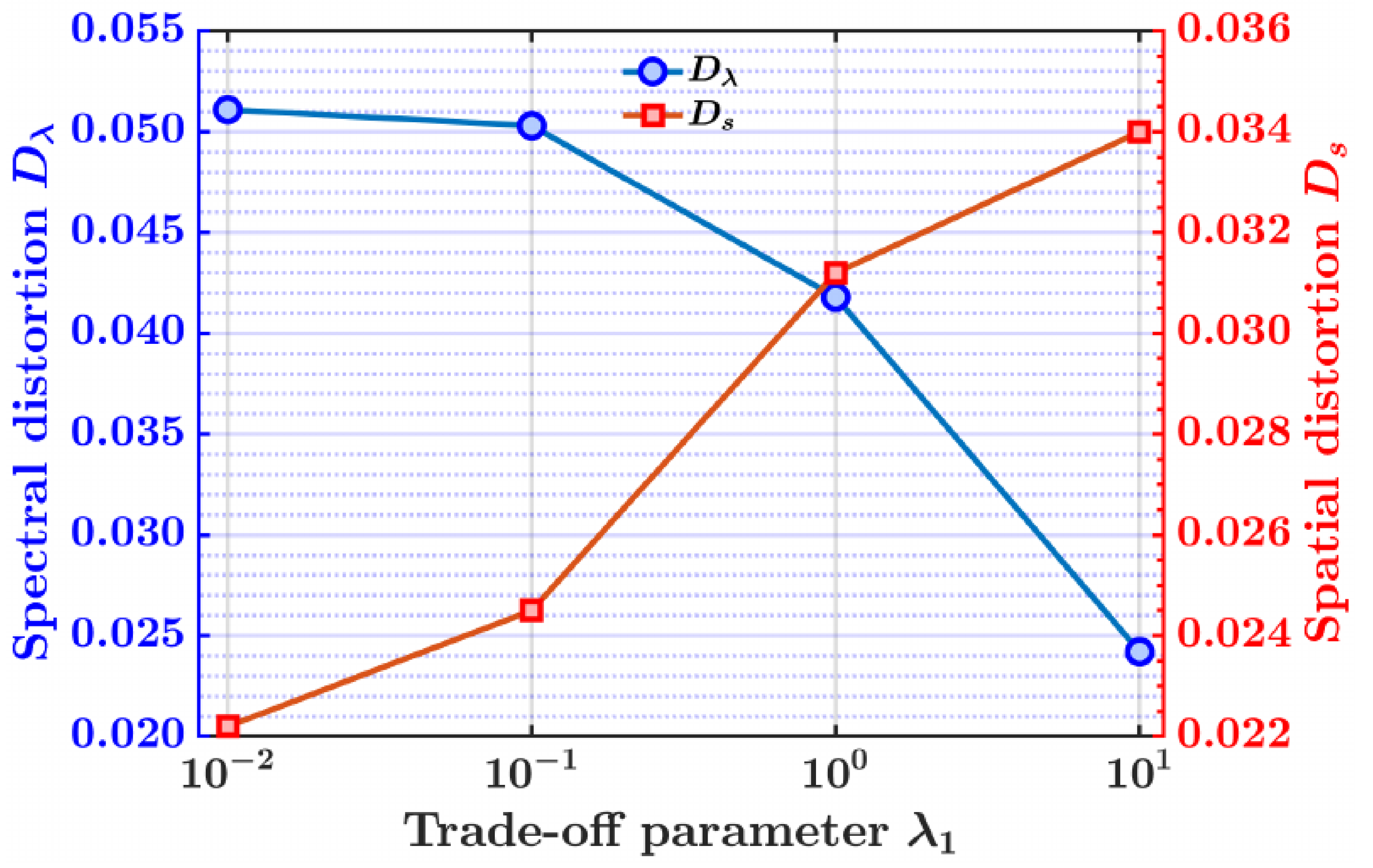

**Figure 6.10.** *Impact of regularization parameter $\lambda_1$.*

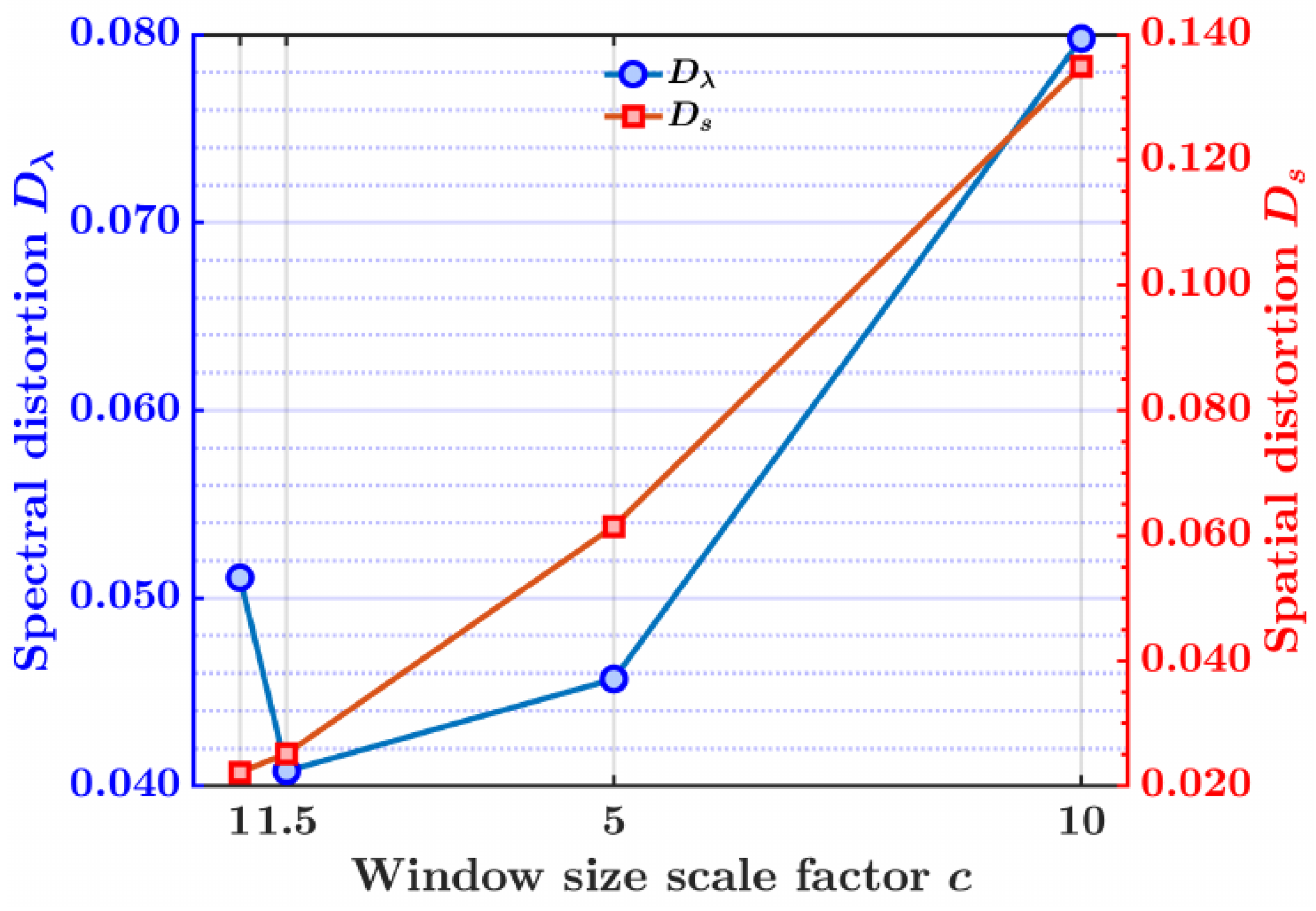

**Figure 6.11.** *Impact of spectral mixing window size $k_i$.*

optimal range for stable reconstruction. In this experiment, the parameter $\beta$ takes values from $\{1, 1, 10, 100, 1000\}$, associated with the values $\{20, 2, 0.2, 0.02, 0.002\}$ of $\mu$, respectively. The reported PSNR values represent the results obtained after 20 iterations of the algorithm.

**6.6.3. Robustness to mismatched, band-dependent noise.** To simulate realistic sensing environments, mismatched Gaussian noise was introduced into the synthetic Washington DC Mall dataset. Specifically, we assigned band-dependent standard deviations, $\sigma_i$, and applied discrepant noise intensities to the HSI and MSI observations. The MSI, characterized by higher spatial resolution, was configured with a band-average SNR of approximately 40 dB, whereas the lower-resolution HSI was set to 30 dB. As evidenced by the results in Table 6.4 and Figure 6.13, Tenfuse demonstrates remarkable robustness. Despite these non-ideal conditions, it consistently outperforms other methods across the majority of evaluation metrics.

**Table 6.4**

*Performance of ten test methods on synthetic Washing DC Mall with mismatched, band-dependent noise.*

| Metrics | **Hypersharpening-based methods** | | | **Learning-based methods** | | | **Model-based methods** | | | |
|---|---|---|---|---|---|---|---|---|---|---|
| | AWLP | SM-SaBC | Nested-GSA | UDALN | SURE | ZSL | CNMF | Hysure | IR-TenSR | Tenfuse |
| PSNR ($\infty$) | 35.61 | 38.80 | 47.25 | 49.29 | 47.89 | 49.99 | 11.38 | 47.36 | 31.86 | **50.22** |
| SAM (0) | 0.39 | 0.15 | 0.05 | 0.05 | 0.08 | 0.04 | 0.21 | 0.08 | 0.21 | **0.04** |
| UIQI (1) | 0.26 | 0.32 | 0.79 | 0.90 | 0.81 | 0.91 | 0.03 | 0.90 | 0.57 | **0.94** |
| ERGAS (0) | 3.81 | 2.28 | 1.78 | 1.21 | 1.75 | 1.10 | 68.47 | 1.14 | 5.67 | **1.08** |
| Time (s) | **0.70** | 1.81 | 1.02 | 3087.57 | 5337.10 | 4012.03 | 5.13 | 83.18 | 146.59 | 74.65 |

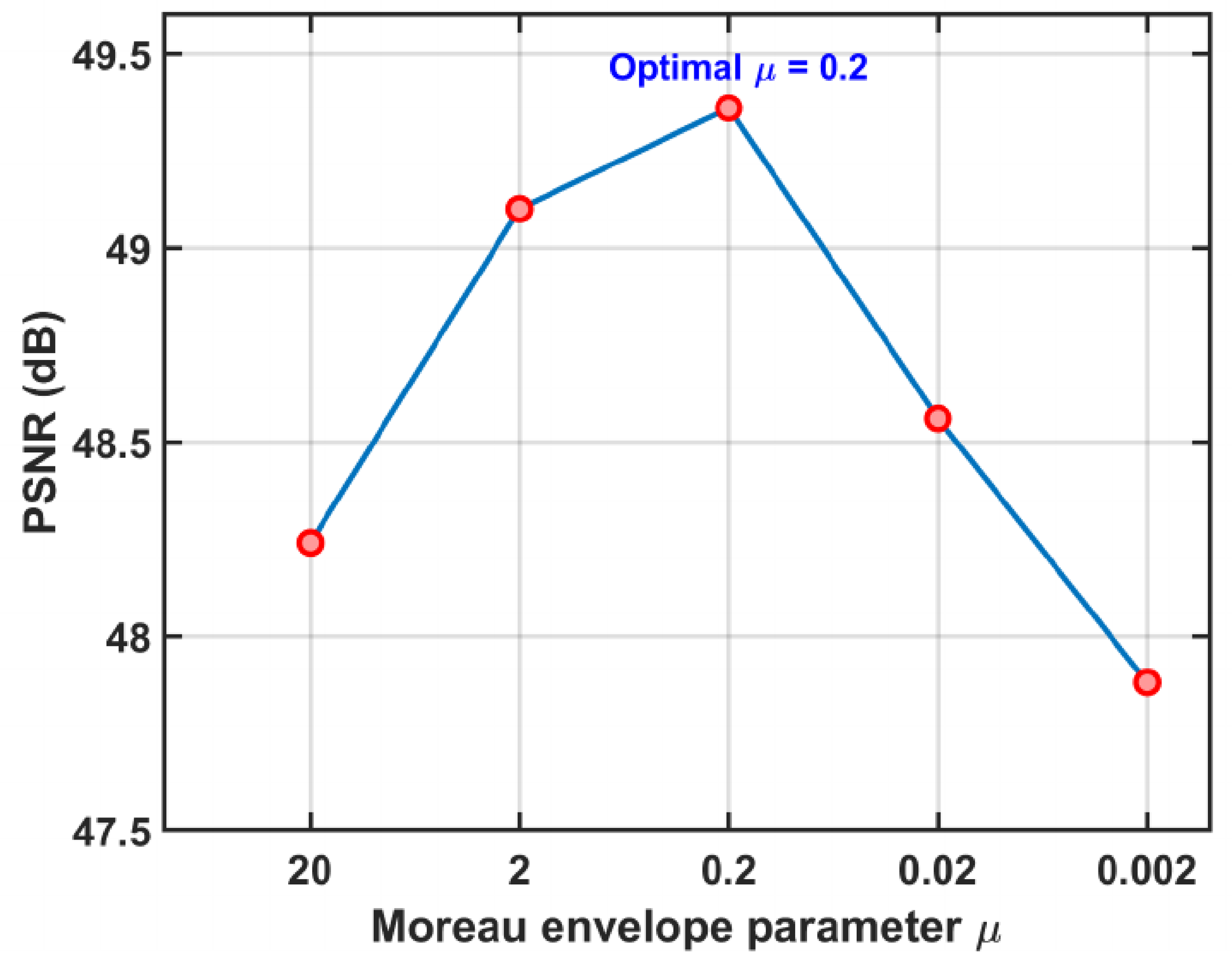

**Figure 6.12.** *Impact of the Moreau envelope parameter $\mu$.*

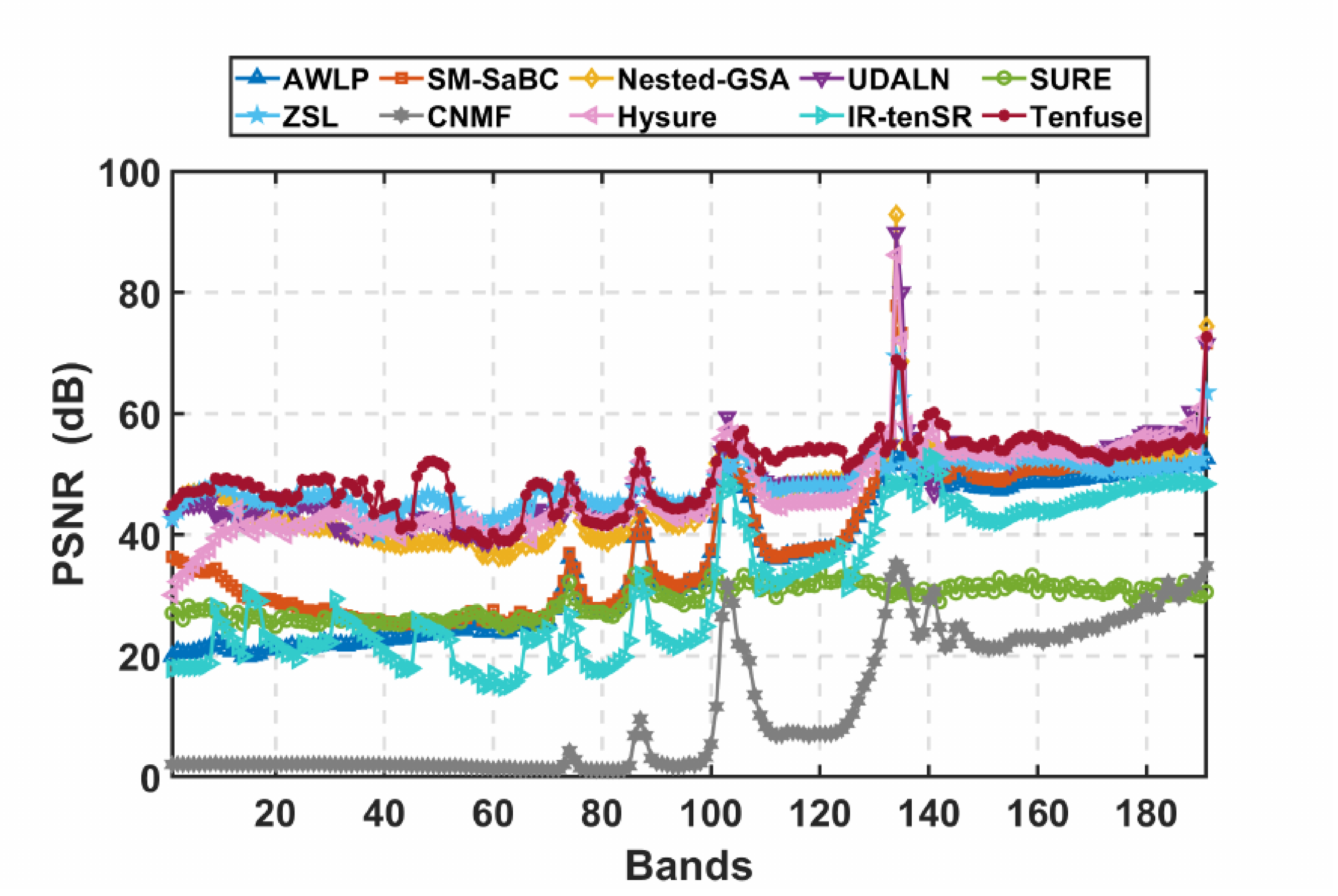

**Figure 6.13.** *Sensitivity to band-varying noise intensity.*

**7. Conclusion.** In this paper, we have formulated the blind HSI-MSI fusion problem as a coupled inverse problem, integrating blind deconvolution in the spatial domain with blind unmixing in the spectral domain. From this novel perspective, we proposed a unified tensor framework that leverages sensor information to achieve flexible self-adjustment and real-time data fusion. Distinctively, our approach is unsupervised and requires no pre-training, making it highly adaptable to real-world scenarios. We further introduced a joint optimization model to simultaneously estimate the target HR-HSI, the PSF, and the SRF. To solve this, a partially linearized ADMM with Moreau envelope smoothing was derived, supported by a rigorous theoretical convergence analysis. Furthermore, we developed an initialization estimator tailored to the specific physical characteristics of the fusion problem. Generally, our framework operates through the integration of sensor-information, physical degradation model, and customized algorithm, which are essential for achieving high-fidelity reconstruction of HR-HSI. Numerical comparisons with state-of-the-art methods on both synthetic and real-world datasets demonstrate the compelling performance of our proposed method.